\documentclass[a4paper]{article}
\usepackage{amsmath,amsthm,amssymb,mathrsfs}
\usepackage{enumerate}
\usepackage{chngcntr}
\usepackage{hyperref}
\usepackage{wrapfig}

\usepackage{graphicx}
\usepackage[top=2.5cm,bottom=2.5cm,left=2.5cm,right=2.5cm]{geometry}
\usepackage{color}
\usepackage{ulem}
\numberwithin{equation}{section}

\newtheorem{dfn}{Definition}[section]
\newtheorem{thm}[dfn]{Theorem}

\newtheorem{lem}[dfn]{Lemma}

\newtheorem{cor}[dfn]{Corollary}

\newtheorem{rem}[dfn]{Remark}
\newtheorem{prop}[dfn]{Proposition}

\newcommand{\LE}{{\rm LE}}
\newcommand{\fraka}{\mathfrak{a}}
\newcommand{\frakb}{\mathfrak{b}}
\newcommand{\frakc}{\mathfrak{c}}
\newcommand{\frakd}{\mathfrak{d}}
\newcommand{\Es}{{\rm Es}}
\newcommand{\mP}{{\mathbb{P}}}

\begin{document}
\title{\bf {\large  Convergence of three-dimensional loop-erased random walk in the natural parametrization}}
\author{Xinyi Li \thanks{Beijing International Center for Mathematical Research, Peking University. Email: \href{mailto:xinyili@bicmr.pku.edu.cn}{xinyili@bicmr.pku.edu.cn}} \and  Daisuke Shiraishi \thanks{Department of Advanced Mathematical Sciences, Graduate School of Informatics, Kyoto University. \newline Email: \href{mailto:shiraishi@acs.i.kyoto-u.ac.jp}{shiraishi@acs.i.kyoto-u.ac.jp}}}
\date{\today}

\maketitle
\begin{abstract}
In this work we consider loop-erased random walk (LERW) and its scaling limit in three dimensions, and prove that 3D LERW parametrized by renormalized length converges to its scaling limit parametrized by some suitable measure with respect to the uniform convergence topology in the lattice size scaling limit.  Our result greatly improves the work \cite{Koz} of Gady Kozma which establishes the weak convergence of the rescaled trace of 3D LERW towards a random compact set with respect to the Hausdorff distance. 
\end{abstract}
\tableofcontents

\section{Introduction}
\subsection{Introduction and main results}\label{sec:1.1}
Understanding the relationship between lattice-based probabilistic models and their limiting processes in $\mathbb{R}^{d}$ is one of the fundamental problems in probability theory. The archetypal (and quite well studied) example is, of course, Donsker's invariance principle which states that the scaling limit of simple random walk (SRW) is Brownian motion. However, many models arising from statistical physics are essentially more complicated than SRW, due to strong interaction of the process with its past. For instance, for random curves with strong self-repulsion (in other words, the curve is not allowed to visit its past),  even the existence of the scaling limit is not trivial at all. 

For models with strong self-repulsion, the most interesting cases are $d=2, 3$. In one dimension, the conditioning on self-repulsion forces the curve to be a straight line and nothing interesting happens. In four or more dimensions, loosely speaking, conditioning on self-avoidance does not have a strong effect on the behavior of the curve. Although this does not follow directly from the fact that Brownian motion is a simple curve in $d \ge 4$, 
nevertheless, in high dimensions, as a general strategy one can analyze rigorously curves with self-repulsion by comparing it with Brownian motion.   

In this article, we will focus on a particular example, the loop-erased random walk (LERW), which is the random  simple  path  obtained  through  erasing  all  loops  chronologically from an SRW path. We refer readers to Section \ref{LERW-intro} for an introduction to the model and a repository of tools needed in this work.

\medskip

LERW was originally introduced by Greg Lawler in \cite{Lawler1980}. Since then, it has been studied extensively both in mathematics and physics literature. Indeed, LERW has a strong connection with other models in statistical physics, especially the uniform spanning tree (UST) which arises in statistical physics in conjunction with the Potts model. Let us also mention that one can interpret LERW not only as the loop-erasure of the SRW but also as a Laplacian random walk (see Section \ref{LERW-intro} for details). Loosely speaking, conditioning the LERW $\gamma$ up to $k$-th step, the transition probability for the next step is given by the solution of a discrete Dirichlet problem on the complement of $\gamma [0, k]$. Although these models are also interesting on a general graph, in this article we shall restrict our attention to LERW in $\mathbb{Z}^{d}$ (and mainly with $d=3$).

As we have already pointed out, LERW for $d \ge 4$ was already well studied - it can be proved without resorting to lace expansion that in this case the scaling limit of LERW is Brownian motion (see Chapter 7 of  \cite{Lawb}). 

The case of $d=2$ is more difficult, yet also already well understood. Schramm first conjectured in \cite{Sch} that LERW has a conformally invariant scaling limit, which is characterized by Schramm-Loewner evolution (SLE) with parameter $2$. This conjecture was sub-sequently confirmed in \cite{LSW}. Since then, there have been substantial progresses on this subject, and thanks to our knowledge on SLE, one can prove very fine results on LERW in two dimensions, in particular convergence in the natural parametrization, see \cite{LV}.

What about $d=3$? Unfortunately, relatively little is known compared with other dimensions. The crucial reason is that we do not have nice tools like SLE to describe the scaling limit of 3D LERW directly. Hence, understanding this limit will most likely remain a hard task for a long time. 
That said, as the study of planar LERW first motivated the discovery of SLE, which fueled the breakthroughs in studies of other planar lattices models from statistical mechanics, understanding the LERW scaling limit and seeking its comprehensive description is, we believe, also vital to the study of other models from statistical physics in three dimensions.

\bigskip

We now turn to the main results in this paper. Let us start by briefly explaining notation as well as some known results on the scaling limit for 3D LERW here (see Section 2 for any missing definition). 

For $n>0$, ({\it not necessarily an integer}), let $S=S^{(n)}$ be a simple random walk on the rescaled lattice $2^{-n} \mathbb{Z}^{3}$ started from the origin\footnote{We use (dyadic) exponential scales throughout this work for the convenience of arguments involving dyadic partitions.}. We write $T=T^{(n)}$ for the first time that $S^{(n)}$ exits from $\mathbb{D} := \{ x \in \mathbb{R}^{3} : |x| <1 \}$ the unit open ball centered at the origin. We then write $\gamma_{n}: = \LE \big( S^{(n)} [0, T] \big)$ for the loop-erasure of $S^{(n)}$ up to the first exit time $T$. 

Let $\big( {\cal H} (\overline{\mathbb{D}} ), d_{\text{Haus}} \big)$ be the space of all non-empty compact subsets in $\overline{\mathbb{D}}$ endowed with the Hausdorff metric $d_{\text{Haus}}$. Regarding $(\gamma_{n})_{n\geq 0}$ as random elements of the metric space $\big( {\cal H} (\overline{\mathbb{D}} ), d_{\text{Haus}} \big)$, Gady Kozma proves in \cite{Koz} that for {\it integer}\footnote{This corresponds to dyadic scales for the underlying lattice.} $n$ there exists a random compact set ${\cal K}$ that $(\gamma_{n})_{n\in\mathbb{Z}^+}$ converges weakly as integer $n \to \infty$ with respect to the Hausdorff metric $d_{\text{Haus}}$ (in fact this convergence can be extended to all positive real $n$; see Section \ref{SCALING} and in particular Corollary \ref{cor:wCONV} for discussions on this). More precisely, he shows that $(\gamma_{n})_{n\geq 0}$ is a Cauchy sequence with respect to the Prokhorov metric. It is also proved that with probability one, ${\cal K}$ is a simple curve (see \cite{SS}) and the Hausdorff dimension of ${\cal K}$ is equal to a deterministic exponent $\beta \in (1, \frac{5}{3} ]$ (see \cite{S2}) where $\beta$ is some deterministic constant. See Section \ref{sec:Es} for more discussion on the exponent $\beta$.

In short, the purpose of this article is to study the scaling limit with respect to a topology stronger than the Hausdorff metric. Our choice of the topology is the supremum distance $\rho$ defined as follows: set $\big( {\cal C} (\overline{\mathbb{D}} ), \rho \big)$ for the space of continuous curves $\lambda : [0, t_{\lambda}] \to \overline{\mathbb{D}} $ with finite time duration $t_{\lambda} \in [0, \infty)$, and for $\lambda_1,\lambda_2\in {\cal C}$, define their distance under $\rho$-metric as 
\begin{equation}\label{eq:rho1st}
\rho (\lambda_{1}, \lambda_{2} ) = |t_{1} - t_{2} | + \max_{0 \le s \le 1} \big| \lambda_{1} ( s t_{\lambda_1} ) - \lambda_{2} ( s t_{\lambda_2} ) \big|. 
\end{equation}
To deal with the scaling limit with respect to the metric $\rho$, we consider the time rescaled LERW $\eta_{n}$ defined by
\begin{equation}\label{etan-1st}
\eta_{n} (t) = \gamma_{n} \big( f_n t\big) \ \text{ for } 0 \le t \le M_{n}/ f_n,
\end{equation}
where $M_{n}$ stands for the length (the number of lattice steps) of $\gamma_{n}$ (we let the walk traverse each edge in unit time and assume linear interpolation of $\gamma_{n}$ here so that $\eta_{n}$ becomes a continuous curve) and 
\begin{equation}\label{eq:fnpreview}
f_n \asymp 2^{\beta n} \mbox{ for all real }n>0    
\end{equation}
stands for a certain scaling factor to be fixed later. See \eqref{eq:fndef} for its precise definition. We also remark that if we restrict ourselves to {\it integer} $n$'s, we can safely pick the explicit scaling factor $f_n=2^{\beta n}$. Here we mention that a recent work \cite{HTLS} by Hern\'{a}ndez-Torres and the authors of this work allows us to pick such explicit scaling factor for {\it all} real $n>0$.  See Remark \ref{rem:mainrem} 4) and below \eqref{eq:fndef} for more discussions. 

When we study $\eta_{n}$ and its limiting process with respect to the distance $\rho$, the first crucial issue is to give a suitable time parametrization for the scaling limit ${\cal K}$. With this in mind, we begin with the definition of the following {\bf occupation measure} $\mu_{n}$ of $\gamma_n$:
\begin{equation*}
\mu_{n} := \big(f_n\big)^{-1} \sum_{x \in \gamma_{n} \cap 2^{-n} \mathbb{Z}^{3} } \delta_{x}, 
\end{equation*}
where $\delta_{x}$ is the Dirac measure at $x$. Note that for each point $x$ lying on the curve $\eta_{n}$, we can compute the exact time that $\eta_{n}$ passes through $x$ by measuring the weight of the sub-path of $\gamma_{n}$ between the origin and $x$ via the measure $\mu_{n}$. In this way, the curve $\eta_{n}$ is obtained by parametrizing $\gamma_{n}$ by the measure  $\mu_{n} $. Thus, it is natural to consider the limiting measure of $\mu_{n}$ and parametrize ${\cal K}$ by this measure. As in the discrete case where we give equal weight to each lattice point hit by $\gamma_n$ (or equivalently, traverse $\eta_n$ in a constant speed), we call this kind of parametrization the {\bf natural parametrization} of $\cal K$.

The first main result guarantees the existence of the joint limit of $\gamma_n$ and the occupation measure $\mu_{n}$. Let ${\cal M} (\overline{\mathbb{D}})$ be the space of all finite measures on $\overline{\mathbb{D}}$ endowed with the weak convergence topology.

\begin{thm}\label{1st}
For $n\in \mathbb{R}^+$, as $ n \to \infty$, the sequence of the joint law $(\gamma_{n}, \mu_{n} )$ converges weakly to some $({\cal K}, \mu)$ with respect to the product topology of ${\cal H} (\overline{\mathbb{D}})$ and ${\cal M} (\overline{\mathbb{D}})$. 
Furthermore, the limit measure $\mu$ is a measurable function of ${\cal K}$.
\end{thm}

In order to parametrize ${\cal K}$ via the measure $\mu$, we need the following basic properties of $\mu$. For a point $x \in {\cal K}$, let ${\cal K}_{x}$ be the simple curve on ${\cal K}$ between the origin and $x$ (recall that ${\cal K}$ is a simple curve almost surely, see Section \ref{SCALING}).

\begin{thm}\label{2nd}
With probability one, the support of the measure $\mu$ coincides with ${\cal K}$. Moreover, it follows that with probability one, for each $x \in {\cal K}$
\begin{equation*}
\lim_{\substack{y \in {\cal K} \\ y \to x}} \mu ( {\cal K}_{y} ) = \mu  ( {\cal K}_{x} ).
\end{equation*}

\end{thm}

We now parametrize ${\cal K}$ through the measure $\mu$. By Theorem \ref{2nd}, it follows that for each $t \in [0, \mu ( {\cal K} )]$, there exists a unique point $x_{t} \in {\cal K}$ satisfying $t = \mu ( {\cal K}_{x_{t}} )$. Define $\eta (t) = x_{t}$ for $t \in [0, \mu ( {\cal K} )]$. It also follows from Theorem \ref{2nd} that $\eta$ is a random continuous curve in ${\cal C}(\overline{\mathbb{D}})$ whose time duration is $\mu ( {\cal K} )$. The next theorem gives the desired convergence in the space $({\cal C}(\overline{\mathbb{D}}),\rho)$.

\begin{thm}\label{3rd}
For $n\in \mathbb{R}^+$, as $n \to \infty$, $\eta_{n}$ converges weakly to $\eta$ in the space $({\cal C}(\overline{\mathbb{D}}),\rho)$.
\end{thm}

As a corollary of these theorems, we can also deal with the scaling limit of the infinite loop-erased random walk (ILERW) as follows. Recall that $S^{(n)}$ stands for SRW on $2^{-n} \mathbb{Z}^{3}$ started at the origin. Since $S^{(n)}$ is transient,  the loop-erasure of $S^{(n)} [0, \infty )$ is well-defined. We then write $\gamma_{n}^{\infty} = \text{LE} \big( S^{(n)} [0, \infty ) \big) $ for the ILERW on $2^{-n} \mathbb{Z}^{3}$. We also consider the time rescaled version defined by 
\begin{equation*}
\eta_{n}^{\infty} (t) = \gamma_{n}^{\infty} (f_n t ) \ \text{ for } t \ge 0,
\end{equation*}
where again we assume the linear interpolation of $\gamma_{n}^{\infty}$ so that $\eta_{n}^{\infty}$ becomes a random element of $\big( {\cal C}, \chi \big) $ the metric space of continuous curves lying in $\mathbb{R}^{3}$ defined on $[0, \infty )$ (see Section \ref{metric} for the space ${\cal C}$), equipped with the metric $\chi$ 
\begin{equation}\label{eq:chi1st}
\chi (\lambda_{1}, \lambda_{2} ) = \sum_{k=1}^{\infty} 2^{-k} \max_{0 \le t \le k} \min \Big\{ \big| \lambda_{1} ( t ) - \lambda_{2} ( t) \big|, 1 \Big\},
\end{equation}
for two continuous curves $\lambda_{1}, \lambda_{2} \in {\cal C}$.  
The next theorem confirms the existence of the scaling limit of ILERW with respect to the metric $\chi$.

\begin{thm}\label{4th}
There exists a random continuous curve $\eta^{\infty}\in {\cal C}$ such that as $n \to \infty$, $\eta_{n}^{\infty} $ converges weakly to $\eta^{\infty}$ with respect to the metric $\chi$.
\end{thm}

We now give some comments on our results and briefly discuss some related open questions. See also Remark \ref{rem:lastrem} for comments on the scaling limit of ILERW.

\begin{rem}\label{rem:mainrem}

{\rm 1)} To our best knowledge, LERW is the only natural model with self-repulsion for which we can prove the existence of the scaling limit in natural parametrization in every dimension. See {\rm \cite[Chapter 7]{Lawb}} for $d \ge 4$ and {\rm \cite{LVradial,LV}} for $d=2$. However, the strategies are quite different for different dimensions. In the case of $d\geq 4$, the proof crucially relies on the fact that the scaling limit is a Brownian motion, while in the case of $d=2$, it is the conformal invariance of critical planar models and the analysis of the Schramm-Loewner evolution that play a key role. It is worth noting that neither proof easily generalizes to the 3D case, while our proof strategy can be adapted to other dimensions.

\medskip

\noindent {\rm 2)} So far we still do not have a nice description of the scaling limit $\eta$ or $\eta^\infty$, despite the convergence results we obtain in this work. In {\rm \cite{SS}}, it is conjectured that $\eta$ is the unique random continuous simple curve that can generate a Brownian motion if loops from Brownian loop soup (see {\rm \cite[Conjecture 1.3]{SS}} for the precise form) are added appropriately. Unfortunately, the present paper does not give any progress for this conjecture. Hence it remains a big challenge to give a ``good" description for $\eta$.

\medskip

\noindent {\rm 3)} As the setup of ILERW is in some sense more natural than that of LERW in a domain, we hope that showing the existence of $\eta^\infty$ would facilitate studies of 3D LERW and related subjects. For instance, results established in this work play a crucial role in the study of 3D uniform spanning trees in the paper {\rm \cite{ACHS}}.  

\medskip

\noindent {\rm 4)} In the case of two dimensions, the scaling limit ${\rm SLE}_2$ is parametrized by its {\bf Minkowski content} in the natural parametrization (see {\rm \cite{LV}}), it is very natural to conjecture that in three dimensions, it is also possible to identify $\mu$ with the Minkowski content of $\eta$. Although in this work we can show that $\mu$ is measurable with respect to $\cal K$ and one can regard $\mu$ as a measure generated from the ``box content'' of $\cal K$ in some sense, it remains a challenge to show that the Minkowski content exists for $\cal K$ and is indeed a deterministic multiple of $\mu$. Recently, a preprint of Hern\'{a}ndez-Torres and the authors on this work answers this question positively; see \cite{HTLS} and the discussions therein. 

\medskip



\noindent {\rm 5)} In {\rm \cite[Theorem 7]{Koz}}, Kozma established a version for the ``universality'' of the scaling limit of LERW under the Hausdorff topology, for graphs that can be ``isotropically interpolated'' to $\mathbb{Z}^3$. It is then a natural question to wonder if this can be extended to convergence in natural parametrization and if so whether it is necessary to pose extra assumptions.

\end{rem}

\subsection{Some words about the proof}\label{sec:1.2}
Theorems \ref{1st} and \ref{3rd}  are the meat of the present article, while Theorems \ref{2nd} and \ref{4th} follow from routine arguments and are neither surprising nor difficult. For Theorem \ref{4th}, let us mention that the distribution of an initial part of $\gamma_{n}$ is similar to that of $\gamma_{n}^{\infty}$ (see \eqref{eq:LERWvsILERW} for a precise statement). Therefore, once we establish Theorem \ref{3rd}, we can construct the scaling limit of ILERW by ``attaching" initial parts of $\eta$ (or some version of $\eta$) appropriately and thus obtain Theorem \ref{4th}. See Section \ref{sec:8} for the precise argument. 

Concerning the proof of Theorem \ref{3rd}, we need to show the tightness of $\eta_{n}$ with respect to the metric $\rho$, then prove the uniqueness of the limit. In a sub-sequent work of the authors \cite{Holder}, the argument for tightness has been improved to obtain the optimal H\"older exponent for the scaling limit. To keep this work as short as possible, we omit the proof for tightness and refer readers to Remark \ref{rem:nocircular} for relevant discussions. As the proof for uniqueness is a bit technical, we postpone the sketch of the proof to Section \ref{sec:guidepara}. So, we will only sketch the proof of Theorems \ref{1st} here.


It follows from \cite[Corollary 1.3]{Escape}  that the sequence of the joint laws $(\gamma_{n}, \mu_{n} )$ is tight with respect to the product topology of  $ {\cal H} (\overline{\mathbb{D}} )$ and  ${\cal M}  (\overline{\mathbb{D}} )$; see also Proposition \ref{tight} of the present article. Thus, we can find a sub-sequence ${n_{k}}$ and some measure $\mu$ such that  $(\gamma_{n_{k}}, \mu_{n_{k}})$ converges weakly to $({\cal K}, \mu)$. (In the sequel we write $n = n_{k}$ for simplicity.) We want to show that $\mu = f ( {\cal K} )$ where $f$ is a deterministic function which does not depend on the chosen sub-sequence, which implies Theorem \ref{1st}. To this end, take a cube $B \subset \mathbb{D}$ with ${\rm dist} (0 \cup \partial \mathbb{D}, B ) > 0$. In order to prove the claim above, the crucial step is to show that $\mu (B)$ is measurable with respect to ${\cal K}$. We will explain how to prove the measurability of $\mu (B)$ in the rest of this subsection, where we assume that $\big\{ (\gamma_{n}, \mu_{n} ) \big\}_{n\geq 0}$ and $({\cal K}, \mu)$ are defined in the same probability space such that $(\gamma_{n}, \mu_{n} )$ converges to $({\cal K}, \mu) $  almost surely.

 Take $\epsilon > 0$ with $2^{-n} \ll \epsilon$. Readers should regard $\epsilon$ as a mesoscopic scale quantity\footnote{In Sections \ref{sec:3} and \ref{sec:4}, we pick $\epsilon=2^{-k^4}$.}, while the mesh size $2^{-n}$ is in the microscopic scale. Decompose $B$ into a collection of cubes $B_{1}, \cdots B_{N_{\epsilon}}$ of side length $\epsilon$.  Let $X_{i}$ be the number of points in $B_{i} \cap 2^{-n} \mathbb{Z}^{3}$ hit by $\gamma_{n}$, and let $Y_{i} = {\bf 1} \{ \gamma_{n} \cap B_i \neq \emptyset \}$ be the indicator function of the event that $\gamma_{n}$ hits $B_i$. By definition, $\mu_{n} (B) \approx 2^{-\beta n} \sum_{i} X_{i}$ (note that some vertices on the boundaries of the cubes may be double-counted or left uncounted).

It will be shown in Corollary \ref{COR} below that we can regard the indicator function $Y_{i}$ as an almost ``measurable" function of  ${\cal K}$. Namely, one has $Y_{i} = Z_{i} := {\bf 1} \{ {\cal K} \cap B_{i} \neq \emptyset \}$ with high probability. With this in mind, we aim to show that the sum of $X_{i}$ is well approximated by that of $\alpha_{0} Y_{i}$ (or equivalently $\alpha_{0}Z_{i}$) with some deterministic $\alpha_{0}$ (see \eqref{Beta} for its precise definition and \eqref{proof-sketch-1} for its connection with $\cal K$) which does not depend on the choice of $i$. 
In other words,
$\alpha_{0} \approx \alpha_{i} := E \big( X_{i} \ \big| \ \gamma_{n} \cap B_{i} \neq \emptyset  \big)$. This approach is justified by proving the following key $L^{2}$-estimate:
\begin{equation}\label{KEY}
E \Big[ \Big( \sum_{i} X_{i} -  \sum_{i} \alpha_{0} Y_{i} \Big)^{2} \Big] \le \xi (n, \epsilon ) \Big( E \Big[ \sum_{i} X_{i} \Big] \Big)^{2},
\end{equation}
where $\xi (n, \epsilon)$ converges to zero when $2^{-n}$ tends to zero first and then $\epsilon$ tends to zero sub-sequently. 
The reader should regard the inequality \eqref{KEY} as a law of large numbers type of estimate, despite that $X_{i}$ is not i.i.d.\ for our model.
The estimate \eqref{KEY} bears some similarity to \cite[Theorem 1.1]{GPS} (which treats various occupation measures for critical planar percolation) in terms of structure.  However, due to the lack of knowledge for 3D LERW and lack of full spatial independence, establishing this estimate is quite difficult compared to the case of critical planar percolation in terms of techniques. The actual proof of \eqref{KEY} is unfortunately quite long and involved, and the whole Section \ref{sec:4} is dedicated to it.  We refer curious readers to the end of Section \ref{sec:sketch} for a sketch of its proof.

To finish the proof of Theorem \ref{1st}, we set $B_0$ for the cube of side length $\epsilon$ centered at $(1/2, 0, 0)$.  We will then prove in Proposition \ref{prop:alpha0}  below (see \eqref{2022-10-8} for the precise estimate) that there exists a universal constant $c_{0}$ which especially does not depend on  the chosen sub-sequence $n = n_{k}$ such that 
\begin{equation}\label{proof-sketch-1}
\alpha_{0} (f_n)^{-1} = \big[ 1 + \widehat{\xi} (n, \epsilon) \big] \, \frac{ c_{0} \epsilon^{3}}{P ( {\cal K} \cap B_0 \neq \emptyset )},
\end{equation} 
where $\widehat{\xi} (n, \epsilon)$ tends to zero in the same way as $\xi(n,\epsilon)$. 
Combining  \eqref{KEY}, \eqref{proof-sketch-1} and the fact that $Y_{i} = Z_{i}$ with high probability, we conclude that 
$$\mu (B) \approx \mu_{n} (B) = (f_n)^{-1}  \sum_{i} X_{i} \simeq  (f_n)^{-1}  \sum_{i} \alpha_{0} Y_{i} \approx  \frac{ c_{0} \epsilon^{3}}{P ( {\cal K} \cap B_0 \neq \emptyset )} \, \sum_{i} Z_{i}$$
where the last quantity is measurable with respect to ${\cal K}$, as desired.

The organization of this article is as follows. We will introduce general notation and recall some basic facts   about simple random walk and metric spaces which will be used throughout this article in Section \ref{sec:2}. We review the definition and various key properties of loop-erased random walk in Sections \ref{LERW-intro} and \ref{sec:LERWlimit}.  A crucial part of this work is the second moment estimate  \eqref{KEY}  on the occupation measure and the number of cubes hit by the LERW. The precise setup and some preliminary estimates will be given in Section \ref{sec:3}. 
The key estimate \eqref{KEY} will be proved in Section \ref{sec:4} (see Proposition \ref{prop4.6} for precise statements). Using this $L^{2}$-estimate, for each cube $B \subset \mathbb{D}$, we will give an $L^{2}$-approximation of $\mu_{n} (B)$ by some measurable quantity with respect to ${\cal K}$ in Section \ref{sec:5}; see in particular Proposition \ref{prop3}. Theorem \ref{1st} will be proved in Section \ref{sec:6} (Theorem \ref{goal}). 
We will prove Theorem \ref{2nd} and Theorem \ref{3rd} in Section \ref{sec:7};  see Propositions \ref{basic2} and \ref{basic3} for Theorem \ref{2nd} and see Proposition \ref{end} for Theorem \ref{3rd}, respectively. Finally, we will study the ILERW in Section \ref{sec:8} and prove Theorem \ref{4th} in Theorem \ref{thm:4rep}. 

\bigskip

\noindent {\bf Acknowledgements:} The authors are grateful to Greg Lawler for numerous helpful and inspiring discussions and to anonymous referees for various useful comments and suggestions on previous versions of this work. XL wishes to thank Kyoto University for its warm hospitality during his visit, when part of this work was conceived. DS wishes to thank The University of Chicago for its warm hospitality during his visit, when part of this work was conceived. The authors would also thank the generous support from the National Science Foundation via the grant DMS-1806979 for the conference ``Random Conformal Geometry and Related Fields'', where part of this work was accomplished. XL is supported by National Key R\&D Program of China (No.\ 2020YFA0712900 and No.\ 2021YFA1002700) and NSFC (No.\ 12071012). DS is supported by a JSPS Grant-in-Aid for Early-Career Scientists, 18K13425 and JSPS KAKENHI Grant Number  17H02849, 18H01123, 21H00989, 22H01128 and  22K03336.

\section{Notation and background}\label{sec:2}
In this section, we will introduce notation and  background facts  that will be useful later.

\subsection{General notation}
We start with set-theoretical notation. For $d\geq 1$, we write $\mathbb{Z}^d$ and $\mathbb{R}^d$ for $d$-dimensional integer lattice and Euclidean space respectively.  In this article, unless otherwise mentioned, we will consider $d=3$ only. 

For a subset $A \subset \mathbb{R}^{d}$, we denote its boundary by $\partial A$. Let $\overline{A} = A \cup \partial A$ be the closure of $A$ and $A^\circ= A\setminus \partial A$ the interior of $A$. For two subsets $A, B \subset \mathbb{R}^{d}$, we write ${\rm dist} (A, B) = \inf_{x \in A, y \in B} |x-y|$ for the distance between $A$ and $B$ where $| \cdot |$ stands for the Euclidean distance in $\mathbb{R}^{d}$. If $A = \{ x \}$, we write ${\rm dist} (\{ x \}, B) = {\rm dist} (x, B)$. Furthermore, we let $A+B = \{ a + b \ | \ a \in A, \ b \in B \}$ be the Minkowski sum of $A$ and $B$. If $A = \{ x \}$ for some $x \in \mathbb{R}^{d}$, we write $x + B $ instead of $\{ x \} + B$. For $r > 0$, we set $r A = \{ r a \ | \ a \in A \}$.
We write $\mathbb{D} = \{ x \in \mathbb{R}^{d} \ | \ |x | < 1 \}$ for the unit open ball in $\mathbb{R}^{d}$ centered at the origin. Denote its closure by $\overline{\mathbb{D}}$. A set $B$ is called a (open) cube if it is a set  of the form $B = \{ y \in \mathbb{R}^{3} \ | \  \|  x-y   \|_{\infty} < r \}$ for some $x \in \mathbb{R}^{3}$ and $r > 0$, where $\|  \cdot \|_{\infty} $ stands for the $L^{\infty}$-norm. A set of the form $ \{ y \in \mathbb{R}^{3} \ | \  \|  x-y   \|_{\infty} \le r \}$  is called a closed cube. When we say $B$ is a cube, it means that $B$ is either open or closed. We also set $\mathbb{B} = \{ x \in \mathbb{R}^{d} \ | \ \| x \|_{\infty} < 1 \}$ for the open cube of side length $2$ in $\mathbb{R}^{d}$ centered at the origin where $\| \cdot \|_{\infty}$ stands for the $L^{\infty}$-norm. We write $\mathbb{D}_r:=r\mathbb{D}$ and $\mathbb{B}_r:=r\mathbb{B}$ for dilations of $\mathbb{D}$ and $\mathbb{B}$ resepectively. 

We now turn to discrete objects. For most of the time, we will focus on the  rescaled lattice $2^{-n} \mathbb{Z}^d$ with mesh size $2^{-n}$ for some $n \ge 0$. For a subset $A \subset  2^{-n}\mathbb{Z}^{d}$, we let $\partial A = \{ x \in  2^{-n}\mathbb{Z}^{d} \setminus A \ | \ \exists y \in A \text{ such that } |x-y| = 2^{-n}  \}$ be the outer boundary of $A$. We write $\partial_{i} A = \{ x \in  A \ | \ \exists y \in 2^{-n}\mathbb{Z}^{d} \setminus A \text{ such that } |x-y| = 2^{-n}  \}$ for the inner boundary of $A$. For $x \in 2^{-n}\mathbb{Z}^{d}$ and $r > 0$, we set $D (x, r) = \{ y \in 2^{-n}\mathbb{Z}^{d} \ | \ |x-y| < r \}$ for the (discrete) ball of radius $r$ centered at $x$. We write $D (r) = D (0, r)$ when $x=0$. 
We also write  $B_\infty(x, r) = \{ y \in 2^{-n}\mathbb{Z}^{d} \ | \|x-y\|_\infty < r \}$ for a cube of side-length $2r$ centered at $x$ and omit $x$ when $x=0$.

\medskip
We now turn to paths. 
As already introduced in the first section,  we will consider the rescaled lattice $2^{-n}\mathbb{Z}^{d}$ along with $\mathbb{Z}^{d}$ for the underlying graph. We will specify the scale of the lattice whenever there is a risk of confusion. 

A {\it path} $\lambda = [\lambda (0), \lambda (1), \cdots , \lambda (m) ] \subset 2^{-n}\mathbb{Z}^{d}$ is a sequence of points lying on $2^{-n}\mathbb{Z}^{d}$ satisfying $| \lambda (i-1)  - \lambda (i) | = 2^{-n}$ for each $1 \le i \le m$. We call $m$ the length of $\lambda$ and denote it by $\Diamond=\Diamond(\lambda)$. For $\lambda$ a path of length $m$, we write $\lambda^{\cal R}=[\lambda (m), \lambda (m-1), \cdots , \lambda (0) ]$ for the time reversal of $\lambda$. If $\lambda (i) \neq \lambda (j) $ for any $i \neq j$, we call $\lambda$ a simple path. For two paths $\lambda = [\lambda (0), \lambda (1), \cdots , \lambda (m) ] \subset 2^{-n}\mathbb{Z}^{d}$ and $\lambda' = [\lambda' (0), \lambda' (1), \cdots , \lambda' (m') ] \subset 2^{-n}\mathbb{Z}^{d}$ with $\lambda (m) = \lambda' (0)$, let $\lambda \oplus \lambda' = [\lambda (0), \lambda (1), \cdots , \lambda (m), \lambda' (1), \cdots , \lambda' (m') ]$ be their  concatenation.

For a path $\lambda$ started in $D(x,r)$, 
let $T_{x, r}= T^{(n)}_{x, r} = \min \{ k \ge 0 \ | \ \lambda(k) \notin D (x, r) \}$ stand for the first time that $\lambda$ exits from $D (x, r)$ and $T_{x,s,r} = \max \{ k \le T_{x,r} \ | \ \lambda(k) \in D (x, s) \}$. Here we use the convention that $\inf \emptyset = + \infty $. Also, we set $T_{r} = T_{r}^{(n)} = T_{0, r}^{(n)} $ for the case that $x = 0$ and set $T_{s,r}$ in the same manner. 


\medskip

We end this subsection with conventions on asymptotics and constants. For two sequences $a_{n}$ and $b_{n}$, we write 
\begin{itemize}
\item $a_{n} \asymp b_{n}$ if $\exists c > 0$ such that $c a_{n} \le b_{n} \le \frac{1}{c} a_{n}$ for all $n$;
\item $a_{n} \sim b_{n}$ if $\lim_{n \to \infty} a_{n}/b_{n} = 1$;
\item $a_{n} \approx b_{n}$ if $\log a_{n} \sim \log b_{n}$.
\end{itemize}
For two functions $f (x)$ and $g (x)$, we write $g(x) = O \big( f(x) \big)$ if $|g (x)| \le c f(x)$ for some universal constant $c > 0 $. If we wish to imply the constant $c$ may depend on another quantity, say $\epsilon $, we write $O_{\epsilon} \big( f(x) \big) $. Similarly, we write $g (x) = o \big( f(x) \big) $ if $g (x ) / f (x) \to 0$.
For $a \in \mathbb{R}$, we write $\lfloor a \rfloor$ for the largest integer less than or equal to $a$.

We use $c,\ C,\  c_{i},\ C_{i},\cdots$ ($i=0, 1,  \cdots $) to denote arbitrary universal positive constants which may change from line to line. If a constant depends on some other quantity, this will be made explicit. For example, if $c$ depends on $\epsilon$, we write $c_{\epsilon}$.

\subsection{Random walk and estimates of Green's function}\label{sec-green}

In this subsection, we will introduce some useful estimates on discrete harmonic functions and the simple random walk which will be used many times in this paper. Note that as for most of the time of the main text we will be working on rescaled lattices, the results introduced here will appear later in a scaled fashion.

We start with a review of the Harnack principle. Let $O \subset \mathbb{R}^{d}$ be a connected open set and $F $ be a compact subset of $O$. Then the Harnack principle states that there exist $C = C (F, O) < \infty$ and $N_0 = N_0 (F, O) < \infty$ such that if $N \ge N_0$,
\begin{equation*}
O_{N} = \{ x \in \mathbb{Z}^{d} \ | \ N^{-1} x \in O \},  \ \ \ \ \ \  F_{N} = \{ x \in \mathbb{Z}^{d} \ | \ N^{-1} x \in F \}
\end{equation*}
and $f : O_{N} \to [0, \infty )$ is discrete harmonic in $O_{N}$, then  
\begin{equation}\label{Harnack}
f (x) \le C f (y) 
\end{equation}
for $x, y \in F_{N}$. See \cite[Theorem 6.3.9]{LawLim} for further details. 

We now turn to random walks. We write $S$ for a  simple random walk (SRW) in $\mathbb{Z}^{d}$. 
We set $P^{x}$ and $E^{x}$ for its probability law and the expectation when it starts from a point $x$. We omit the superscript if $x = 0$.  

Let $1 \le M < N $. Write $A =A_{M,N} = \{ x \in \mathbb{Z}^{3} \ | \ M \le |x| \le N \}$. Set $\tau = \inf \{ k \ge 0 \ | \  S(k) \notin  A  \}$ for the first exit time of the discrete annulus. Then it is proved in \cite[Proposition 1.5.10]{Lawb} that for all $x \in A$, 
\begin{equation}\label{srwbound}
P^{x} \big( |S_{\tau} | \le M \big) = \frac{ |x|^{-1} - N^{-1}  + O (M^{-2} )  }{M^{-1} - N^{-1}}.
\end{equation}

We often want to consider the case that $M=1$ and $|x|$ is large. In that case, the estimate \eqref{srwbound} is not useful because the error term $O (M^{-2} )$ is much bigger than $|x|^{-1}$. To deal with this issue, we introduce the Green's function $G : \mathbb{Z}^{3} \times \mathbb{Z}^{3} \to [0, \infty )$ by 
\begin{equation}\label{green-entire}
G (x, y) = E^{x} \Big( \sum_{j =0}^{\infty} {\bf 1} \{ S (j) = y \} \Big)  \text{ for } x, y \in \mathbb{Z}^{3}.
\end{equation} 
Since $S$ is transient, $G (x, y) < \infty$ for each $x, y \in \mathbb{Z}^{3}$. We write $G (x)  = G (0, x)$.   \cite[Theorem 1.5.4]{Lawb} shows that there exists a universal constant $a > 0$ such that 
\begin{equation}\label{lawler-green}
G (x) = a |x|^{-1} + O \big( |x|^{-2} \big) \text{ as } |x| \to \infty.
\end{equation}
Suppose $M =1$. Note that $|S_{\tau} | \le M \Leftrightarrow S_{\tau} = 0$ in this case. By considering a bounded martingale $ G \big( S (j \wedge \tau) \big)$, we have 
\begin{equation}\label{srwbound-2}
P^{x} \big( S_{\tau} = 0 \big) = \frac{ a |x|^{-1} - a N^{-1}  + O (|x|^{-2} )  }{G (0) - a N^{-1}}.
\end{equation}

We now consider a domain $A \subset \mathbb{Z}^{d}$.  For $x,y \in \mathbb{Z}^{d}$, we write $G_{A} (x, y) $ for the Green's function on $A$. More precisely, letting $\tau = \{ k \ge 0 \ | \ S(k) \notin  A \}$, we have
\begin{equation}\label{grin}
G_{A} (x, y) = E^{x} \Big( \sum_{k = 0}^{\tau} {\bf 1} \{ S(k) = y \} \Big). 
\end{equation}
Let us also mention that later in this work we will slightly abuse the same symbols $G ( x, y )$ and $G_{A} (x, y)$ to stand for the Green's functions for a simple random walk $S = S^{(n)}$ on $2^{-n} \mathbb{Z}^{3}$, and restricted to a subset $A \subset 2^{-n} \mathbb{Z}^{3}$  between $x,y\in 2^{-n} \mathbb{Z}^{3}$. 

Concerning the Green's function on $D (N)$,
we will need the following lemma in Section \ref{sec:4}.
\begin{lem}\label{Green-estimate}
There exists a universal constant $C < \infty$ such that for all $0 < \epsilon_{1} < \epsilon_{2} < 10^{-2}$ and $N \ge 1$ with $\epsilon_{1} < \epsilon_{2}^{4}$ and $1/N < \epsilon_{1}^{4}$, 
\begin{equation}
 \Big| \frac{G_{D(N)} (x, y) }{G_{D(N)} (x', y') } - 1 \Big| \le C   \epsilon_{1}^{3/4} \quad \mbox{ and as a special case, } \quad  \Big| \frac{G_{D(N)} (0, x) }{G_{D(N)} (0, x') } - 1 \Big| \le C  \epsilon_{1}^{3/4}, \label{green-1}
 \end{equation}
as long as four points $x, x', y, y' \in D(N) $ satisfy the following conditions: 
\begin{itemize}
\item[{\rm (i)}] ${\rm dist} \big( x,  \partial D(N) \big) \ge \epsilon_{2} N $ and ${\rm dist} \big( y,  \partial D(N) \big) \ge \epsilon_{2} N $,

\item[{\rm (ii)}] $|x- y| \ge \epsilon_{2} N$ and $|x - y| \le 10 \min \{ N - |x|, N - |y|  \}$, 

\item[{\rm (iii)}] $|x- x'| \le \epsilon_{1} N $ and $|y - y' | \le \epsilon_{1} N $.
\end{itemize}
Moreover, if $x,y\in D(N)$ satisfies conditions (i) and (ii) above, then
\begin{equation}\label{green-2}
G_{D(N)}(x,y) \asymp |x-y|^{-1}.
\end{equation}
\end{lem}
\begin{proof}
We first note that it suffices to work with $(x,y)$ and $(x,y')$. Recall that $T_{N}$ stands for the first time that $S$ exits from $D(N)$.
We then observe that
$$
G_{D(N)}(x,y)=G(x,y)-H_{D(N)}(x,y)\mbox{ where }H_{D(N)}(x,y):=\sum_{z\in\partial D(N)} P^x(S (T_{N}) =z) G(z,y),
$$
and a similar decomposition exists for $G_{D(N)}(x,y')$. 
By \eqref{lawler-green},
$$
\left|\frac{G(x,y)}{G(x,y')}-1\right|\leq O\left( \frac{\epsilon_1}{\epsilon_2}\right)\mbox{ and }\left|\frac{G(z,y)}{G(z,y')}-1\right|\leq O\left( \frac{\epsilon_1}{\epsilon_2}\right) \mbox{ uniformly for all $z\in\partial D(N)$}.
$$
Hence, $|H_{D(N)}(x,y)/H_{D(N)}(x,y')-1|\leq O(\epsilon_1/\epsilon_2)$ as well.
Also, by the  Harnack principle \eqref{Harnack} and the assumption on the location of $x,y$, we see that
$G_{D(N)}(x,y)\geq c G(x,y)$ for some universal constant $c>0$.
The claim \eqref{green-1} hence follows as an easy corollary.
The claim \eqref{green-2} is a quick corollary of the asymptotics of Green's function in the whole space (see \eqref{green-entire}) and Harnack principle \eqref{Harnack}.
\end{proof}

\subsection{Metric spaces and weak convergence}\label{metric}
In this subsection, we briefly recall various metric spaces and review the weak convergence of random measures and random sets.

We let $\big( {\cal H} (\overline{\mathbb{D}} ), d_{\text{Haus}} \big)$ be the space of all non-empty compact subsets of $\overline{\mathbb{D}}$ endowed with the Hausdorff metric 
\begin{equation}\label{Haus}
d_{\text{Haus}} (A, B) := \max \Big\{ \sup_{a \in A} \inf_{b \in B} |a-b|,  \  \sup_{b \in B} \inf_{a \in A} |a-b| \Big\} \  \text{ for } A, B \in {\cal H} (\overline{\mathbb{D}} ).
\end{equation}
It is well known that $\big( {\cal H} (\overline{\mathbb{D}} ), d_{\text{Haus}} \big)$ is a compact, complete metric space (see \cite{Henr} for example).

Take two continuous curves $\lambda_{i} : [0, t_{\lambda_{i}}] \to \mathbb{R}^{3}$ ($i=1,2$) where $t_{\lambda_{i}} \in [0, \infty)$ stands for the time duration of $\lambda_{i}$.  We then define their  {\it supremum distance} $\rho (\lambda_{1}, \lambda_{2} )$ by 
\begin{equation}\label{rho-metric}
\rho (\lambda_{1}, \lambda_{2} ) = |t_{\lambda_{1}} - t_{\lambda_{2}} | + \max_{0 \le s \le 1} \big| \lambda_{1} ( s t_{\lambda_{1}} ) - \lambda_{2} ( s t_{\lambda_{2}} ) \big|.
\end{equation}
Using this metric $\rho$, we set $\big( {\cal C} (\overline{\mathbb{D}} ), \rho \big)$ for the space of continuous curves $\lambda : [0, t_{\lambda}] \to \overline{\mathbb{D}} $ with finite time duration $t_{\lambda} \in [0, \infty)$.
It is easy to show that $\big( {\cal C} (\overline{\mathbb{D}} ), \rho \big)$ is a separable metric space (see  \cite[Section 2.4]{KS}).

We next define $\big( {\cal C} , \chi \big)$ as the space of continuous curves $\lambda : [0, \infty)  \to \mathbb{R}^{3} $ with $\lambda (0) = 0$ and $\lim_{ t \to \infty} | \lambda (t) | = \infty$ where the metric $\chi$ is defined by
\begin{equation}\label{eq:chirep-1-1}
\chi (\lambda^{1}, \lambda^{2} ) = \sum_{r=1}^{\infty} 2^{-r} \min \Big\{ \rho \big( \lambda^{1}_{r}, \lambda^{2}_{r} \big) , 1 \Big\} \  \text{ for } \lambda^{i} \in {\cal C},
\end{equation}
with $\rho$ defined as \eqref{rho-metric} and $\lambda^{\bullet}_{r}$, $r > 0$  defined as follows:
\begin{itemize}

\item For $r > 0$  and $\lambda \in {\cal C}$, we write $\tau^{\lambda}_{r}=\tau_{r}(\lambda)$ for the first time that $\lambda$ exits from $\mathbb{B}_{r} $ the open cube with side length $2r$ centered at the origin. Let 
\begin{equation}\label{trunc}
\lambda_{r} : t \in [0, \tau^{\lambda}_{r}] \mapsto \lambda (t) \in \overline{\mathbb{B}}_{r}
\end{equation}
be the truncation of $\lambda$ up to $\tau^{\lambda}_{r}$.

\item For $r  > 0$, we let 
\begin{equation*}
\text{${\cal C}_{r} =$  the space of continuous curves $\lambda : [0, t_{\lambda}] \to  \overline{\mathbb{B}}_{r}$ with $t_{\lambda} \in [0, \infty)$}.
\end{equation*}
 Note that $\lambda_{r} \in {\cal C}_{r}$ for $\lambda \in {\cal C}$.

\end{itemize}
It is easy to check that both $({\cal C}, \chi)$ and $({\cal C}_{r}, \rho)$ are separable metric spaces  (see \cite[Section 2.4]{KS} for this). We also mention that $({\cal C}_{r}, \rho)$ is not complete since the sequence of curves $\{ \lambda^{n} \}_{n \ge 1} \subset {\cal C}_{1}$ defined by 
$$\lambda^{n} (t) = (2^{n} t, 0, 0) \ \ \text{ for } \ 0 \le t \le t_{\lambda^{n}} := 2^{-n},$$
is Cauchy while not convergent.

\medskip

We now turn to weak convergence of random measures. 
Let $(M, d)$  be a metric space with its Borel sigma algebra ${\cal B} (M)$. We denote the space of all probability measures on $ \big( M,  {\cal B} (M) \big) $ by ${\cal P} (M)$, where the space ${\cal P} (M)$ is equipped with the topology of the weak convergence. For a subset $A \subset M$ and $\epsilon > 0$, we write $A_{\epsilon} = \{  x \in M \ | \ \exists y \in A \text{ such that } d(x,y) < \epsilon \}$ for the $\epsilon$-neighborhood of $A$. The Prokhorov metric $\pi : {\cal P} (M)^{2} \to [0,1]$ is defined by 
\begin{equation}\label{prokhorov}
\pi (\mu, \nu) = \inf \Big\{ \epsilon > 0 \ | \ \mu (A) \le \nu (A_{\epsilon} ) + \epsilon \text{ and } \nu (A) \le \mu (A_{\epsilon} ) + \epsilon \text{ for all } A \in {\cal B} (M) \Big\}.
\end{equation}
It is well known that if $(M, d)$ is separable, convergence of measures in the Prokhorov metric is equivalent to weak convergence of measures. Therefore, $\pi$ is a metrization of the topology of weak convergence on ${\cal P} (M)$. 
It is also well known that if $(M, d)$ is a compact metric space, the metric space $\big( {\cal P} (M), \pi \big) $ is compact.
Finally, if $(M,d)$ is a complete, separable metric space, the metric space $\big( {\cal P} (M), \pi \big) $ is complete.

Finally, we discuss random subsets of $\mathbb{R}^d$. To characterize the speed of weak convergence of random closed sets with respect to the $d_{\text{Haus}}$-metric, we introduce the so-called  L\'evy metric as follows. For $X,Y$  two random closed subset of $\mathbb{D}$, writing $T_X(K)=P(X\cap K\neq\emptyset)$ for any closed $K\subset\mathbb{D}$, we define their distance under  L\'evy metric as
\begin{equation}\label{eq:Levy}
d_{\text{L\'evy}}(X,Y)=\inf \Big\{\epsilon\ \Big| \ T_X(K)\leq T_{Y+D(\epsilon)}(K)+\epsilon,\ T_Y(K)\leq T_{X+D(\epsilon)}(K)+\epsilon, \;\forall\mbox{ closed }K\subset\mathbb{D}\Big\}.
\end{equation}
It is well-known that the convergence under L\'evy metric is equivalent to the weak convergence of random closed sets with respect to the $d_{\text{Haus}}$-metric; see e.g.\ \cite{Molchanov} for reference.

\section{Loop-erased random walk}\label{LERW-intro}
In this section, we review the definition and various properties of loop-erased random walk (esp.\ in three dimensions).
\subsection{Definition and domain Markov property}\label{DMPLEW}
We begin with the definition of the chronological loop-erasure of a path. Take a path $$\lambda = [\lambda (0), \lambda (1), \cdots , \lambda (m) ] \subset\mathbb{Z}^{d}$$ with finite length $m$. We define its loop-erasure $\text{LE} (\lambda )$ in the following way. Let $t_{0} = \max \{ k \le m \ | \ \lambda (k) = \lambda (0) \}$ and for each $i \ge 1$, write $t_{i} = \max \{ k \le m \ | \ \lambda (k) = \lambda ( t_{i-1} + 1 ) \}$. Set $l = \min \{ i \ | \ t_{i}  = m \}$. Then, $\text{LE} (\lambda)$ the loop-erasure of $\lambda$ is defined through
\begin{equation}\label{LEP}
\text{LE} (\lambda) = [\lambda  (t_{0}), \lambda (t_{1} ), \cdots , \lambda  (t_{l} ) ].
\end{equation} 
Note that $\text{LE} (\lambda) \subset \lambda$ is a simple path satisfying that $\text{LE} (\lambda) (0) = \lambda (0)$ and $\text{LE} (\lambda) (l) = \lambda (m) $.

We now turn to the loop erasure of a random walk. Let $S$ be the SRW on $\mathbb{Z}^{d}$ and write $T < \infty$ for some (random or non-random) time. We call such loop-erasure $\text{LE} \big( S[0, T] \big) $ a loop-erased random walk (LERW). If $d \ge 3$, since SRW is transient, we can consider the loop-erasure of $S[0, \infty)$. In this case, we call $\text{LE} \big( S[0, \infty) \big)$ the infinite loop-erased random walk (ILERW).

\medskip

The domain Markov property (DMP) is a key feature of LERW. Let $\gamma$ be the loop-erasure of a simple random walk started from somewhere in a domain $A\subset\mathbb{Z}^3$ and stopped at $\tau_A$, the time the SRW exits $A$.  Take two simple paths $\lambda[0,m]\subset A$ and  $\lambda'[0,r] \subset A \cup \partial A$ satisfying that $ \lambda' (0) = \lambda (m)$, $ \lambda'[0,r-1] \subset A$,  $\lambda' (r) \in \partial A$ and $\lambda \oplus \lambda'$ is still a simple path. Now let $R$ be the random walk started from $\lambda (m)$ and conditioned that $R[1, \tau_{A}] \cap \lambda = \emptyset$, then 
\begin{equation}\label{DMP}
P \Big( \gamma = \lambda \oplus \lambda'  \ \Big| \ \gamma [0, m] = \lambda \Big) = P \Big( \text{LE} \big( R [0, T_{A} ] \big) = \lambda' \Big).
\end{equation}
(See e.g.\ \cite[Proposition 7.3.1]{Lawb} for reference.) Namely, the DMP states that conditioned on $\gamma [0, m] = \lambda$, the distribution of the rest of $\gamma$ is given by the loop-erasure of $R$, a random walk conditioned to avoid $\lambda$.

{The loop-erased random walk also satisfies various reversibility properties. For brevity, we only recall the following variant. Let $A$ be a finite subset of $\mathbb{Z}^d$ and $x,y\in A$. Let $S_{x,y}$ and $S_{y,x}$ be two random walk started from $x$ and $y$ respectively conditioned to hit $y$ and $x$ resp.\ (and stopped there) before exiting $A$. Then one has}
\begin{equation}\label{eq:reversibility}
\LE(S_{x,y})\overset{\rm d}{=}\big(\LE(S_{y,x})\big)^{\cal R}.
\end{equation}

\subsection{Properties of LERW}\label{sec:basicproperties}
In this subsection, we recall various key properties of the loop-erased random walk. Note that for brevity we will restrict our review to 3 dimensions.

\medskip


We start with asymptotic independence, which heuristically means that the behavior of the beginning part of an LERW is almost independent of what happens far away. 
In this work, we will repeatedly use the following versions of this property. 



\begin{lem}[Asymptotic Independence]\label{asymp-indep}
Let $S$ be a simple random walk on $\mathbb{Z}^3$. 
For $N\geq 1$, consider the LERW $\lambda=\LE(S[0,T_N])$ and the ILERW $\lambda^\infty=\LE(S[0,\infty))$. For all $M\leq N/4$, 
\begin{equation}\label{eq:LERWvsILERW}
P\big[\lambda[0,T_M]\neq  \lambda^\infty[0,T_M]\big] = 1+O(M/N),
\end{equation}
and hence for any path $\eta$ from the origin to $\partial D(M)$,
\begin{equation}
P\big[\lambda[0,T_M]=\eta\big]= (1+O(M/N))P\big[ \lambda^\infty[0,T_M]=\eta\big].
\end{equation}

Fix $v,w\in D(N)$ and let $R$ be a random walk started at $w$ conditioned to hit $v$ (and stop there) and before exiting $D(N)$. We write $P^{\circ}$ for its probability measure. Write $l = |v -w|$ and pick $\epsilon \in (0, \frac{1}{9})$. For a path started at $w$ and ended at $v$, write $\sigma_v:=\max\{t\leq \Diamond: \lambda(t)\in \partial D (v, \epsilon l)\}$ for the last time that it passes through $\partial D (v, \epsilon l)$. 
 For any paths $\eta^{1}$ and $\eta^{2}$, it follows that 
\begin{equation}\label{eq:asymp-indep}
P^{\circ} \Big( R [\sigma_{v},\Diamond] =  \eta^{1},\LE \big( R  \big) [0, T_{w,\epsilon l}]  = \eta^{2} \Big) 
= P^{\circ}  \Big( R [\sigma_{v},\Diamond] = \eta^{1} \Big)P^{\circ} \Big(  \LE \big( R  \big) [0, T_{w,\epsilon l}]  = \eta^{2}  \Big) \big( 1 + O ( \epsilon ) \big). 
\end{equation}
Moreover, 
\begin{equation}\label{eq:asymp-variant}
     P^{\circ}\big[\big(  \LE \big( R  \big) [0, T_{w,\epsilon l}]  = \eta^{2}  \big)\big]= P\big[\lambda^\infty[0,T_{\epsilon l}]=\eta^{2}\big](1+O(\epsilon)).
\end{equation}
All claims are still true if first exit times of balls (i.e., $T_\cdot$) are replaced by those of cubes.
\end{lem}
We omit the proof as the first claim \eqref{eq:LERWvsILERW} is a direct corollary of the domain Markov property of LERW and the hitting probability estimate of SRW \eqref{srwbound} and the proof of the second claim \eqref{eq:asymp-indep} is quite similar to that of Proposition 4.6 of \cite{Mas}.

\medskip
We now turn to the separation lemma for LERW. The separation lemma is a key tool for the study of critical models (LERW, critical planar percolation, cut points of Brownian motion, etc.). Roughly speaking, its moral is that conditioning on the event that the object of interest (say the curve of LERW or the interface of percolation exploration process) passes through a vertex or some mesoscopic region in a certain manner, then there is a uniformly positive probability that its behavior on and near the boundary of a larger region is ``regular'' in the sense that the distance between both ends of the curve is of the same order as the diameter of the region. See also \cite{DGLZ}, \cite{GLQ}, \cite{GPS},  \cite{LawlerNI}, \cite{Lawcut}, \cite{LawBM} and \cite{Lawrecent}  for instances of the separation lemma on various models.


 To study the law of the LERW near a point it crosses, it is very natural to consider a pair of walks, one LERW and one SRW and then condition that they do not intersect up to some scale. We consider the space of traces of a pair of LERW and SRW with initial configurations in three dimensions. More precisely, take $ k \ge 1$ and define 
\begin{equation}\label{LAMBDAbar-k}
{\Lambda}_{k} := \{ (\gamma, \lambda , \Theta_\gamma, \Theta_\lambda) \ | \ (\gamma, \lambda,  \Theta_\gamma, \Theta_\lambda) \text{ satisfies the conditions (i), (ii), (iii) and (iv)}  \},
\end{equation} 
where 
\begin{itemize}
\item[(i)] $\gamma$ is a simple path in $\mathbb{Z}^{3}$ satisfying that $\gamma [0, \Diamond -1 ] \subset D (k) $ and $\gamma \big( \Diamond \big) \in \partial D (k) $. 

\item[(ii)] $\lambda$ is a path in $\mathbb{Z}^{3}$ satisfying that $\lambda [0, \Diamond -1 ] \subset D (k) $ and $\lambda \big( \Diamond \big) \in \partial D (k) $.

\item[(iii)] 
$\text{$\gamma [0, \Diamond ] \cap \lambda [1, \Diamond] = \emptyset$.}$

\item [(iv)] $\Theta_\gamma$ and $\Theta_\lambda$ are subsets of $D(k)$, such that
$$
\gamma \cap \Theta_\gamma=\big\{\gamma(0)\big\};\quad\quad\lambda\cap \Theta_\lambda=\big\{\lambda(0)\big\}.
$$
%
\end{itemize}

Take $(\gamma, \lambda, \Theta_\gamma, \Theta_\lambda)  \in \Lambda_{k}$. Write $x = \gamma \big(\Diamond \big) $ and $y =  \lambda \big( \Diamond \big)$ for their endpoints. Let 
\begin{itemize}
\item $R^{1} = R^{1}_{\gamma}$ be the random walk on $\mathbb{Z}^{3}$ started at $x$ and conditioned that $R^{1} [1, \infty ) \cap (\gamma \cup \Theta_\gamma) = \emptyset$;
\item  $R^{2} = R^{2}_{\lambda}$ be the random walk started at $y$ conditioned that $R^{2} [1, \infty ) \cap (\{y\} \cup \Theta_\lambda) = \emptyset$.
\end{itemize}
With slight abuse of notation we still write $P$ for the joint law.

 Let $m \geq 2k$ (here 2 is can be any prefixed arbitrary  constant $C>1$). Write $T^{1}_{m}$ for the first time that $\text{LE} ( R^{1} )$ exits from $D (m)$. Similarly, set $T^{2}_{m}$ for the first time that $R^{2}$ exits from $D (m)$. Let
 \begin{equation}\label{NON-INT-EVENT}
 A_{(\gamma, \lambda,  \Theta_\gamma, \Theta_\lambda)}^{m} = \Big\{ \Big( \gamma \oplus \text{LE} (R^{1}_{\gamma}) [0, T^{1}_{m} ] , \lambda \oplus R^{2}_{\lambda} [0, T^{2}_{m}] , \Theta_\gamma, \Theta_\lambda \Big) \in \Lambda_{m} \Big\}
 \end{equation}
 be the non-intersection event of interest.

%
%
%
%
%
%
%
%
In fact, conditioned on the event $A^{m}_{(\gamma, \lambda,\Theta_\gamma, \Theta_\lambda)}$ (recall that $m\geq 2k$) have a good chance of being reasonably far apart even if the end point of $\gamma$ is very close to that of $\lambda$. To be more precise, let 
\begin{equation}
D_{m} (\gamma, \lambda ) =  m^{-1} \min \Big\{ {\rm dist} \Big( \gamma \big( \Diamond\big), \lambda \Big), {\rm dist} \Big( \lambda \big( \Diamond \big), \gamma \Big)  \Big\},
\end{equation}
for $(\gamma, \lambda ,\Theta_\gamma, \Theta_\lambda) \in \Lambda_{m}$. We then have
\begin{lem}[The Separation Lemma]\label{lem:mainseparationlemma}   
There exists a universal constant $c > 0$ such that for all $m \ge 2k$ and all $(\gamma, \lambda ,\Theta_\gamma, \Theta_\lambda ) \in \Lambda_{k}$
\begin{equation}\label{SEPARATION-LEMMA-1}
P \Big\{ D_{m} \Big( \gamma \oplus \text{LE} (R^{1}_{\gamma}) [0, T^{1}_{m} ], \lambda \oplus R^{2}_{\lambda} [0, T^{2}_{m}] \Big) \ge c \ \Big| \ A^{m}_{(\gamma, \lambda ,\Theta_\gamma, \Theta_\lambda)} \Big\} \ge c.
\end{equation}
\end{lem}
This lemma is proved by applying essentially the same argument as the proof of \cite[Theorem 6.5]{S}, which treats the special case with empty initial configurations. We hence omit the proof.

For convenience of notation, we will call $(\gamma, \lambda,\Theta_\gamma, \Theta_\lambda ) \in \Lambda_{m}$ {\bf well-separated}, if $D_m(\gamma, \lambda )>c$ for $c>0$ from \eqref{SEPARATION-LEMMA-1}. This estimate will be used many times in the present article.


\medskip

We now turn to couplings. 
%
%
Take $m > l \ge 2k$. Let 
 \begin{equation}\label{induced-meas}
 \mu^{k,l,m}_{(\gamma, \lambda, \Theta_\gamma, \Theta_\lambda)} (\eta^{1}, \eta^{2}  ) := P \Big( \text{LE} (R^{1}_{\gamma}) [T^{1}_{l}, T^{1}_{m}] = \eta^{1}, \  R^{2}_{\lambda} [T^{2}_{l},  T^{2}_{m}] = \eta^{2}  \ \Big| \ A_{(\gamma, \lambda,\Theta_\gamma, \Theta_\lambda)}^{m}  \Big)
 \end{equation}
 be the distribution of $\big( \text{LE} (R^{1}_{\gamma} ), R^{2}_{\lambda} \big)$ conditioned on the event $ A_{(\gamma, \lambda,\Theta_\gamma, \Theta_\lambda)}^{m}$. 
The following bound on the total variation distance on measures with different initial configurations can be regarded as a generalization of the asymptotic independence and allows us to construct couplings of pairs of non-intersecting LERW and SRW with different initial configurations. Such couplings will be intensively used in Section \ref{sec:4} when we analyze point-crossing and cube-crossing events.


\begin{lem}[Bound on the total variation]\label{lem:maincoupling}
    There exist universal constants $c, \delta > 0$ such that for all $k < l < m$ and any, and $(\gamma_i, \lambda_i, \Theta_{\gamma_i}, \Theta_{\lambda_i})\in\overline{\Lambda}_k$, $i=1,2$,
\begin{equation}\label{coupling2}
\parallel \mu^{k,l,m}_{(\gamma_1, \lambda_1, \Theta_{\gamma_1}, \Theta_{\lambda_1})} - \mu^{k,l,m}_{(\gamma_2, \lambda_2, \Theta_{\gamma_2}, \Theta_{\lambda_2})} \parallel \le c \Big( \frac{k}{l} \Big)^{\delta}.
\end{equation}
\end{lem}
This lemma is proved by applying essentially the same argument as in  \cite[Section 4]{Escape} (see also Remarks 4.4 and 4.8 therein for discussions and \cite{Lawrecent} for the first proof of this type of coupling), which treats the special case of $\Theta_{\gamma_i}=\{\gamma_i(0)\}$ and $\Theta_{\lambda_i}=\lambda_i(0)$, $i=1,2$. We hence omit the proof. 
%
%


\subsection{Escape probabilities}\label{sec:Es}
In this subsection, we discuss the asymptotics of the probability that independent LERW and SRW do not intersect up to some scale.

We consider two independent SRW's $S^{1}$ and $ S^{2}$ on $\mathbb{Z}^{3}$ started at the origin. Let $T^{i}_{N}$ be the first time that $S^{i}$ exits from $D(N)$. Write 
\begin{equation}\label{lambda-1-n}
\lambda^{1}_{N} =\text{LE} \big(  S^{1} [0, T^{1}_{N} ] \big).
\end{equation}
In \cite{Escape}, the following non-intersection probabilities are considered.  For $1 \le M \le N$, we define
  \begin{align}
  & \Es (N) := P \Big( \lambda^{1}_{N} \cap S^{2} [1, T^{2}_{N} ] = \emptyset \Big), \notag \\
  & \Es (M, N) := P \Big( \lambda^{1}_{N} [T_{M,N}, \Diamond] \cap S^{2} [0, T^{2}_{N} ] = \emptyset \Big), \label{ESCAPE}
  \end{align}
 (recall that $T_{M,N}$ denotes the last time that the LERW $\lambda^{1}_{N}$ stays in $D(M)$). 
As a shorthand, later in this work when we are working on the rescaled lattice $2^{-n}\mathbb{Z}^3$, we will write 
\begin{equation}\label{eq:rescaledEs}
\Es_n(a,b):=\Es(a2^{n}, b2^{n}).
\end{equation}
  
  It is proved in  \cite[Theorem 1.2]{Escape} that there exist universal constants $c, \delta > 0$ such that for all {\bf integer} $n \ge 1$,
  \begin{equation}\label{COR-escape}
  \Es ( 2^{n} ) = c 2^{- ( 2- \beta ) n } \Big( 1 + O \big( 2^{- \delta n } \big) \Big),
  \end{equation}
for some constant
\begin{equation}\label{growth}
\beta \in \left(1, \frac{5}{3} \right]
\end{equation}  
which is called the {\bf growth exponent} of 3D LERW. (In fact, \eqref{COR-escape} is true for all real $n\geq 0$, just as the case for the asymptotics of $f_n$; see the discussion below \eqref{eq:fnpreview} and below \eqref{eq:fndef} for more discussions.) 
Furthermore, in  \cite[Corollary 1.3]{Escape}, it is proved that for all positive reals $N\geq M\geq 1$
  \begin{equation}\label{ES}
  \Es (N) \asymp N^{- (2 - \beta )}, \ \Es (M, N) \asymp \Big( \frac{N}{M} \Big)^{- ( 2- \beta ) }.
  \end{equation}
  
  These non-intersection probabilities are defined for the infinite loop-erased random walk similarly. We write $\lambda_{\infty}^{1} = \LE \big( S^{1} [0, \infty ) \big) $ for the infinite LERW. Take $1 \le M \le N$. We set 
 \begin{align*}
  & \Es^{\infty} (N) := P \Big(  \lambda_{\infty}^{1} [0, T^1_N] \cap S^{2} [1, T^{2}_{N} ] = \emptyset \Big), \\
  & \Es^{\infty} (M, N) := P \Big( \lambda_{\infty}^{1} [T_{M,N}, T^1_N] \cap S^{2} [0, T^{2}_{N} ] = \emptyset \Big). 
  \end{align*}
By adapting a similar technique used in the proof of \cite[Propositions 6.8 and 6.10]{S} where it is proved $\Es (N) \asymp \Es (M) \Es (M,N)$, one can prove that $\Es^{\infty} (N) \asymp \Es^{\infty} (M) \Es^{\infty} (M,N)$.
Combining this observation with \cite[Proposition 6.7]{S} and \eqref{ES}, we have 
\begin{equation}\label{ES-infty}
\Es^{\infty} (N) \asymp N^{- (2 - \beta )}, \ \Es^{\infty} (M, N) \asymp \Big( \frac{N}{M} \Big)^{- ( 2- \beta ) }.
\end{equation}

%

\subsection{Hittability and quasi-loops}\label{qloops}
Take a discrete path $\lambda \subset \mathbb{Z}^{3}$. Let $R$ be the SRW on $ \mathbb{Z}^{3}$. Suppose that the distance between $R (0)$ and $\lambda$ is much smaller than the diameter of $\lambda$. Does $R$ intersect with $\lambda$ before exiting a large ball centered at $R (0)$ with high probability? The answer is ``no" in general in three dimensions. For example, if we take $\lambda$ as a straight line, it is unlikely for $R$ to hit $\lambda$. The reader should recall that Brownian motion does not hit a straight line in three dimensions.   
However, since the growth exponent of 3D LERW is strictly greater than 1, if $\lambda$ in the previous paragraph is the trace of LERW, we expect that it is very likely for $R$ to intersect $\lambda$ immediately. In other words, with high probability the trace of LERW is a ``{\bf hittable}'' set.  

\begin{figure}[h]
\begin{center}
\includegraphics[scale=0.6]{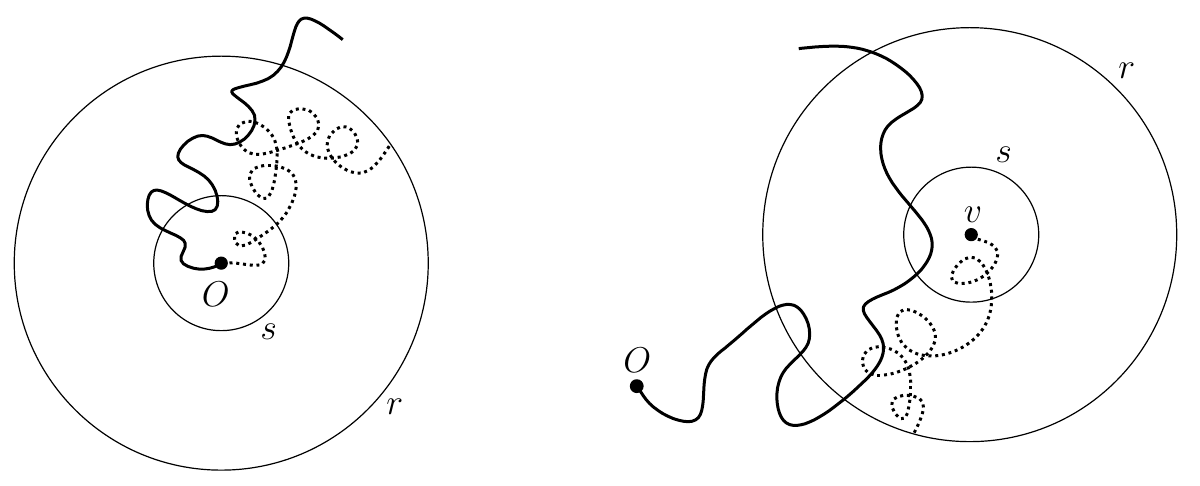}
\caption{Illustration for the hittability of the LERW. The thick solid curve stands for $\lambda^T$ while the thick dotted curve represents $R \big[ 0, T_{0, r}^{R} \big] $ (left) and $R \big[ 0, T_{v, r}^{R} \big] $ (right).}\label{hit-hit}
\end{center}
\end{figure}

We now explain this in a more precise manner. Recall that $S$ stands for the SRW started at the origin on $\mathbb{Z}^{3}$ and let 
 $\lambda^T:=\text{LE}  \big( S [0, T] \big)$ for any $T>0$. 
 It is proved in \cite[Lemma 4.6]{Koz} that for all $K \ge 1$ there exist $a > 0$ and $C < \infty$ that for any $1\leq s< r$ and 
\begin{align}\label{holder-1}
&P \bigg(  \text{ There exists } T \text{ such that }  \lambda^T \not\subset D(r), \text{ and }  \notag \\
& \ \ \ \ \ P_{R}^{0} \Big( R \big[ 0, T_{0, r}^{R} \big]  \cap \lambda^T \cap \big( D(r) \setminus D(s) \big) = \emptyset \Big) > \Big( \frac{s}{r} \Big)^{a} \bigg) \le C \Big( \frac{s}{r} \Big)^{K},
\end{align}
 where we  write $R$ for a simple random walk on $\mathbb{Z}^{3}$ which is independent of $S$, and $P_{R}^{x}$ for the probability law of $R$ assuming that $R(0) = x$. Loosely speaking, the inequality \eqref{holder-1} ensures that the loop-erased random walk $\text{LE}  \big( S^{(n)} [0, T] \big)$ restricted to $ \big( D ( r) \setminus D(s) \big) $ is a ``{hittable}'' set in the sense that $R$ hits it before leaving $B (0, r)$ with high probability; see the left picture of Figure \ref{hit-hit}.
 
What about hittability of the LERW at other parts? 
To answer this question, \cite[Lemma 3.3]{SS} shows that  for all $K \ge 1$ there exist $a > 0$ and $C < \infty$ that for any $1\leq s< r$, $v \notin D(r)$ 
\begin{align}\label{holder-2}
&P \bigg(  \text{ There exists } T \text{ such that }  \lambda^T \cap D(v, s) \neq \emptyset, \text{ and }  \notag \\
& \ \ \ \ \ P_{R}^{v} \Big( R \big[ 0, T_{v, r}^{R} \big]  \cap \lambda^T [0, T_{v,s}] \cap \big( D (v, r) \setminus D (v, s) \big) = \emptyset \Big) > \Big( \frac{s}{r} \Big)^{a} \bigg) \le C \Big( \frac{s}{r} \Big)^{K}.
\end{align}
This inequality \eqref{holder-2} guarantees that the LERW restricted to $\big( D (v, r) \setminus D (v, s) \big)$ is also hittable with high probability assuming it hits $D (v, s)$; see the right picture of Figure \ref{hit-hit}.

We will next recall the notion of a {\bf quasi-loop}. Take $1\leq s < r$ and $N \ge 1$. For a path $\lambda $ in $\mathbb{Z}^{3}$ and $w \in  \mathbb{Z}^{3}$, we say $\lambda$ has a $(s, r)$-quasi-loop at  $w$  if there exist $i_{1} < i_{2}$ such that $\lambda (i_{1}), \lambda (i_{2}) \in D (w, s)$ and $\lambda [i_{1}, i_{2} ] \not\subset D (w, r ) $.  We set ${\rm QL} (s, r ; \lambda )$ for the set of all such $w$'s. 
For $N \ge 1$, let $\lambda_N:=\LE(S[0,T_N])$. It is shown in \cite[Theorem 6,1]{SS} that
\begin{lem}[Rareness of quasi-loops]\label{lem:noquasiloop}
   There exist $0<C,L < \infty$ and $b > 0$ such that for all $\theta >0$ and $N \ge 1$, 
\begin{equation}\label{eq:noquasiloop}
P \Big( {\rm QL} \big(\theta^{L}, \theta; \lambda_{N} \big) \neq \emptyset \Big) \le C \theta^{b}.
\end{equation} 
\end{lem}

\section{Scaling limit of 3D LERW}\label{sec:LERWlimit}
In this section we briefly review our knowledge on the scaling limit of 3D LERW and discuss asymptotics of the one-point function and tightness of the LERW in terms of occupation measure and curve under natural parametrization. 

From here onwards, for the convenience of discussion, we always work on the rescaled lattice $2^{-n} \mathbb{Z}^{3}$ for some {\bf real} $n\geq 0$. Let $S^{(n)}$ be the simple random walk on $2^{-n} \mathbb{Z}^{3}$ and let $T$ be the first time $S$ exits $\mathbb{D}$. 
Let 
\begin{equation}
    \gamma_{n}:=\LE(S^{(n)}[0,T])
\end{equation} stand for the LERW on $2^{-n} \mathbb{Z}^{3}$ in $\mathbb{D}$ 
and regard $\gamma_{n}$ as a random element of ${\cal H} (\overline{\mathbb{D}} )$ (recall its definition in 
in \eqref{Haus}).

\subsection{Existence of the scaling limit}\label{SCALING}
  In \cite{Koz}, Kozma proves that $\gamma_n$ converges 
  with a polynomial speed in L\'evy metric, along dyadic scales (this corresponds to taking integer $n$'s). More precisely, it is proved (see Theorems 5 and 6, ibid.) that for $n\in\mathbb{Z}^+$, there exists ${\cal K}$, a random closed subset of $\mathbb{D}$, and $0<c,C<\infty$ such that
\begin{equation}\label{eq:convdyadic}
d_{\text{L\'evy}}(\gamma_{n},{\cal K})\leq C2^{-cn}.
\end{equation}
In the same work, it was also proved (see Section 6.1, ibid.) that the ${\cal K}$ is rotation-invariant as well as scaling-invariant in the following sense. For $r>0$ write 
$$
{\cal K}^r :=r{\cal K}
$$
for the $r$-blowup of ${\cal K}$ and let $D_{r,n}=2^{-n}\mathbb{Z}^3\cap \mathbb{D}$ be the discretization of $\mathbb{D}_r$ and write $\gamma^{r}_{n}$ for the LERW in $D_{r,n}$ started from the origin. Then 
there exist constants $c(r),C(r)\in(0,\infty)$, such that for $n\in \mathbb{Z}^+$,
\begin{equation}\label{eq:convrdya}
d_{\text{L\'evy}}(\gamma^{r}_{n},{\cal K}^r)\leq C2^{-cn},
\end{equation}
In other words, ${\cal K}^r$ is the scaling limit of LERW in the domain $\mathbb{D}_r$ along dyadic scales.

It is noteworthy that the scaling invariance of ${\cal K}$ actually implies that the convergence in \eqref{eq:convdyadic} follows not only along dyadic scales, but also as $m\to\infty$ continuously, which we now state and justify:
\begin{prop}\label{prop:Hanyscale}
There exists $c,C\in (0,\infty)$ such that as $m\to \infty$
\begin{equation}
d_{\text{L\'evy}}(\gamma_{m},{\cal K})\leq Cm^{-c}.
\end{equation}
\end{prop}
\begin{proof}
We observe that it suffices to show the constants $c,C$ in \eqref{eq:convrdya} can be taken uniformly for all $r\in[1,2]$. 

We return to the proofs in \cite{Koz}. Citation of equation numbers, theorems, lemmas all refer to those of that work. Note that
when applying\footnote{Note that in Theorem 5 it is required that ${\cal D}\subset[\frac{1}{4},\frac{3}{4}]^3$, but one can easily replace this requirement by, say, $[-3,3]^3$ through translation and dilation, which will enable us to take $\cal D$ as $\mathbb{D}_r$ (in our notation) for any $r\in[1,2]$,  while affecting constants by a universal factor  only.} Theorems 5 and 6 with $\cal D$ (in the notation of \cite{Koz}, same below) as a ball\footnote{In \cite{Koz} it is also required that $\cal D$ is a polyhedron, however this requirement is only relevant in Lemma 2.13, which is sub-sequently only applied in Lemmas 5.9 and 5.10. As we can easily build a version of Lemma 2.13, for $\cal D$ being a ball, one can also relax this requirement, allowing $\cal D$ to be a ball for the domain on which we consider scaling limits of LERW.} $\mathbb{D}_r$, the constants may depend on the domain $\cal D$ and the starting point $a$ (in the notation of \cite{Koz}). However, the dependence only originates from (132), (134) and (138). For $r\in[1,2]$, we can easily make uniform choices for constants in (132). For constants in (134), we see that for $r\in[1,2]$, one can write down a version of Lemma 2.3 for $\mathbb{R}^3\setminus \mathbb{D}_r$ for all $r\in[1,2]$ with universal constants. As $\mathbb{R}^3\setminus \mathbb{D}_r$ has only one connected component, this implies that the constants in (135) (and hence  (137)) can be chosen uniformly, which leads to a uniform choice of constants in (134). The same argument works for (138).
\end{proof}
We remark that similar arguments for the universality of constants in \cite{Koz} can be found below (5.51) of \cite{Escape}.

The following corollary hence follows immediately thanks to the equivalence of weak convergence and L\'evy metric (see the end of Section \ref{metric}).
\begin{cor} \label{cor:wCONV}For $m\in\mathbb{R}^+$,
$$
\gamma_m\overset{\rm w}{\to}{\cal K}\mbox{ as $m\to\infty$}
$$
in $({\cal H}(\overline{\mathbb{D}}),d_{\rm Haus})$.
\end{cor}
  Not much is known about the scaling limit $\cal K$. In  \cite[Theorem 1.2]{SS}, it is shown that ${\cal K}$ is a simple curve almost surely and by adding loops from an independent Brownian loop soup appropriately to $\cal K$ one recovers the trace of a Brownian motion.  In \cite{S2}, it is proved that the Hausdorff dimension of ${\cal K}$ is equal to $\beta$ almost surely. See Remark \ref{rem:mainrem} 2) for more discussions.

\subsection{One-point function}\label{sec:ONE}
In this subsection, we focus on the the asymptotics of the probability that LERW crosses a certain point, colloquially referred to as the ``one-point function''.
Take a point $x \in \mathbb{D} \setminus \{ 0 \}$. We set $x_{n}$ for one of the nearest points from $2^{n} x$ among $\mathbb{Z}^{3}$. Then \cite[Theorem 1.1]{Escape} shows that there exist absolute constants $c, \delta > 0$ and a quantity depending on $x \in \mathbb{D} \setminus \{ 0 \}$ \begin{equation}\label{one-point-def}
c_{x}=c(x) > 0
\end{equation}
such that for all {\bf integer} $n \ge 1$ and $x \in \mathbb{D} \setminus \{ 0 \}$
  \begin{equation}\label{one-point}
  P \Big( x_{n} \in \gamma_n \Big) = c_{x} 2^{- (3 - \beta ) n} \Big( 1 + O \big( d_{x}^{-c} 2^{- \delta n} \big) \Big) \   \ \text{as } n \to \infty,
  \end{equation}
  where 
  $d_{x} = \min \{ |x|, 1- |x| \}$. 
  For positive real $n$ we have that 
\begin{equation}\label{eq:weakeronepoint}
    P \Big( x_{n} \in  \gamma_n \Big) \asymp c_x 2^{- (3 - \beta ) n}.
\end{equation}
In order to properly rescale various quantities in this work, we now define a ``reference scaling factor'' through one-point function as follows. 
Letting $\widehat{x}=(1/2,0,0)$, for any positive real $n$, we write
\begin{equation}\label{eq:fndef}
f_n := 2^{3n} P \big( \widehat{x}_{n} \in  \gamma_n  \big),
\end{equation}
It is obvious that 
\begin{equation}\label{eq:fnasymp}
f_n \asymp 2^{\beta n}
\end{equation} and by \eqref{one-point}, for integer $n$'s, $f_n = C 2^{n} \big(1+O(2^{-cn})\big)$ for some $C,c>0$.
As mentioned below \eqref{eq:fnpreview}, in a recent follow-up work \cite{HTLS}, it is proved that the above claim is true for all real $n>0$. In order to not create circular argument, in this paper we only work with estimates available in this subsection instead of newest ones from \cite{HTLS}. 

\subsection{Tightness of rescaled LERW}\label{sec:tight}
 In this subsection, we briefly review the controls that give tightness of 3D LERW  under two different notions: the renormalized occupation measure and the curve under natural parametrization. These results allow us to take sub-sequential limits of 3D LERW and characterize such limits. 

As introduced in Section \ref{sec:1.1}, let $\mu_n$ be the rescaled occupation measure of $\gamma_n$ and $M_n=\mu_n(\mathbb{D})$ be the total number of steps of $\gamma_n$.
The following (stretched-)exponential tail bounds on $M_n$ is shown in \cite{S}.
\begin{lem}[Exponential tail bound]There exist $c_{1}, c_{2}, \delta > 0$ such that for all $n \ge 1$ and $b \ge 1$ ,
 \begin{equation}\label{exp-tail}
 P \Big(  M_{n} / E (M_{n} ) \in [1/b , b] \Big) \ge 1 - c_{1} \exp \big( - c_{2} b^{\delta} \big).
 \end{equation}
 \end{lem}
As a corollary of \eqref{eq:weakeronepoint}, one has  
 \begin{equation}\label{improve}
 E (M_{n} ) \asymp 2^{\beta n}.
 \end{equation}
  Combining this with \eqref{exp-tail} and weak convergence of  $\gamma_{n}$ with respect to the Hausdorff distance (see Section \ref{SCALING} for this) we see that
\begin{prop}\label{tight}
The family of variables  $ \{ (\gamma_{n}, \mu_{n}) \}_{n \ge 1}$ is tight with respect to the product topology of $({\cal H} (\overline{\mathbb{D} } ), d_{\text{Haus}} )$ and the topology of the weak convergence on ${\cal M}  (\overline{\mathbb{D} } )$.
\end{prop}


We now turn to natural parametrization.
We recall that $\eta_{n} (t )$ 
stands for the time-rescaled LERW as introduced above Theorem \ref{3rd}. 
To prove that $\{ \eta_{n} \}_{n \ge 1}$ is tight with respect to the topology in the metric space $\big( {\cal C} (\overline{\mathbb D}), \rho \big)$ (see Section \ref{metric} for definition), it suffices to give the following bound on the ``modulus of continuity'' for $\{ \eta_{n} \}$, which is a paraphrase of \cite[Proposition 3.7]{Holder}, whose proof we omit for brevity (see Remark \ref{rem:nocircular} for more discussions).
\begin{prop}\label{HC-upper}
Let 
\begin{equation}
    K_n=K_{n,\delta,\epsilon}:=\Big\{\exists x=\eta_n(s), y=\eta_n(t)\mbox{ s.t. } 0\leq t-s\leq \delta\mbox{ and }|x-y|\geq \delta^{1/\beta-\epsilon}\Big\}.
\end{equation}
There exist universal constants $c,C\in(0,\infty)$ such that for all $\delta, \epsilon \in (0,1/10)$ with $\delta^{-\epsilon}>100$.
\begin{equation}\label{eq:HC-upper}
    P(K_n^c) \geq 1- C\delta^{c\epsilon}.
\end{equation}
\end{prop}
Proposition \ref{HC-upper} shows the equicontinuity of $\{ \eta_{n} \}_{n \ge 1 }$. Combined with \eqref{exp-tail}, 
 the tightness of $\eta_n$ is a direct consequence of the Arzel\'{a}-Ascoli theorem (see \cite[Theorem 7.3]{Bil} for this):
\begin{cor}\label{0626}
As a sequence of $\big( {\cal C} (\overline{\mathbb D}), \rho \big)$-valued random variables, $\{ \eta_{n} \}_{n \ge 1 }$ is tight. 
\end{cor}


\begin{rem}\label{rem:nocircular}
Although \cite{Holder} assumes Theorem \ref{3rd} of this work, namely the scaling limit of 3D LERW in natural parametrization,   the proof of \cite[Proposition 3.7]{Holder} nevertheless is purely an analysis of the discrete object and hence there is no risk of circular argument. 

\end{rem}

\section{Setup and preliminary estimates for the key $L^2$-estimate \eqref{KEY} }\label{sec:3}
In this section, we will introduce the setup for the crucial $L^2$-estimate \eqref{KEY}, and give preliminary up-to-constants second moment estimates and prove stability of the one-point crossing probability.

\subsection{Up-to-constants $L^2$-bounds}\label{sec:uptocstL2}

We start with the partition of the space.

\begin{dfn}\label{def:Cubes}
Fix $n\geq1$ (not necessarily an integer) and take an integer $ k \ge 1$. Fix $\overline{x} \in  \mathbb{D}$ and consider the (closed) cube $B =\overline{B}_\infty(\overline{x},2^{-k})= \{ x \in \mathbb{R}^{3}  \ | \ \| x - \overline{x} \|_{\infty}  \le 2^{-k} \}$. 
We assume that $B \subset \mathbb{D}$ and $B$ is ``typical'' in the following sense:
\begin{equation}\label{eq:Breq}
{\rm dist} \big( B, \{0 \} \cup \partial \mathbb{D} \big) > 2^{-k}.
\end{equation}  
With small abuse of notation, we also use $B$ to stand for $B\cap 2^{-n} \mathbb{Z}^3$ as a subgraph of $2^{-n}\mathbb{Z}^3$, whenever this does not cause any confusion. Let \begin{equation}
 M_k=2^{3(k^4-k)}.
\end{equation}
We ``roughly divide'' $B$ into smaller  cubes $B_{1}, B_{2}, \cdots , B_{m_{k}}$ 
in the following way: let
$$B_{i} = \overline{B}_\infty(x_i,2^{-k^4})
$$
stand for the (discrete) cube of side length $2^{-k^{4}+1}$ centered at $x_{i}$, with $x_i\in 2^{-n}\mathbb{Z}^3$, $i=1,\cdots,M_k$, chosen in a way such that $B\subset\bigcup_{i=1}^{m_{k}} (\overline{B_{i}})$ and for any $i\neq j$, $B_i^\circ\cap B_j^\circ=\emptyset$.
Note that such a collection exists but is not necessarily unique. If there are multiple choices, pick one arbitrarily.

Throughout this section and the next one, we write
\begin{equation}\label{eq:rDEF}
\begin{split}
&\epsilon =   2^{-k^{4}+1}  \text{ for side length of  } B_i,\\
&\text{and } r=\min \{ |\overline{x}|,1-|\overline{x}| \} > 2^{-k} \text{ for the distance between $\overline{x}$ and $\{ 0 \} \cup \partial \mathbb{D}$}.
\end{split}
\end{equation}
\end{dfn}

To simplify notation, in this section we write  $S= S^{(n)}$ and $\gamma=\gamma_n$ for short.

Let $X_i$ be the number of points in $B_i \cap 2^{-n} \mathbb{Z}^{3}$ hit by $\gamma$. Note that since $\gamma$ is a simple path, $X_i$ is also the time $\gamma$ spends in $B_i \cap 2^{-n} \mathbb{Z}^{3}$. Then, we let $Y_{i}$ be the indicator function of the event that $\gamma$ hits $B_i$. We set 
\begin{equation}\label{X}
X = \sum_{i=1}^{M_k } X_{i}
\end{equation}
for a ``rough count''\footnote{On the boundaries of $B_i$ there may be double-counts and left-outs.} on the total number of points in $B \cap 2^{-n} \mathbb{Z}^{3}$ hit by $\gamma$ and set
\begin{equation}\label{Y}
Y = \sum_{i=1}^{M_k } Y_{i}
\end{equation}
for the number of cubes among $B_1,\ldots,B_{M_{k}}$ hit by $\gamma$. 

As discussed in Section \ref{sec:1.2}, we are going to give an $L^2$-approximation of $X$ using $Y$ 
in the form of \eqref{KEY} (see Proposition \ref{prop4.6} below for the precise form). To this end, we need to give some preliminary up-to-constants bounds for various second moment estimates.

We write
\begin{align}
& \Big\{ 0 \xrightarrow{\gamma} x \xrightarrow{\gamma} y \Big\} := \Big\{ \gamma \text{ first hits } x \text{ and then hits } y \Big\}\mbox{ and } \label{eq:crosspointdef-0} \\
&\Big\{ 0 \xrightarrow{\gamma} A \xrightarrow{\gamma} A' \Big\} := \Big\{ \gamma \text{ first hits } A \text{ and then hits } A' \Big\},\label{eq:crossboxdef-0-1} 
\end{align}
for points $x \neq y$ and for sets $A, A'$ with $A \cap A' = \emptyset$ as shorthand for point-crossing events and cube-crossing events respectively, and similarly define ``hybrid'' point-cube-crossing events $\{0 \xrightarrow{\gamma} v \xrightarrow{\gamma} A \}$ and $\{0 \xrightarrow{\gamma} A \xrightarrow{\gamma} v \}$ for a point $v$ and a set $A$ such that $v\notin A$.
\begin{prop}\label{prop:bound-box-hit}
For all $k \ge 1$, and $B\subset\frac{2}{3}\mathbb{D}$ satisfying ${\rm dist} \big( B, \{0 \} \big) > 2^{-k}$, it follows that  
\begin{equation}\label{M7}
 P (0 \xrightarrow{\gamma} B_i \xrightarrow{\gamma} B_j ) \asymp \Big( \frac{\epsilon}{r} \Big)^{3-\beta } \Big( \frac{\epsilon}{l} \Big)^{3-\beta },
\end{equation}
where  $l = |x_{i} - x_{j} | $.
Similarly, for $v\in B_i$ and $w\in B_j$, 
 \begin{equation} \label{M7prime} 
P (0 \xrightarrow{\gamma} v \xrightarrow{\gamma} w )\asymp \Big( \frac{2^{-n}}{r} \Big)^{3-\beta } \Big( \frac{2^{-n}}{l} \Big)^{3-\beta }; 
\end{equation}
\begin{equation} \label{M7primeprime}
P (0 \xrightarrow{\gamma} v \xrightarrow{\gamma} B_j )  \asymp P (0 \xrightarrow{\gamma} B_i \xrightarrow{\gamma} w) \asymp \Big( \frac{2^{-n}}{r} \Big)^{3-\beta } \Big( \frac{\epsilon}{l} \Big)^{3-\beta }. 
\end{equation}
If $B \subset  \mathbb{D} \setminus \frac{1}{3} \mathbb{D}$, then
\begin{equation}
 P (0 \xrightarrow{\gamma} B_i \xrightarrow{\gamma} B_j ) \asymp \Big( \frac{\epsilon^{3-\beta } }{r^{1-\beta } } \Big) \Big( \frac{\epsilon}{l} \Big)^{3-\beta }.
\end{equation}
Similar bounds exist for point-crossing and point-cube-crossing events.
\end{prop}
We also record here corresponding classical first moment bounds for comparison. For $B \subset  \frac{2}{3} \mathbb{D}$, $B_i$ one of the sub-cubes of $B$ and $v\in B_i$,
  \begin{equation} \label{20221118}
  P ( v \in \gamma ) \asymp \Big( \frac{2^{-n}}{r} \Big)^{3-\beta };  
 \end{equation}
 \begin{equation}\label{20221118-1} 
 P ( \gamma  \cap B_i \neq \emptyset) \asymp \Big( \frac{\epsilon}{r} \Big)^{3-\beta }.
 \end{equation}




\begin{proof} 
We only prove \eqref{M7} as all other bounds can be proved in the same fashion.

We start with the upper bound, which follows easily from \cite[Theorem 3.1]{S2} (see also \cite[Remark 3.2]{S2}). There it is proved that there exists universal constant $C > 0$ and constant $N_{k}$ depending only on $k$ such that for all $k \ge 1$ and $n \ge N_{k}$
\begin{equation}\label{1116-1}
 P (0 \xrightarrow{\gamma} B_i \xrightarrow{\gamma} B_j  ) \le C P ( \gamma \cap B_i \neq \emptyset ) P^{x_{i}} ( \gamma \cap B_j \neq \emptyset ).
\end{equation}
A variant of \eqref{20221118-1} (see e.g.\  \cite[Lemma 7.1]{SS} for a proof) gives that for large $n \ge N_{k}$  
\begin{equation}\label{1116-2}
 \ P^{x_{i}} ( \gamma \cap B_j \neq \emptyset ) \le C \frac{\epsilon}{l} \Es_n ( \epsilon, l ),
\end{equation}
where $\Es_n (\cdot, \cdot)$ is defined as in \eqref{eq:rescaledEs}.
Combined with \eqref{20221118-1}, by \eqref{COR-escape}, we see that the right hand side of \eqref{1116-1} is bounded above by $ C\big(\epsilon/r\big)^{3 - \beta} \big( \epsilon/l \big)^{3 - \beta} $ which gives the upper bound of  \eqref{M7}.

We now turn to the lower bound. We recall that $T$ stands for the first time that $S$ exits from $\mathbb{D}$. Let $T^{i} = \max \{ t \le T \ | \ S (t) \in B_{i} \}$ be the last time that $S$ stays in $B_i$. We define $T^{j}$ similarly. Decomposing  $S [0, T]$ at $T^{i}$ and $T^{j}$ as in the proof of \cite[Theorem 3.1]{S2},  the probability $ P (0 \xrightarrow{\gamma} B_i \xrightarrow{\gamma} B_j )$ can be bounded from below by 
\begin{equation}\label{M8}
\sum_{y \in \partial_i B_i} \sum_{z \in \partial_i B_j} P^{y} \big( 0 \in S[0, T] \big) P^{y} \big( z \in S[0, T] \big) P (F_{y,z}),
\end{equation}
where for each $y \in \partial_i B_i$ and $z \in \partial_i B_j$ the event $F_{y, z}$ will be defined as in \eqref{eq:Fyzdef} below.

To define the event $F_{y, z}$, we first introduce three independent random walks $Y^{1}$, $Y^{2}$ and $Y^{3}$ such that 
\begin{itemize}
\item $Y^{1}$ starts from $y$ and is conditioned that it hits the origin  before exiting $\mathbb{D}$. We denote the first hitting time to the origin for $Y^{1}$ by $\tau_{0}$, and let 
$$
\lambda :=\LE (Y^{1} [0, \tau_{0} ])
$$
and $$u=\max\{t\leq \Diamond \ | \ \lambda(t)\in B_i\big\}.$$

\item $Y^{2}$ starts from $y$ and is conditioned that it hits $z$  before exiting $\mathbb{D}$. We denote the first hitting time to $z$ for $Y^{2}$ by $\tau_{z}$. We also write
$$s = \min \{t \ge 0 \ | \ \LE \big(Y^{2} [0, \tau_{z} ] \big)(t)\in B_j\big\}.$$

\item $Y^{3}$ is the simple random walk started at $z$ stopped at exiting $\mathbb{D}$. We write $T$ for the first time that  $Y^{3}$ exits from $\mathbb{D}$.

\end{itemize} 
Then,
\begin{equation}\label{eq:Fyzdef}
F_{y,z}:= {\cal A}\cap {\cal B}\cap {\cal C}\cap {\cal D}\cap {\cal E},\mbox{ where}
\end{equation}
\begin{itemize}
\item[]  ${\cal A}:=\big\{u<\inf\{t\geq 0;\lambda(t) \notin B (x_{i}, \frac{l}{8} )\}\big\} \cap \big\{  y \notin Y^{1} [1, \tau_{0} ] \big\}$;
\item[] ${\cal B}:=\big\{\lambda(u,\Diamond] \cap \big( Y^{2} [0, \tau_{z}] \cup Y^{3} [0, T] \cup B (x_{j}, \frac{l}{3} ) \big) = \emptyset \big\}$;

\item[] ${\cal C}:=\big\{Y^{2} [1, \tau_{z} ]  \cap B_i = \emptyset\big\}\cap\big\{Y^{3} [1, T ]  \cap \big( B (x_{i}, \frac{l}{3} ) \cup B_j \big) = \emptyset\big\}$;

\item[] ${\cal D}:=\big\{\mbox{$\exists t_{0} \in [T_{x_{i}, \frac{l}{8}}, T_{x_{i}, \frac{l}{4}}]$ such that  $\big( Y^{2} [0, t_{0}]  \cup B(x_{i}, \frac{l}{8}) \big) \cap  Y^{2} [t_{0}+1, \tau_{z}] = \emptyset$}\big\}$;

\item[] ${\cal E}:=\big\{\LE ( Y^{2} [0, \tau_{z} ] ) [0, s] \cap Y^{3} [0, T] = \emptyset\big\}$.

\end{itemize}
Here $T_{x_{i}, r}$ in the event ${\cal D}$ stands for the first time that $Y^{2}$ exits from $B (x_{i}, r)$. See Figure \ref{1117-1-0} for these events.  

\begin{figure}[h]
\begin{center}
\includegraphics[scale=0.55]{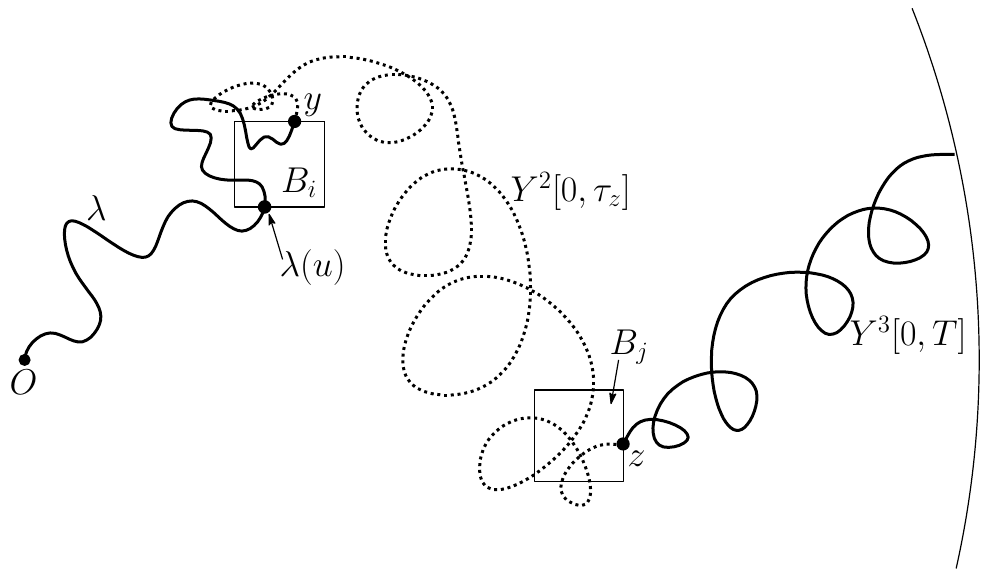}
\caption{Setup for the event $F_{y,z}$. The solid simple curve starting at $y$ stands for $\lambda$. 
The dotted curve is $Y^{2} [0, \tau_{z}]$ which also  starts at $y$. The remaining curve started at $z$ represents $Y^{3} [0, T]$.}\label{1117-1-0}
\end{center}
\end{figure}

As in the proof of \cite[Proposition 3.10]{S2}, decomposing the paths near $B_i$ and $B_j$ and applying the separation lemma (Lemma \ref{lem:mainseparationlemma}) to see that there is a uniformly positive (conditional) probability that the paths near these cubes are well-separated conditioned on non-intersection so that we can produce a lower bound on the probability of $F_{y,z}$,  we have 
\begin{equation*}
P (F_{y,z} ) \ge c \big( \epsilon 2^{n} \big)^{-2} \Es_n (\epsilon , l  ) \Es_n ( \epsilon , r  ) \asymp \big( \epsilon 2^{n} \big)^{-2} \Big( \frac{\epsilon}{l} \Big)^{2-\beta} \Big( \frac{\epsilon}{r} \Big)^{2-\beta}.
\end{equation*}
Since the number of points in $ \partial_i  B_i \cap 2^{-n} \mathbb{Z}^{3}$ is comparable to $\big( \epsilon 2^{n} \big)^{2}$, using the fact that $P^{y} \big( 0 \in S[0, T] \big) \asymp \big( r 2^{n} \big)^{-1}$ and $P^{y} \big( z \in S[0, T] \big) \asymp \big( l 2^{n} \big)^{-1}$ (see \eqref{srwbound-2} for these estimates on hitting probabilities), it follows that 
\begin{align*}
\sum_{y \in \partial_i B_i} \sum_{z \in \partial_i B_j} P^{y} \big( 0 \in S[0, T] \big) P^{y} \big( z \in S[0, T] \big) P (F_{y, z} ) &\ge c \big( \epsilon 2^{n} \big)^{4} \big( r 2^{n} \big)^{-1} \big( l 2^{n} \big)^{-1} \big( \epsilon 2^{n} \big)^{-2} \Big( \frac{\epsilon}{l} \Big)^{2-\beta} \Big( \frac{\epsilon}{r} \Big)^{2-\beta} \\
&= c \Big( \frac{\epsilon}{r} \Big)^{3-\beta } \Big( \frac{\epsilon}{l} \Big)^{3-\beta },
\end{align*}
which completes the proof.\end{proof}

As a quick corollary of Proposition \ref{prop:bound-box-hit}, the next proposition gives an upper bound on the second moment of $Y$, whose proof we omit for brevity.

\begin{prop}\label{2nd-mom-mom}
There exists a universal constant $C > 0$ such that for all $k \ge 1$, $n \ge N_{k}$ and $B\subset\frac{2}{3}\mathbb{D}$ satisfying ${\rm dist} \big( B, \{0 \} \big) > 2^{-k}$, it follows that  
\begin{equation}\label{2nd-mom}
E [ Y^{2} ] \le C \epsilon^{- 2 \beta} 2^{- (3 + \beta ) k} r^{\beta- 3 },
\end{equation}
where $\beta$ is the growth exponent.
If $B \subset  \mathbb{D} \setminus \frac{1}{3} \mathbb{D}$, then 
\begin{equation}
E[Y^2] \leq C \epsilon^{- 2 \beta} 2^{- (3 + \beta ) k} r^{\beta- 1 }.
\end{equation}
\end{prop}

\subsection{Stability of one-point function}\label{sec:stability}
In this subsection, we show that one-point probability is stable under small perturbations. As a corollary, this also implies the continuity of the one-point function $c(x)$ from \eqref{one-point-def} on $\mathbb{D}\setminus\{0\}$. 

\begin{lem}\label{lem:stablepertubation}
 {Pick $x_0\in \mathbb{D}\setminus\{0\}$. For $n,k>0$ such that $2^{-n}<2^{-k^4+1}< d_x$, one has}
\begin{equation}\label{TRANSINVA}
P \big( x \in \gamma_{n} \big) = \Big( 1 + O \big( 2^{-c k^{3}} \big) \Big) P \big( x_{0} \in \gamma_{n} \big),
\end{equation}
uniformly in $x \in B_0 \cap 2^{-n} \mathbb{Z}^{3}$ where $B_0=B_\infty(x_0,2^{-k^4})$.
\end{lem}
\begin{proof} We refer readers to Table \ref{symbols-2} at the end fo this proof for a list of symbols.

To start with, we let $X$ be the random walk on $2^{-n} \mathbb{Z}^{3}$ conditioned to hit the origin 
before exiting $\mathbb{D}$. We also consider $R$ the simple random walk on $2^{-n} \mathbb{Z}^{3}$ which is independent of $X$. Write $Q^{x,x}$ for their joint law.
Then by  \cite[Lemma 5.1]{Escape}, we see that 
\begin{equation}\label{atui}
P ( x \in \gamma_{n} ) = G_{\mathbb{D}} \big( 0, x \big) Q^{x, x} \Big( \text{LE} (X [0, \tau_{0} ] ) \cap R [1, T] = \emptyset \Big) 
\end{equation}
Also, we note that $\tau_{0}$ (resp. $T$) represents the first time that $X$ (resp. $R$) hits the origin (resp. exits from $\mathbb{D}$). As in \eqref{green-1}, we see that 
\begin{equation*}
G_{\mathbb{D}} \big( 0,  x \big) = G_{\mathbb{D}} \big( 0,  x_{0} \big) \big( 1 + O (2^{-k^{4}} ) \big),
\end{equation*}
uniformly in $x \in B_{0} \cap 2^{-n} \mathbb{Z}^{3}$.

Thus, writing 
\begin{equation}\label{0712-a-1}
h(x) = Q^{x, x} \big( \text{LE} (X [0, \tau_{0} ] ) \cap R [1, T] = \emptyset \big), 
\end{equation}
we need to show that the difference between $h(x) $ and $ h(x_{0})$ is smaller than $C \cdot h (x_{0} ) \, 2^{- c k^{3}}$ uniformly in 
$x \in B_{0} \cap 2^{-n} \mathbb{Z}^{3}$. For that purpose, let $R^{1}$ and $R^{2}$ be the independent simple random walks on  $2^{-n} \mathbb{Z}^{3}$ and denote their joint law by $P^{x,x}$. Then, by definition of $X$ and $R$,
\begin{equation}\label{0712-a-0}
h(x) = \frac{P^{x, x} \Big( \text{LE} (R^{1} [0, \tau^{1}_{0} ] ) \cap R^{2} [1, T^{2}] = \emptyset, \  \tau^{1}_{0} < T^{1} \Big)}{P^{x} ( \tau^{1}_{0} < T^{1} ) }
\end{equation}
where $\tau^{1}_{0}$ stands for the first time that $R^{1}$ hits the origin and $T^{i}$ is the first time that $R^{i}$ exits from $\mathbb{D}$. However, by  \eqref{srwbound-2}, it follows that 
\begin{equation*}
 P^{x} ( \tau^{1}_{0} < T^{1} ) = P^{x_{0}} ( \tau^{1}_{0} < T^{1} )  \big( 1 + O (2^{-k^{4}} ) \big),
\end{equation*}
uniformly in $x \in B_{0} \cap 2^{-n} \mathbb{Z}^{3}$. Therefore, setting
\begin{equation}\label{0712-a-2}
q (x) = P^{x, x} \Big( \text{LE} (R^{1} [0, \tau^{1}_{0} ] ) \cap R^{2} [1, T^{2}] = \emptyset, \  \tau^{1}_{0} < T^{1} \Big)
\end{equation}
for the numerator of the  fraction  in the right hand side of \eqref{0712-a-0}, in order to prove \eqref{TRANSINVA}, it suffices to show that 
\begin{align}\label{atui1}
q (x) = q (x_{0} )  \big( 1 + O (2^{-c k^{3}} ) \big),
\end{align}
uniformly in $x \in B_{0} \cap 2^{-n} \mathbb{Z}^{3}$. To show this, we will follow the same spirit of  \cite[Lemma 5.3]{Escape}. Loosely speaking, it is proved in \cite[Lemma 5.3]{Escape} that taking $r$ sufficiently small the shape of $ \text{LE} (R^{1} [0, \tau^{1}_{0} ] ) \cap B_\infty (r) $ is not important to compute $q (x)$. To be more precise, we let $T^{i}_{r}$ be the first time that $R^{i}$ hits $\partial B_\infty (r)$. It is then proved in \cite[Lemma 5.3]{Escape} that 
\begin{align}\label{atui2}
q (x) =P^{x, x} \Big( \text{LE} (R^{1} [0, T^{1}_{2^{-k^{3}}} ] ) \cap R^{2} [1, T^{2}] = \emptyset, \  \tau^{1}_{0} < T^{1} \Big)  \big( 1 + O (2^{-c k^{3}} ) \big), \ \ \text{uniformly in $x \in B_{0}$,}
\end{align}
for some universal constant $c > 0$. Namely, in the sense of \eqref{atui2}, the loop-erasure of the end part of $R^{1}$ is not important to compute $q (x)$ as defined in  \eqref{0712-a-2} (see the proof of  \cite[Lemma 5.3]{Escape} for the details). 

On the other hand, it follows from \eqref{srwbound-2} that for all $z \in \partial B ( 2^{- k^{3}} )$
\begin{equation}\label{ctilde}
P^{z} (\tau^{1}_{0} < T^{1} ) = \widetilde{c} \, 2^{k^{3}} 2^{-n} \big( 1 + O (2^{- k^{3}} ) \big),
\end{equation}
for some universal constant $\widetilde{c} > 0$. Using this and the strong Markov property, we have 
\begin{align}\label{0712-b}
&P^{x, x} \Big( \text{LE} (R^{1} [0, T^{1}_{2^{-k^{3}}} ] ) \cap R^{2} [1, T^{2}] = \emptyset, \  \tau^{1}_{0} < T^{1} \Big) \notag  \\
&= E^{x, x} \Big[ {\bf 1} \Big\{ T^{1}_{2^{-k^{3}}} < T^{1}, \  \text{LE} (R^{1} [0, T^{1}_{2^{-k^{3}}} ] ) \cap R^{2} [1, T^{2}] = \emptyset \Big\} \cdot P^{R^{1} \big( T^{1}_{2^{-k^{3}}} \big)} ( \tau^{1}_{0} < T^{1} ) \Big] \notag  \\
&= \widetilde{c} \, 2^{k^{3}} 2^{-n} P^{x, x } \Big( T^{1}_{2^{-k^{3}}} < T^{1}, \  \text{LE} (R^{1} [0, T^{1}_{2^{-k^{3}}} ] ) \cap R^{2} [1, T^{2}] = \emptyset  \Big) \big( 1 + O (2^{- k^{3}} ) \big),
\end{align}
uniformly in $x \in B_{0} \cap 2^{-n} \mathbb{Z}^{3}$. With this in mind, we set 
\begin{equation}\label{0712-b-1}
p (x) = P^{x, x } \Big( T^{1}_{2^{-k^{3}}} < T^{1}, \  \text{LE} (R^{1} [0, T^{1}_{2^{-k^{3}}} ] ) \cap R^{2} [1, T^{2}] = \emptyset  \Big)
\end{equation}
for the probability appeared in the right hand side of \eqref{0712-b}.

Take $x \in B_{0} \cap 2^{-n} \mathbb{Z}^{3}$. We want to compare $p (x) $ and $p (x_{0})$ where $p ( \cdot )$ is as defined in \eqref{0712-b-1}.
This will be done by a technique using the translation invariance of the LERW (we note that a similar technique was used in the proof of Proposition \ref{PROP2}). 
Let $y = x -x_{0}$. Note that $|y| \le 3 \cdot 2^{-k^{4}}$. Define 
\begin{align}\label{0713-a-1}
 \widehat{B} (2^{-k^{3}}) = y+ B_\infty (2^{-k^{3}}) 
 \ \text{ and } \ \  \  \widehat{\mathbb{D}} = y+ \mathbb{D}.
\end{align}
We write 
\begin{itemize}
\item $\widehat{T}^{1}_{2^{-k^{3}}}$ for the first time that $R^{1}$ hits $\partial  \widehat{B} (2^{-k^{3}})$;

\item $\widehat{T}^{i}$ for the first time that $R^{i}$ exits from $\widehat{\mathbb{D}}$.
\end{itemize}
Also, we set 
\begin{equation}\label{0713-e-1}
\widehat{\gamma} = \text{LE} (R^{1} [0, \widehat{T}^{1}_{2^{-k^{3}}} ] ).
\end{equation}

By the translation invariance, we see that 
\begin{align}\label{0712-u}
p (x_{0}) = P^{x, x } \Big( \widehat{T}^{1}_{2^{-k^{3}}} < \widehat{T}^{1}, \ \widehat{\gamma}  \cap R^{2} [1, \widehat{T}^{2}] = \emptyset  \Big),
\end{align}
where $p (\cdot )$ is as defined in \eqref{0712-b-1}. We want to compare the right hand side of \eqref{0712-u} with $p (x)$.
For that purpose, we set 
\begin{equation}\label{0713}
\widehat{p} (x) =   P^{x, x } \Big( \widehat{T}^{1}_{2^{-k^{3}}} < \widehat{T}^{1}, \ \widehat{\gamma}  \cap R^{2} [1, \widehat{T}^{2}] = \emptyset  \Big)
\end{equation}
for the right hand side of \eqref{0712-u}. We also define
\begin{equation}\label{0713-1}
\widehat{F} = \big\{ \widehat{T}^{1}_{2^{-k^{3}}} < \widehat{T}^{1} \big\}  \ \ \text{ and }  \ \ \widehat{H} = \big\{ \widehat{\gamma}  \cap R^{2} [1, \widehat{T}^{2}] = \emptyset \big\}
\end{equation}
so that $\widehat{p} (x) = P^{x, x }  ( \widehat{F} \cap \widehat{H} )$.

To compare $p (x) $ and $\widehat{p} (x)$, 
we need to show that the diameter of $R^{2}$ from $\widehat{T}^{2} \wedge T^{2}$ to $\widehat{T}^{2} \vee T^{2}$ is 
small. To be more precise, we write
\begin{equation}\label{0713-2}
G^{2} = \Big\{ \text{diam} \Big( R^{2} \big[  \widehat{T}^{2} \wedge T^{2},  \widehat{T}^{2} \vee T^{2} \big] \Big) \le 2^{- \frac{k^{4}}{2}} \Big\}.
\end{equation}
We define $\widehat{u}^{1}$ and $u^{2}$ by 
\begin{equation}\label{0713-3}
\widehat{u}^{1} = \inf \Big\{ j \ge 0 \ \Big| \  \widehat{\gamma}  (j) \notin B_\infty (x, 10^{-1} ) \Big\} \ \ \text{ and } \ \  u^{2} =  \inf \Big\{ j \ge 0 \ \Big| \  R^{2} (j) \notin B_\infty (x, 10^{-1} ) \Big\}.
\end{equation}

We are now ready to deal with the event $G^{2}$. Using the strong Markov property for $R^{2}$, we have that 
\begin{align}\label{0713-4}
P^{x, x }  \big( \widehat{H} \cap (G^{2} )^{c} \ \big| \ \widehat{F} \big) \le P^{x, x }  \Big( \widehat{\gamma}  [0, \widehat{u}^{1} ] \cap R^{2} [1, u^{2} ] = \emptyset   \ \Big|  \  \widehat{F} \Big) \cdot  \max_{z \in \partial B (x, \frac{1}{10} )}  P^{z} \big( (G^{2})^{c} \big). 
\end{align}
However, since $\text{dist} \big( \partial \mathbb{D},  \partial \widehat{\mathbb{D}} \big) \le |y| \le 3 \cdot 2^{-k^{4}}$ by definition as in \eqref{0713-a-1}, it follows from the gambler's ruin estimate (see \cite[Proposition 5.1.6]{LawLim} for this) that $P^{z} \big( (G^{2})^{c} \big)$ is bounded above by $C 2^{- \frac{k^{4}}{2}}$ uniformly in $z \in \partial B_\infty (x, \frac{1}{10} )$. Furthermore, \cite[Proposition 4.4]{Mas} guarantees that conditioning on the event $\widehat{F}$, the conditional distribution of $\widehat{\gamma}  [0, \widehat{u}^{1} ]$ is approximated by that of the infinite LERW stated at $x$ until the first time that it exits from $B_\infty (x, 10^{-1} )$. More precisely, if we write $\gamma'$ for the infinite LERW started at $x$ and set $u'$ for the first time that $\gamma'$ exits from  $B_\infty (x, 10^{-1} )$, then there exists a universal constant $c > 0$ such that 
\begin{equation}\label{0714-a-1}
c \, P \big( \gamma' [0, u'] = \lambda \big)     \le  P \big( \widehat{\gamma}  [0, \widehat{u}^{1} ] = \lambda \ \big| \ \widehat{F} \big) \le c^{-1} \, P \big( \gamma' [0, u'] = \lambda \big) 
\end{equation}
for each path $\lambda$. Thus, the right hand side of \eqref{0713-4} is smaller than
\begin{equation}\label{0714-a-2}
C 2^{- \frac{k^{4}}{2}} \, P^{x,x} \big( \gamma' [0, u']  \cap R^{2} [1, u^{2} ] = \emptyset \big)  \asymp 2^{- \frac{k^{4}}{2}} \, \Es (2^{n}).
\end{equation}
Here we used \cite[Proposition 6.7]{S}  in the last equation above, (see Section \ref{sec:Es} for the definition and property of $\Es (\cdot )$). 
Combining this with the fact that $P^{x, x }  \big( \widehat{H}  \ \big| \ \widehat{F} \big)$ is comparable to $\Es (2^{n} ) $, we have that the left hand side of \eqref{0713-4} is bounded above by $C 2^{- \frac{k^{4}}{2}} P^{x, x }  \big( \widehat{H}  \ \big| \ \widehat{F} \big)$.

With this in mind, setting 
\begin{equation}\label{0713-b-1}
\widehat{H}' = \big\{ \widehat{\gamma}  \cap R^{2} [1, T^{2}] = \emptyset \big\},
\end{equation}
we want to replace $\widehat{H}$ by $\widehat{H}' $. Note that the difference between $\widehat{H}$ and
 $\widehat{H}' $ comes from that of $\widehat{T}^{2}$ and $T^{2}$, see \eqref{0713-1} and \eqref{0713-b-1} for this. Suppose 
 that the event $\widehat{H} \cap (\widehat{H}' )^{c} \cap G^{2}$ occurs. This implies that $\widehat{T}^{2} < T^{2}$ and 
 $\widehat{\gamma} $ hits $R^{2} [ \widehat{T}^{2}, T^{2} ]$; see Figure \ref{0713-f}. But the event $G^{2}$ ensures that the diameter of $R^{2} [ \widehat{T}^{2}, T^{2} ]$ is smaller than $2^{- \frac{k^{4}}{2}}$. Thus, writing
 \begin{equation}\label{0713-b-3}
 A = \Big\{ \widehat{\gamma}  [0, \widehat{u}^{1} ] \cap R^{2} [1, u^{2} ] = \emptyset \Big\}
 \end{equation}
 for the event appeared in the right hand side \eqref{0713-4} (see Figure \ref{0713-f} for this event), it follows from the domain Markov property for the LERW (see Section \ref{DMPLEW} for this) that 
 \begin{align}\label{0713-b-2}
 P^{x, x }  \big( \widehat{H} \cap (\widehat{H}' )^{c} \cap G^{2} \ \big| \ \widehat{F} \big) \le P^{x, x }  \big( A \ \big| \ \widehat{F} \big) \,  \max_{\lambda \in \Lambda, \ w \in \partial \widehat{\mathbb{D}}}  P \Big( Z [0,a] \cap B_\infty (w,   2^{- \frac{k^{4}}{2}} ) \neq \emptyset \Big).
 \end{align}
Here we set
\begin{itemize}
\item  $\Lambda = \Big\{ \lambda \ \Big| \ \lambda \text{ is a simple path in $2^{-n} \mathbb{Z}^{3}$ with } P^{x} \big( \widehat{\gamma}  [0, \widehat{u}^{1} ] = \lambda \ \big| \ \widehat{F} \big) > 0 \Big\},$

\item  $Z$ for the random walk on $2^{-n} \mathbb{Z}^{3}$ started at $\lambda ( \text{len} (\lambda ))$ and conditioned that $Z [1, a] \cap \big(\lambda \cup \partial \widehat{\mathbb{D}} \big) = \emptyset$ where $a$ is the first time that $Z$ hits $\widehat{B} (2^{- k^{3}} )$, see \eqref{0713-a-1} for $\widehat{\mathbb{D}}$ and $\widehat{B} (2^{- k^{3}} )$.
\end{itemize}

\begin{figure}[h]\label{supp-fig-3}
\begin{center}
\includegraphics[scale=0.55]{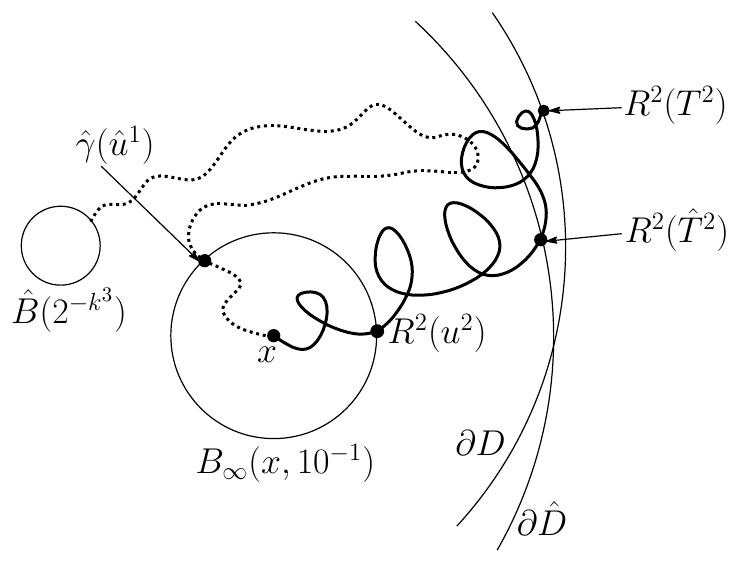}
\caption{Illustration for the events $\widehat{H} \cap (\widehat{H}' )^{c} \cap G^{2}$ and $A$. The dotted curve stands for $\widehat{\gamma}$ and the solid curve represents $R^{2}$. Also, we let $D = \mathbb{D}$ and $\hat{D} = \widehat{\mathbb{D}}$ in the picture.
}\label{0713-f}
\end{center}
\end{figure}

We want to show that the maximum in the right hand side of \eqref{0713-b-2} is smaller than $C 2^{- \frac{k^{4}}{2}}$. 
To see this, we fix $\lambda \in \Lambda$ and $w \in  \partial \widehat{\mathbb{D}}$. Let $z_{0} = \lambda ( \text{len} ( \lambda ) )$ be the end point of $\lambda$, which is the starting point of $Z$. Then we have that 
\begin{equation}\label{0714}
P \Big( Z [0,a] \cap B_\infty (w,   2^{- \frac{k^{4}}{2}} ) \neq \emptyset \Big) = \frac{P^{z_{0}} \Big( R^{1} [0, \widehat{T}^{1}_{2^{-k^{3}}}] \cap B (w,   2^{- \frac{k^{4}}{2}} ) \neq \emptyset, \ \widehat{F}, \  R^{1} [1, \widehat{T}^{1}_{2^{-k^{3}}}] \cap \lambda = \emptyset  \Big)}{P^{z_{0}} \Big( \widehat{F}, \  R^{1} [1, \widehat{T}^{1}_{2^{-k^{3}}}] \cap \lambda = \emptyset  \Big)},
\end{equation}
where $\widehat{F}$ is as defined in  \eqref{0713-1}. Writing $b$ for the first time that $R^{1}$ exits from $B_\infty (x, \frac{1}{5} )$, it follows from \eqref{srwbound} that the denominator of the fraction in the right hand side of \eqref{0714} is bounded below by $c \, 2^{- k^{3}}  P^{z_{0}} \big( R^{1} [1, b] \cap \lambda = \emptyset \big)$ for some universal constant $c > 0$. On the other hand, using \eqref{srwbound} again, the  numerator  of the fraction is bounded above by $C \, 2^{- k^{3}} \, 2^{- \frac{k^{4}}{2}} P^{z_{0}} \big( R^{1} [1, b] \cap \lambda = \emptyset \big)$ for some universal constant $C < \infty$. Thus, the maximum in the right hand side of \eqref{0713-b-2} is smaller than $C 2^{- \frac{k^{4}}{2}}$, as desired.

From this and the fact that $P^{x, x }  \big( A \ \big| \ \widehat{F} \big)$ is comparable to $\Es (2^{n} )$ (see \eqref{0714-a-1} and \eqref{0714-a-2} for this), we have that the left hand side of \eqref{0713-b-2} is bounded above by $C 2^{- \frac{k^{4}}{2}} P^{x, x }  \big( \widehat{H}  \ \big| \ \widehat{F} \big)$. Consequently, it follows that 
\begin{equation}\label{0713-b-3-2}
P^{x, x }  \big( \widehat{H} \cap (\widehat{H}' )^{c}  \ \big| \ \widehat{F} \big) \le C 2^{- \frac{k^{4}}{2}} P^{x, x }  \big( \widehat{H}  \ \big| \ \widehat{F} \big).
\end{equation}
 Changing the roles of $\widehat{H} $ and $\widehat{H}'$, we also have that $P^{x, x }  \big( \widehat{H}^{c} \cap \widehat{H}'  \ \big| \ \widehat{F} \big) $ is bounded above by $C 2^{- \frac{k^{4}}{2}} P^{x, x }  \big( \widehat{H}  \ \big| \ \widehat{F} \big)$. Therefore, we can replace $\widehat{H}$ by $\widehat{H}'$, i.e., it follows that 
 \begin{equation}\label{0713-c-1}
 \widehat{p} (x) = P^{x, x }  \big( \widehat{F} \cap  \widehat{H} \big) = P^{x, x }  \big( \widehat{F} \cap  \widehat{H}' \big) \,  \Big[  1 + O \big( 2^{- \frac{k^{4}}{2}} \big) \Big].
 \end{equation}

In order to replace $\widehat{T}^{1}$ and $\widehat{T}^{1}_{2^{-k^{3}}}$  with $T^{1}$ and $T^{1}_{2^{-k^{3}}}$ respectively, we can follow the completely same arguments as above that we used to replace $\widehat{T}^{2}$ by $T^{2}$. Modifying the arguments above, one can prove that 
\begin{equation}\label{0713-d-1}
\widehat{p} (x) = p (x) \Big[  1 + O \big( 2^{- \frac{k^{4}}{2}} \big) \Big],
\end{equation}
 where $p (x)$ and $\widehat{p} (x) $ are as defined in \eqref{0712-b-1} and \eqref{0713}. Easy modifications to get \eqref{0713-d-1} are left to the reader.

Consequently, it follows that 
\begin{align*}
q (x_{0}) &\quad\;\stackrel{\eqref{atui2}}{=}\quad P^{x_{0}, x_{0}} \Big( \text{LE} (R^{1} [0, T^{1}_{2^{-k^{3}}} ] ) \cap R^{2} [1, T^{2}] = \emptyset, \  \tau^{1}_{0} < T^{1} \Big) \,  \big[  1 + O (2^{-c k^{3}} ) \big] \\
& \stackrel{\eqref{0712-b},  \eqref{0712-b-1}}{=}  \widetilde{c} \, 2^{k^{3}} 2^{-n} p(x_{0} )   \,  \big[  1 + O (2^{-c k^{3}} ) \big] \\
& \stackrel{\eqref{0712-u},  \eqref{0713}}{=}   \widetilde{c} \, 2^{k^{3}} 2^{-n}     \widehat{p} (x)   \,  \big[  1 + O (2^{-c k^{3}} ) \big] \\ 
&\quad\;\stackrel{\eqref{0713-d-1}}{=} \quad \widetilde{c} \, 2^{k^{3}} 2^{-n} p (x)  \,  \big[  1 + O (2^{-c k^{3}} ) \big]   \\ 
&\stackrel{\eqref{0712-b},  \eqref{0712-b-1}}{=}  P^{x, x} \Big( \text{LE}  (R^{1} [0, \tau^{1}_{0} ] ) \cap R^{2} [1, T^{2}] = \emptyset, \  \tau^{1}_{0} < T^{1} \Big)   \,  \big[  1 + O (2^{-c k^{3}} ) \big]  \\
& \quad\,\stackrel{\eqref{atui2}}{=} \quad\,q(x)   \,  \big[  1 + O (2^{-c k^{3}} ) \big],
\end{align*}
which gives \eqref{atui1} and hence \eqref{TRANSINVA}.\end{proof}

\begin{table}[hbtp]
  \centering
  \begin{tabular}{|c|c||c|c|}
    \hline \hline
    Symbol  & Meaning   &    Symbol  & Meaning  \\
    \hline \hline
    $T^{i}$ & First time that $R^{i}$ exits from $\mathbb{D}$ &
    $T^{i}_{r}$  & First time that $R^{i}$ hits $\partial B (r)$   \\ \hline
    $p (x) $  &    Probability def'd in \eqref{0712-b-1}&
    $y$ & $x - x_{0}$ \\ \hline
    $\widehat{B} (r)$ & $ y+ B_\infty (r) $ &
    $\widehat{\mathbb{D}}$ & $y+ \mathbb{D} $ \\ \hline
    $\widehat{T}^{i}$ & First time that $R^{i}$ exits from $\widehat{\mathbb{D}}$ &
    $\widehat{T}^{i}_{r}$  & First time that $R^{i}$ hits $\partial \widehat{B} (r)$   \\ \hline
    $\widehat{\gamma}$ &  $\text{LE} (R^{1} [0, \widehat{T}^{1}_{2^{-k^{3}}} ] )$ &
    $\widehat{p} (x) $  & Probability def'd in \eqref{0713} \\ \hline
    $\widehat{F}$ & Event $ \big\{ \widehat{T}^{1}_{2^{-k^{3}}} < \widehat{T}^{1} \big\} $ &
    $\widehat{H} $ &  Event $ \big\{ \widehat{\gamma}  \cap R^{2} [1, \widehat{T}^{2}] = \emptyset \big\}$ \\ \hline
    $\widehat{H}' $ & Event $ \big\{ \widehat{\gamma}  \cap R^{2} [1, T^{2}] = \emptyset \big\} $ & & \\ 
   \hline \hline
  \end{tabular}
  \caption{List of symbols used in the proof of Lemma \ref{lem:stablepertubation}}
  \label{symbols-2}
\end{table}

\section{Proof of the key $L^2$-estimate \eqref{KEY}}\label{sec:4}
In this section we focus on the proof of the key two-point estimate \eqref{KEY}, which we restate as Proposition \ref{prop4.6} below, through which we obtain an $L^2$-approximation of the renormalized occupation measure of the LERW by counting the number of mesoscopic cubes the LERW crosses, which in turn can be well-approximated by the corresponding quantity for the scaling limit. 

\subsection{Statement, preliminary reductions and comments}\label{sec:sketch}
We first recall notation introduced at the beginning of Section \ref{sec:3}, in particular the definition of the cube $B$ and its partition into smaller cubes $B_i$'s in Definition \ref{def:Cubes}. Let $x_{0}$ stand for the discretization of $(\frac{1}{2}, 0, 0 )$ in $2^{-n}\mathbb{Z}^3$ (i.e., any closest point on $2^{-n}\mathbb{Z}^3$). Note that this corresponds to the choice of $\widehat{x}_n$ in the definition of the scaling factor $f_n$ in \eqref{eq:fndef}. We then consider $B_{0}$, a ``reference'' closed cube of side length $2^{-k^{4}+1}$ (which is of the same size as $B_i$'s in Definition \ref{def:Cubes}) centered at $x_0$. We set $X_{0}$ for the number of points in $B_{0} \cap 2^{-n} \mathbb{Z}^{3}$ which is passed through by $\gamma$. Let $Y_{0}$ be the indicator function of the event that $\gamma$ hits $B_{0}$. Let 
\begin{equation}\label{Beta}
\alpha_0=\alpha_0(n,k) := E \Big( X_{0} \ \Big| \ Y_{0} = 1 \Big).
\end{equation}
Using \eqref{20221118} and \eqref{20221118-1}, we see that 
\begin{equation}
\alpha_{0} \asymp (\epsilon 2^{n} )^{\beta}.
\end{equation}
Recall that $X$ and $Y$ are (roughly speaking) the total number of vertices and cubes passed by $\gamma$ defined  in \eqref{X} and \eqref{Y} respectively.

We are now ready to state the key estimate  \eqref{KEY} in its precise form.
\begin{prop}\label{prop4.6}
There exist universal constants $c > 0$, $C < \infty$ and a constant $N_{k}$ depending only on $k$ such that for all $k \ge 1$, $n \ge N_{k}$ and any cube $B$ satisfying the requirements in Definition \ref{def:Cubes}, we have
\begin{equation}\label{eq}
E \Big( \big( X - \alpha_0 Y \big)^{2} \Big) \le C 2^{- c k^{2} } \big[ E (X) \big]^{2}.
\end{equation}
\end{prop}
The crux of this proposition is the following estimates that relate the conditional expectation of $X_i$'s to $\alpha_0$, which we summarize in the following proposition. 

We start with notation. To distinguish the typical case (where $i$ and $j$ are distant from each other) from the atypical case  (whose contribution towards the RHS of \eqref{eq1} below is negligible), we say that
\begin{equation}\label{goodbad}
\mbox{the pair $(i, j)$ is {\bf good}, if $l:=|x_i-x_j| \ge 2^{- k^{2}}$, and {\bf bad} otherwise}.
\end{equation}
As a convention, for two quantities $p$ and $q$ which depend on $k$ (and possibly also on $n\gg k$), we write $p \simeq q$ if 
 \begin{equation}\label{eq:simeqconv}
 p = \big( 1 + O (2^{-k^{2}} ) \big) q
 \end{equation} uniformly for all $n\geq N_k$ where $N_k$ depends on $k$.

\begin{prop}\label{prop:goodpair} Under the same setup as in Proposition \ref{prop4.6}, for a good pair $(i, j)$, such that for any $v\in B_i$, $w\in B_j$,
\begin{equation}\label{eq:MAIN}
|B_{0}|^{2} P \Big( 0 \xrightarrow{\gamma} v \xrightarrow{\gamma} w  \Big) \simeq \alpha_0^{2} P \Big( 0 \xrightarrow{\gamma} B_{i} \xrightarrow{\gamma} B_{j} \Big) 
\end{equation}
In a similar fashion,
\begin{equation}\label{eq:MAINprime}
|B_{0}| P \Big( 0 \xrightarrow{\gamma} x_{i} \xrightarrow{\gamma} B_{j}  \Big) \simeq \alpha_0 P \Big( 0 \xrightarrow{\gamma} B_{i} \xrightarrow{\gamma} B_{j} \Big)
\end{equation}
and the same is true if the LHS  of \eqref{eq:MAINprime} above is replaced by $|B_{0}| P \Big( 0 \xrightarrow{\gamma} B_{i} \xrightarrow{\gamma} x_{j}  \Big)$.
\end{prop}

Because the proof of this proposition is quite long (in fact almost the entire section is dedicated to it), we will postpone it to Section \ref{sec:4.6}. A short explanation on the proof and the organization of the rest of this section is given at the end of this subsection.

As a quick corollary of Proposition \ref{prop:goodpair} we obtain the following conditional expectations of $X_i$'s, which feeds into the proof of Proposition \ref{prop4.6}.
\begin{cor}\label{cor:goodpair}
Under the same setup as in Proposition \ref{prop:goodpair},
\begin{align}
E \Big( X_{i} X_{j} \ \Big| \ Y_{i} = Y_{j} = 1 \Big) &= \Big( 1 + O \big(  2^{-k^{2} } \big) \Big) \alpha_0^{2}; \label{eq3} \\
E \Big( X_{i}  \ \Big| \ Y_{i} = Y_{j} = 1 \Big) &= \Big( 1 + O \big(  2^{-k^{2} } \big) \Big) \alpha_0, \label{eq4} 
\end{align}
\end{cor}
%

We now turn back to Proposition \ref{prop4.6}. 
\begin{proof}[Proof of Proposition \ref{prop4.6}]
We rewrite the LHS of \eqref{eq} as 
\begin{equation}\label{eq1}
E \Big( \big( X - \alpha_0 Y \big)^{2} \Big) = \sum_{1 \le i, j \le M_k} E (Y_{i} Y_{j} ) E \Big( \big( X_{i} - \alpha_0 \big) \big( X_{j} - \alpha_0 \big)  \ \Big| \ Y_{i} = Y_{j} = 1 \Big).
\end{equation} 
Corollary \ref{cor:goodpair}  implies that for good $(i,j)$'s,
\begin{equation}\label{rewrite1}
E \Big( \big( X_{i} - \alpha_0 \big) \big( X_{j} - \alpha_0 \big)  \ \Big| \ Y_{i} = Y_{j} = 1 \Big) = \alpha_0^{2} O \big( 2^{-k^{2}} \big). 
\end{equation}
For bad $(i, j)$'s, by Proposition \ref{prop:bound-box-hit} (in particular \eqref{M7} and \eqref{M7primeprime}), it follows that
\begin{equation}\label{eq6-2}
E \Big( X_{i} X_{j} \ \Big| \ Y_{i} = Y_{j} = 1 \Big) \le C \alpha_0^{2};\quad E \Big( X_{i}  \ \Big| \ Y_{i} = Y_{j} = 1 \Big), E \Big( X_{j}  \ \Big| \ Y_{i} = Y_{j} = 1 \Big) \le C \alpha_0
\end{equation}
for some universal constant $C < \infty$. Therefore,
\begin{equation}
E \Big( \big( X - \alpha_0 Y \big)^{2} \Big) 
\le C 2^{-k^{2}} \sum_{(i, j) \text{ good}} E (Y_{i} Y_{j} )  \alpha_0^{2}  + C \sum_{(i, j) \text{ bad}} E (Y_{i} Y_{j} ) \alpha_0^{2}.  \label{eq9}
\end{equation}
We first bound $E(Y_iY_j)$ for bad $(i,j)$'s. As a direct consequence of Proposition \ref{prop:bound-box-hit}, we have 
\begin{equation}\label{eq10}
\sum_{(i, j) \text{ bad}} E (Y_{i} Y_{j} )  \le C 2^{-c k^{2}} \sum_{(i, j) \text{ good}} E (Y_{i} Y_{j} ),
\end{equation}
for some universal constants $0 < c, C < \infty$. Thus, we have
\begin{equation}\label{eq11}
E \Big( \big( X - \alpha_0 Y \big)^{2} \Big) \le C 2^{-c k^{2}}  \sum_{(i, j) \text{ good}} E (Y_{i} Y_{j} )  \alpha_0^{2}.
\end{equation}
Moreover, by \eqref{eq3}, we see that
\begin{equation}\label{eq12}
\sum_{(i, j) \text{ good}} E (Y_{i} Y_{j} )  \alpha_0^{2} \le C \sum_{(i, j) \text{ good}} E (X_{i} X_{j} ) \le C  \sum_{1 \le i, j \le M_k}  E (X_{i} X_{j} ) 
\end{equation}
Consequently, we have 
\begin{equation}\label{eq13}
E \Big( \big( X - \alpha_0 Y \big)^{2} \Big) \le C 2^{-c k^{2}} \sum_{1 \le i, j \le M_k}  E (X_{i} X_{j} ) = C 2^{-c k^{2}} E\big( X^{2} \big).
\end{equation}
When $B \subset \frac{2}{3}\mathbb{D}$, it follows that
\begin{equation} \label{eq15}
{E (X) \asymp 2^{-3k} 2^{\beta n} r^{\beta-3} \ge c2^{-3k} 2^{\beta n}} \quad\mbox{ and }\quad E\big( X^{2} \big) \asymp E(X) \Big( 2^{-k} \cdot 2^{n} \Big)^{\beta}.
\end{equation}
Similar bounds hold when $B\subset \mathbb{D}\setminus\frac{1}{3}\mathbb{D}$.
 Thus, we have
\begin{equation}\label{eq17}
E \big( X^{2} \big) \le C 2^{(3-\beta ) k} \Big( E(X) \Big)^{2},
\end{equation}
which gives
\begin{equation}\label{eq18}
E \Big( \big( X - \alpha_0 Y \big)^{2} \Big) \le C 2^{-c k^{2} + (3-\beta ) k } \Big( E(X) \Big)^{2} \le C 2^{-c' k^{2} } \Big( E(X) \Big)^{2}.
\end{equation}
This gives \eqref{eq} and finishes the proof.\end{proof}

We now briefly comment on the proof of Proposition \ref{prop:goodpair} and the organization of this section.

The heart of the matter is the analysis of the cube-crossing event on the RHS of \eqref{eq:MAIN}. In Section \ref{sec:decomp}, we rewrite this event through the last-exit decomposition of the underlying simple random walk, and  in Proposition \ref{prop:decomp} we express (in an asymptotically equivalent sense) the probability of this event, conditioned on the location of the last exits from these cubes, into that of an event (see the definition of ${\cal F}$ in \eqref{eq:Hdef}) involving the non-intersection in a certain manner of (truncated versions of) three independent random walks (with some conditioning) and their loop erasures. We refer readers to the end of Section \ref{sec:decomp} for more discussion on this decomposition.

We then concentrate on the non-intersection event $\cal F$ in Section \ref{sec:decoupling}. We claim in Proposition \ref{prop:Lem11} (whose proof we postpone till Section \ref{sec:endingpart}) that (again in an asymptotically equivalent sense) in expressing the probability of $\cal F$ by conditioning on the configuration of paths near both cubes, one can replace the law of these conditioned configurations by that of the ``end parts'' only, sampled from a reference measure which do not depend on the choice of $i$ and $j$ in Proposition \ref{prop:goodpair} above. We refer readers to the beginning of Section \ref{sec:decoupling} for more discussions.

After treating the counterpart point-crossing event (i.e., the LHS of \eqref{eq:MAIN}) and obtaining a similar decomposition in Section \ref{sec:pt}, we proceed with the proof of Proposition \ref{prop:goodpair} with the help of couplings introduced in Section \ref{sec:basicproperties}. The fact that the conditional probability largely only depends on the ``end parts'' of paths near the cubes, which we have established in Section \ref{sec:decoupling}, plays a crucial role here.

In Section \ref{sec:4.7}, we lay out the proof of Proposition \ref{prop:decomp} postponed from Section \ref{sec:decomp}. It consists of multiple steps of fine tuning the definition of cube-crossing event through the analysis on the atypical behavior or LERW, such as forming ``quasi-loops'' at certain segments. 

Finally in Section \ref{sec:endingpart}, we give the proof of Proposition \ref{prop:Lem11} of Section \ref{sec:decoupling}, where we first decouple the events near both cubes and then show that it is indeed legitimate to look at the ``end parts'' of the configurations near both cubes, as in the discussion above for Section \ref{sec:decoupling}  More detailed explanations can be found at the beginning of  Section \ref{sec:endingpart}.

\subsection{Path decomposition for the  cube-crossing event} \label{sec:decomp}
As discussed at the end of last subsection, in this subsection, we are going to lay down the first step of the proof of Proposition \ref{prop:goodpair} by rewriting the cube-crossing event on the RHS of \eqref{eq:MAIN} into a form which is much easier to analyze. In Proposition \ref{prop:decomp}, which is the crux of this subsection, we will show that for a pair of typical cubes the probability that they are crossed consecutively can be approximated by Green's functions and the probability of  some non-intersection event of three independent random walks (and loop erasures thereof); see the paragraph after  Proposition \ref{prop:decomp} for a brief explanation on the decomposition. The proof is postponed to Section \ref{sec:4.7} as it involves multiple reductions that are both long and technical.

Pick $i,j\in \{1,2,\ldots,M_k\}$  and keep them fixed throughout this subsection. 
From now on until the end of Section \ref{sec:pt}, we assume $(i, j)$ is good in the sense of \eqref{goodbad}.

Pick $y \in \partial B_{i}$ and $z \in \partial B_{j}$ and consider three independent stopped random walks:
\begin{itemize}
\item  $X^{1}$:  the random walk started at $y$ and conditioned to hit the origin (and then stopped there) before leaving $\mathbb{D}$;
\item $X^{2}$: the walk started at $z$ and conditioned to hit $y$ (and then stopped there) before leaving $\mathbb{D}$; 
\item $X^{3}$: the simple random walk started at $z$ stopped upon leaving $\mathbb{D}$.
\end{itemize}
Write $\mP^{y,z,z}$ for their joint law.

From now on till the end of this section, we use the following notation and conventions. For a path $\lambda$, we write
\begin{equation}\label{eq:Udef}
U_{i,-m}=U_{i,-m}(\lambda):=\min\big\{t\geq 0 \mid \lambda(t)\in \partial B_\infty \big( x_{i}, 2^{-m} \big)\big\}
\end{equation} 
and define $U_{j,-m}$ 
in the same fashion. Recall that if in addition $\lambda$ is finite, we write 
\begin{equation}\label{eq:Diamond}
\Diamond=\Diamond(\lambda)
\end{equation}
for its ending time. As a convention, we always omit the dependence on the specific path in notation whenever there is no confusion. 

We then consider 
\begin{equation}\label{eq:gammadef}
\gamma^1=\LE(X^1)\mbox{ and }\gamma^{2}=\LE(X^2[0,\sigma])
\end{equation} where we write
\begin{equation}\label{eq:sigmadef}
\sigma:=\max\big\{t \mid X^2(t) \in \partial B_\infty(x_{i},2^{-k^3})\big\}
\end{equation}
and let  
\begin{equation}\label{eq:V1V2def}
 V_1= \max  \big\{ t \le U_{i, -k^4+k^3}(\gamma^1)\ \Big| \ \gamma^{1} (t) \in B_{i} \big\}\mbox{, and } V_2= \max  \big\{ t \le U_{j, -k^4+k^3}(\gamma^2) \ \Big| \ \gamma^{2} (t) \in B_{j} \big\}.
\end{equation}
Note that in the definition above, the reason for requiring $V_1 \leq U_{i, -k^4+k^3}$ and $V_2 \leq U_{j, -k^4+k^3}$ is that we want to let the initial configurations chosen in the coupling argument in \eqref{eq:couplinginit} in the proof of Proposition \ref{prop:goodpair} in Section \ref{sec:4.6} satisfy the requirement of the results on couplings recalled in Section \ref{sec:basicproperties}.

We are now ready to state the decomposition proposition. As mentioned at the beginning of this subsection, the proof will be postponed to the end of Section \ref{sec:4.7}. 
\begin{prop}\label{prop:decomp}
For good pairs $(i,j)$, one has (recall the end of Section \ref{sec-green} for convention on the use of $G_{\mathbb{D}} := G_{\mathbb{D} \cap 2^{-n} \mathbb{Z}^{3}} $)
\begin{equation}\label{eq:decomp}
P\Big({0\overset{\gamma}{\rightarrow}B_i \overset{\gamma}{\rightarrow} B_j}\Big) 
\simeq\sum_{y\in\partial_i B_i} \sum_{z\in\partial_i B_j} G_{\mathbb{D}} (0,y) G_{\mathbb{D}} (y,z) \mP^{y,z,z}\left({\cal F}
\right),
\end{equation}
where the event ${\cal F}$ is defined by (recall that $\Diamond$ stands for the ending time of a path) 
\begin{equation}\label{eq:Hdef}
{\cal F}={\cal F}(y,z):=\left\{\begin{aligned}
&\gamma^1[V_1,\Diamond] \cap \big((X^2)^{\cal R}[1, U_{j,-k^3+k^2}]\cup X^3\big) =\emptyset; \; (X^2)^{\cal R}[1, U_{j,-k^3+k^2}]\cap B_i=\emptyset; \\
& \gamma^2[V_2,\Diamond
] \cap X^3[1,\Diamond]=\emptyset; \qquad\qquad\qquad\quad\quad\;\;
 X^3[1,\Diamond]\cap (B_i \cup B_j) =\emptyset
\end{aligned}\right\}.
\end{equation}
\end{prop}

\begin{figure}[h]
\begin{center}
\includegraphics[scale=0.6]{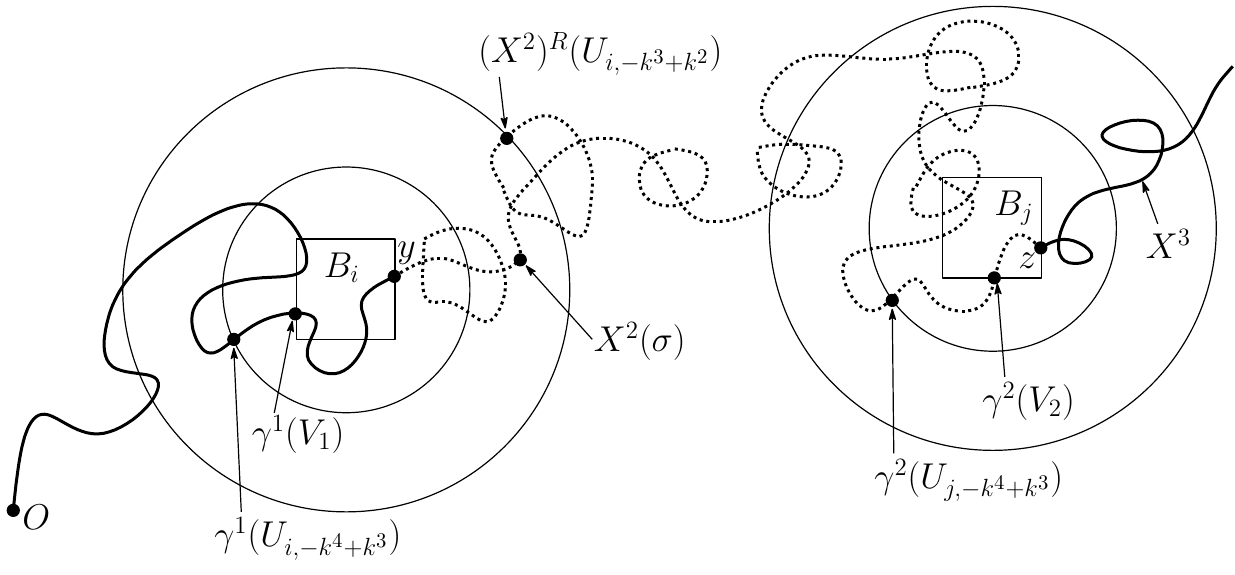}
\caption{Illustration for  the event ${\cal F}$ in \eqref{eq:Hdef}. The solid curves starting from $y$ and $z$ stand for $\gamma^{1}$ and $X^{3}$ respectively. The dotted curve from $z$ to $y$ represents $X^{2}$, while $\gamma^{2}$ is the loop erasure of $X^{2} [0, \sigma]$ defined as in \eqref{eq:gammadef}.}\label{fig28}
\end{center}
\end{figure}

See Figure \ref{fig28} for the event ${\cal F}$. We now briefly explain the moral of this decomposition. When a random walk passes cubes $B_i$ and $B_j$ in order, we can partition it into three parts at last exits of $B_i$ and $B_j$ that corresponds (in distribution) to $(X^1)^{\cal R}$, $(X^2)^{\cal R}$ and $X^3$ defined above with the requirement that $X^2$ does not enter $B_i$ and $X^3$ stays away from $B_i \cup B_j$ (except at the end and beginning respectively). In addition, to characterize the event that its loop erasure also passes these two cubes, one needs to pose extra constraints which roughly state as (in the notation above)  $\text{``$\gamma^1$''}\cap (X^2 \cup X^3)=\text{``$\gamma^2$''}\cap X^3=\emptyset$ (here ``$\gamma^\cdot$'' means a properly truncated version of $\gamma^\cdot$). However, in order to create some independence to facilitate the application of coupling results when we analyze the probability of non-intersection, in the actual definition of $\cal F$, we need to further truncate $X^2$ and $X^3$ as in \eqref{eq:Hdef} (also note that by definition $\gamma^2$ is already the loop-erasure of a truncated $X^2$). As we are going to see in Section \ref{sec:4.7}, these nuances in the definition  incur an error of order $O(2^{-ck^2})$, whence the $\simeq$ symbol in \eqref{eq:decomp}.

\subsection{Decoupling of paths }\label{sec:decoupling}
With the decomposition in Proposition \ref{prop:decomp} in the last subsection in mind, the main goal of this subsection is Proposition  \ref{prop:boxdecomp}  in which we rewrite (in an approximate sense) the cube-crossing probability into a ``decoupled'' form that facilitates the comparison with the corresponding rewritten form of the point-crossing probability given in Proposition \ref{prop:pointdecomp} in sub-sequent subsections. We first treat the probability of the event $\cal F$ from Proposition \ref{prop:decomp} conditioned on non-intersection events in the vicinity of $B_i$ and $B_j$, which we denote by ${\cal G}^1$ and ${\cal G}^2$ below in \eqref{eq:G1def} and \eqref{eq:G2def} respectively and show that this conditional probability can be approximately expressed through a function of configurations near $B_i$ and $B_j$; see Proposition \ref{eq:Hdecomp1} for more details. Then, as an important intermediate step, in Proposition \ref{prop:Lem11} we show that only the ``end part'' (from inside $B_i$ and $B_j$ to outside) of these configurations matters for calculation of the conditional probability of $\cal F$. As discussed at the end of Section \ref{sec:sketch}, this is a crucial fact that allows us to apply coupling results in Section \ref{sec:4.6}.
As the proof of Proposition \ref{prop:Lem11} is long and technical, we postpone it (along with that of Lemma \ref{lem:FFocompare}) to Section \ref{sec:endingpart}.

 Recall the definition of the (conditioned) random walks,  $U_{i,-m}$, $U_{j,-m}$,  $\gamma^1$, $\gamma^2$ and $V_{1}, V_{2}$ in \eqref{eq:Udef},   \eqref{eq:gammadef} and \eqref{eq:V1V2def} respectively. With the event ${\cal F}$ defined in \eqref{eq:Hdef} in mind, we now define for $y\in\partial_i B_i$ and  $z\in\partial_i B_j$ the following events 
\begin{equation}\label{eq:G1def}
{\cal G}^1:={\cal G}^1(y)=\Big\{ \gamma^1[V_1,U_{i,-k^3}]\cap
(X^2)^{\cal R} [1,U_{i,-k^3}]=\emptyset,\; (X^2)^{\cal R} [1,U_{i,-k^3}]\cap B_i=\emptyset\Big\}
\end{equation}
and
\begin{equation}\label{eq:G2def}
{\cal G}^2:={\cal G}^2(z)=\Big\{ \gamma^2[V_2,U_{j,-k^3}]\cap X^3[1, U_{j,-k^3}]=\emptyset,\; X^3[1, U_{j,-k^3}]\cap B_j=\emptyset\Big\},
\end{equation}
that dictate the non-intersection of SRW and LERW near $B_i$ and $B_j$.

From now on in this subsection we fix $y\in\partial_i B_i$ and  $z\in\partial_i B_j$ and abbreviate $\mP^{y,z,z}$ as $\mP$. Note that by definition ${\cal F}(y,z)\subset{\cal G}^1\cap {\cal G}^2$. 
We observe that by Lemma \ref{asymp-indep}, 
\begin{equation}\label{eq:G12localization}
\mP({\cal G}^1\cap {\cal G}^2)\simeq \mP({\cal G}^1) \mP({\cal G}^2).
\end{equation}

We now let 
$$
\Gamma_i=\left\{ (\fraka,\frakb) \Bigg| \begin{split}
\fraka:\mbox{a simple path,} \; &\fraka[0,\Diamond)\subset B(x_i, 2^{-k^3}), \; \fraka(\Diamond)\in \partial B(x_i, 2^{-k^3});\\
 \frakb:\mbox{a path},\;\;\;\;\quad\quad\; & \frakb[0,\Diamond) \subset B(x_i,2^{-k^3}), \;b(\Diamond) \in \partial B(x_i, 2^{-k^3}).
\end{split} \right\},
$$
and define $\Gamma_j$ similarly by replacing $\fraka,\frakb,i$ by $\frakc,\frakd,j$ respectively in the definition above. 
Then, consider a set of pairs of paths $(\fraka,\frakb)$:
$$
\Delta_i=\Delta_i(y):=\left\{(\fraka,\frakb)\in\Gamma_i \mid \mP\Big( \Big\{\big(\gamma^1[0,U_{i,-k^3}],(X^2)^{\cal R}[0,U_{i,-k^3}]\big)=(\fraka,\frakb)\Big\}\cap{\cal G}^1\Big)>0\right\},
$$
and define $\Delta_j=\Delta_j(z)$ similarly. Note that $\Delta_{i}$ is just the subset of $\Gamma_i$ containing pairs of paths $(\fraka,\frakb)$ starting from $y$, such that $\fraka$ and $\frakb$ are non-intersecting and $\frakb$ avoids $B_i$. 
We then define events
\begin{itemize}
\item ${\cal G}_{\fraka,\frakb}:=\Big\{\big(\gamma^1[0,U_{i,-k^3}],(X^2)^{\cal R}[0,U_{i,-k^3}]\big)=(\fraka,\frakb)\Big\}$ for $(\fraka,\frakb)\in\Delta_i$;
\item ${\cal G}_{\frakc,\frakd}:=\Big\{\big(\gamma^2[0,U_{j,-k^3}],X^3[0,U_{j,-k^3}]\big)=(\frakc,\frakd)\Big\}$ for $(\frakc,\frakd)\in\Delta_j$.
\end{itemize}
Similar to \eqref{eq:G12localization}, by Lemma \ref{asymp-indep}  we also have
\begin{equation}\label{eq:GabGcd}
\mP({\cal G}_{\fraka,\frakb}\cap {\cal G}_{\frakc,\frakd})\simeq\mP({\cal G}_{\fraka,\frakb})\mP({\cal G}_{\frakc,\frakd}).
\end{equation}
 Then, it follows that
\begin{equation}\label{eq:represent}
\mP[{\cal F}| {\cal G}^1\cap {\cal G}^2]
\overset{\eqref{eq:G12localization}}{\underset{\eqref{eq:GabGcd}}{\simeq}}\sum_{(\fraka,\frakb)\in \Delta_i} \sum_{(\frakc,\frakd)\in \Delta_j}\mP({\cal F} \big| {\cal G}_{\fraka,\frakb}\cap {\cal G}_{\frakc,\frakd})\times \frac{\mP[{\cal G}_{\fraka,\frakb}]\mP[{\cal G}_{\frakc,\frakd}]}{\mP[{\cal G}^1]\mP[{\cal G}^2]} .
\end{equation}

We now rewrite $\mP({\cal F} \big| {\cal G}_{\fraka,\frakb}\cap {\cal G}_{\frakc,\frakd})$ into a function of the quadruple $(\fraka,\frakb,\frakc,\frakd)$.

\begin{figure}[h]
\begin{center}
\includegraphics[scale=0.55]{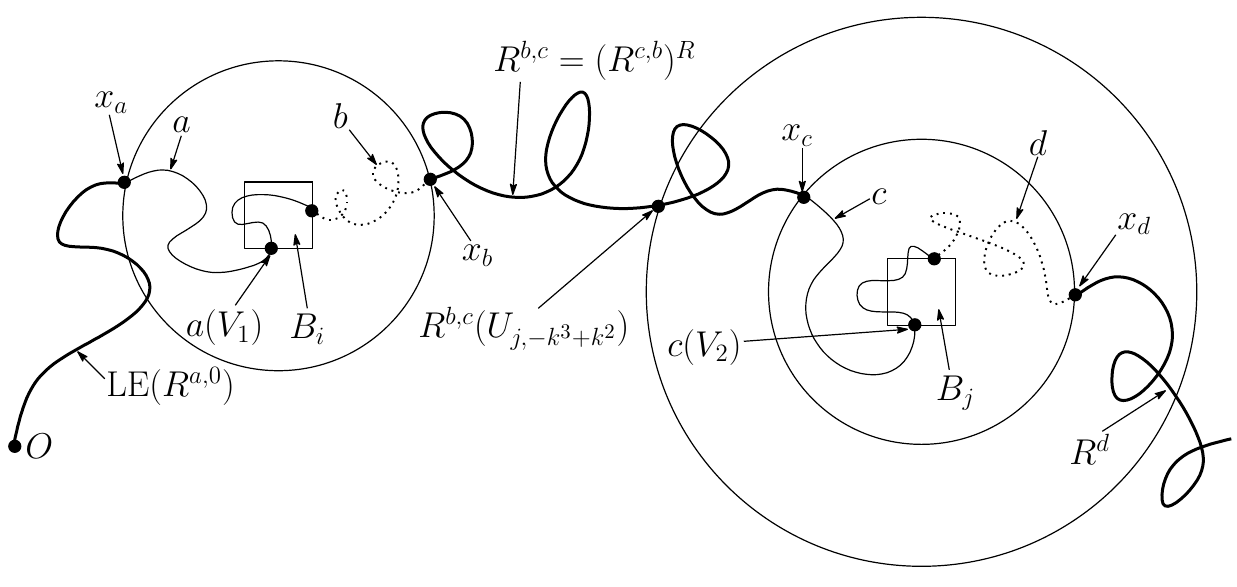}
\caption{Illustration for  the setup in Definition \ref{defnf}. The thin solid curves starting from the boundary of $B_{i}$ and $B_{j}$ stand for $\fraka$ and $\frakc$ respectively, while the thin dotted curves starting from the boundary of $B_{i}$ and $B_{j}$ stand for $\frakb$ and $\frakd$ respectively. The thick curves starting from $x_{\fraka}$, $x_{\frakb}$ and $x_{\frakd}$ represent $\LE ( R^{\fraka,0} )$, $R^{\frakb,\frakc}$ and $R^{\frakd}$ respectively. For technical reasons the paths $\fraka$ through $\frakd$ are marked in roman font $a$ through $d$ respectively.}\label{fig35}
\end{center}
\end{figure}

\begin{dfn}\label{defnf}
For $(\fraka,\frakb)\in\Gamma_i$ and $(\frakc,\frakd)\in\Gamma_j$ with $x_{\aleph}$, $\aleph=\fraka,\frakb,\frakc,\frakd$, the respective ending point of these paths, let $R^{\fraka,0}$, $R^{\frakc,\frakb}$ and $R^\frakd$ be independent random walks such that 
\begin{itemize}
\item $R^{\fraka,0}$  starts from $x_\fraka$ and is conditioned to hit the origin (and stopped there) before hitting $\fraka$ and $\partial \mathbb{D}$;
\item $R^{\frakc,\frakb}$ starts from $x_\frakc$ and is conditioned to hit $x_\frakb$ before leaving $\mathbb{D}$ truncated at the last time it hits $x_\frakb$ and further conditioned that $R^{\frakc,\frakb}[1,\Diamond]\cap \frakc =\emptyset$;
\item $R^\frakd$ starts from $x_\frakd$ stopped at first exit from $\mathbb{D}$.
\end{itemize}
For convenience we also write $R^{\frakb,\frakc}=(R^{\frakc,\frakb})^{\cal R}$ for the reverse of $R^{\frakc,\frakb}$. Write $\mathbb{Q}$ for their joint law.

Recall the definition of random times $V_1$ and $V_2$ for $\gamma_1$ and $\gamma_2$ in \eqref{eq:V1V2def}. With slight abuse of notation, we also write $V_1$ and $V_2$ for similarly defined times for $\fraka$ and $\frakc$, with the convention that $\max\emptyset=0$.
We now let
\begin{equation}\label{eq:fdef}
f(\fraka,\frakb,\frakc,\frakd):=\mathbb{Q}\left(
\begin{split}
& \big(\fraka[V_1,\Diamond]\cup \LE(R^{\fraka,0})\big)\bigcap \big(R^{\frakb,\frakc}[0,U_{j,-k^3+k^2}]\cup \frakb\cup R^\frakd\cup \frakd\big)=\emptyset; \\
& \big(\LE(R^{\frakc,\frakb}) \cup \frakc[V_2,\Diamond]\big)\bigcap (R^\frakd\cup \frakd\big)=\emptyset 
\end{split}
\right).
\end{equation}
See Figure \ref{fig35} for the setup. 
\end{dfn}

By the domain Markov property of LERW and the strong Markov property of the simple random walk, for $(\fraka,\frakb) \in \Delta_i$ and $(\frakc,\frakd)\in \Delta_j$, we have 
\begin{equation}\label{eq:Hrep}
f(\fraka,\frakb,\frakc,\frakd)=\mP({\cal F}^o\big| {\cal G}_{\fraka,\frakb} \cap {\cal G}_{\frakc,\frakd}),
\end{equation}
where ${\cal F}^o={\cal F}^o(y,z)$ is defined through
\begin{equation}\label{eq:Fodef}
{\cal F}^o:=\left\{\begin{aligned}
&\gamma^1[V_1,\Diamond] \cap \big((X^2)^{\cal R}[1, U_{j,-k^3+k^2}]\cup X^3\big) =\emptyset; \; (X^2)^{\cal R}[1, U_{i,-k^3}]\cap B_i=\emptyset; \\
& \gamma^2[V_2,\Diamond
] \cap X^3[1,\Diamond]=\emptyset; \qquad\qquad\qquad\qquad\;\;\;
 X^3[1,U_{j,-k^3}]\cap B_j =\emptyset
\end{aligned}\right\}.
\end{equation}
Note that the difference between the definition of ${\cal F}^o$ and $\cal F$ (see \eqref{eq:Hdef}) is that for the former we only require the part of $(X^2)^{\cal R}$ and $X^3$ before exiting a cube of size $2^{-k^3}$ 
not to intersect $B_i$ and $B_j$ respectively (instead of  (almost) the whole path for the latter). 

We now claim that this tiny difference does not change the (conditional) probability drastically. More precisely,
\begin{lem}\label{lem:FFocompare}
For $(\fraka,\frakb) \in \Delta_i$ and $(\frakc,\frakd)\in \Delta_j$,
\begin{equation}\label{eq:FFocompare}
\mP({\cal F}\big| {\cal G}_{\fraka,\frakb} \cap {\cal G}_{\frakc,\frakd})\simeq\mP({\cal F}^o\big| {\cal G}_{\fraka,\frakb} \cap {\cal G}_{\frakc,\frakd}).
\end{equation}
\end{lem}
We postpone the proof to Section \ref{sec:endingpart}, as it is a bit technical (and tedious if we spell out all the details), and moreover in some sense reminiscent of some part of the proof of Proposition \ref{prop:Lem11}, which will also be laid out in that subsection.

Combining \eqref{eq:represent}, \eqref{eq:Hrep} and \eqref{eq:FFocompare}, we obtain that 
\begin{prop}\label{prop:Hdecomp1} 
\begin{equation}\label{eq:Hdecomp1}
\mP({\cal F} \big| {\cal G}^1\cap {\cal G}^2) \simeq \sum_{(\fraka,\frakb)\in \Delta_{i}} \sum_{(\frakc,\frakd)\in \Delta_{j}} f(\fraka,\frakb,\frakc,\frakd) \frac{\mP({\cal G}_{\fraka,\frakb})}{\mP({\cal G}^1)}\frac{\mP({\cal G}_{\frakc,\frakd})}{\mP({\cal G}^2)}.
\end{equation}
\end{prop}

\medskip

We now turn to decoupling. We start with the definition of the law of a pair of paths which can be heuristically understood as the ``location-free'' local law of the LERW passing through a cube.
\begin{itemize}
\item Let 
$\gamma^y_\infty$, $S^y_2$, $\gamma^z_\infty$ and $S^z_3$ be independent ILERW's and SRW's started at $y$ and $z$ respectively. We write $\mP_\infty$ for their joint law.
We now define the following events (compare them with the definition of ${\cal G}^1$ and ${\cal G}^2$ at \eqref{eq:G1def} and \eqref{eq:G2def}):
\begin{equation}\label{eq:G1infdef}
{\cal G}^1_\infty={\cal G}^1_\infty(y):=\left\{ \gamma^y_\infty[V_1,U_{i,{-k^3}}]\cap S^y_2[1,U_{i,{-k^3}}]=\emptyset,\; S^y_2[1,U_{i,{-k^3}}]\cap B_i=\emptyset
\right\};
\end{equation}
\begin{equation}\label{eq:G2infdef}
{\cal G}^2_\infty={\cal G}^2_\infty(z):=\left\{ \gamma^z_\infty[V_2,U_{j,{-k^3}}]\cap S^z_3[1,U_{j,{-k^3}}]=\emptyset,\; S^z_3[1,U_{j,{-k^3}}]\cap B_j=\emptyset
\right\},
\end{equation}
where we again abuse the notation $V_1$ and $V_2$ originally defined for $\gamma^1$ and $\gamma^2$ in \eqref{eq:V1V2def} and use them to denote the corresponding random times for $\gamma^y_\infty$ and $\gamma^z_\infty$.  
By \eqref{eq:LERWvsILERW}, one has for any $y\in \partial_i B_i$ and $w\in \partial_i B_j$
\begin{equation}\label{eq:G12asymp}
\mP\big({\cal G}^1(y)\big)\simeq \mP\big({\cal G}^1_\infty(y)\big)\quad \mbox{ and }\quad \mP\big({\cal G}^2(z)\big) \simeq  \mP\big({\cal G}^2_\infty(z)\big). 
\end{equation}
\item We then define $\nu^y(\cdot,\cdot)$ for the probability measure induced by $(\gamma^y_\infty[0,U_{i,-k^3}], S^y_2[0,U_{i,-k^3}])$
conditioned on the event ${\cal G}^1_\infty$,  and define $\nu^z(\cdot,\cdot)$ similarly with conditioning on ${\cal G}^2_\infty$. Note that the distribution of paths (modulo starting points) under $\nu^y$ and $\nu^z$ only depends on the relative location of $y$ and $z$ with respect to $x_i$ and $x_j$, and not on the choice of $i$ or $j$ as long as the cube $B$ we consider in this section satisfies \eqref{eq:Breq}.
\end{itemize}

We now claim that \eqref{eq:Hdecomp1} can be rewritten as follows:
\begin{lem}  \label{lem:asymp}
 \begin{equation}\label{eq:Hdecomp2}
\mP({\cal F} \big| {\cal G}^1\cap {\cal G}^2)\simeq \sum_{(\fraka,\frakb)\in \Gamma_i} \sum_{(\frakc,\frakd)\in \Gamma_j} f(\fraka,\frakb,\frakc,\frakd) \nu^y(\fraka,\frakb)\nu^z(\frakc,\frakd).
\end{equation}
\end{lem}
\begin{proof}
By \eqref{eq:asymp-variant}, the Radon-Nikodym derivative of the law of $\gamma^y_\infty[0,U_{i,{-k^3}}]$ 
with respect to that of $\gamma^1[0,U_{i,-k^3}]$, is uniformly $1+O(2^{-ck^3})$ 
while the law of 
$(X^2)^{\cal R}[0,U_{i,-k^3}]$ with respect to $S^y_2[0,U_{i,-k^3}]$ is also uniformly $1+O(2^{-ck^3})$ by basic applications of the Harnack principle. In a similar fashion, same bounds on Radon-Nikodym derivatives also exist for $\big(\gamma^2[0,U_{j,-k^3}],X^3[0,U_{j,-k^3}]\big)$ and $\big(\gamma^z_\infty[0,U_{j,{-k^3}}],S^z_3[0,U_{j,{-k^3}}]\big)$.
Therefore for any $(\fraka,\frakb)\in \Delta_i$ and $(\frakc,\frakd)\in \Delta_j$, 
$$
{\mP({\cal G}_{\fraka,\frakb})}/{\mP({\cal G}^1)}\simeq \nu^y(\fraka,\frakb)\;\mbox{ and }\;{\mP({\cal G}_{\frakc,\frakd})}/{\mP({\cal G}^2)}\simeq\nu^z(\frakc,\frakd).$$
Hence \eqref{eq:Hdecomp2} follows from \eqref{eq:Hdecomp1}.
\end{proof}

We now claim that, loosely speaking, $\mP\big({\cal F} \big| {\cal G}^1\cap {\cal G}^2\big)$ is almost a function of the ``end part'' of $(\fraka,\frakb,\frakc,\frakd)$ only. To be more precise, we fix some notation for the truncation of paths.
More precisely, for $(\fraka,\frakb)\in\Gamma_i$ and $(\frakc,\frakd)\in\Gamma_j$, let 
\begin{equation}\label{eq:pidef}
\pi(\fraka,\frakb):= \big( a[U_{i,{-10 k^3}}, \Diamond], b[U_{i,{-10 k^3}}, \Diamond]\big);\quad \pi(\frakc,\frakd):=\big(c[U_{j,{-10 k^3}}, \Diamond],d[U_{j,{-10 k^3}}, \Diamond]\big)
\end{equation}
be the end part of $(\fraka,\frakb)$ and $(\frakc,\frakd)$ respectively. We also define $\pi(\fraka,\frakb,\frakc,\frakd)$ in the same fashion. We also define the space of truncated paths as follows. 
\begin{equation}\label{eq:Gammabardef}
\overline{\Gamma}_i:=\left\{ (\overline{\fraka},\overline{\frakb}) \big|  (\overline{\fraka},\overline{\frakb})\in\Gamma_i,\;  \overline{\fraka}(0),\overline{\frakb}(0)\in \partial B(x_i, 2^{-10k^3}) \right\},\mbox{ and define $\overline{\Gamma}_j$  similarly}.
\end{equation}
We then write 
\begin{itemize}
\item $\overline{\nu}^{y}$ for the probability measure on $\overline{\Gamma}_i$ induced by $
\pi\big(\gamma^y_\infty[0,U_{i,{-k^3}}], S^y_2[0,U_{i,{-k^3}}]\big)
$ conditioned on ${\cal G}^{1}_\infty$;
\item $\overline{\nu}^{z}$ for the probability measure on $\overline{\Gamma}_j$ induced by $\pi\big(\gamma^z_\infty[0,U_{j,{-k^3}}], S^z_3[0,U_{j,{-k^3}}]\big)$
conditioned on ${\cal G}^{2}_\infty$.
\end{itemize}

We are ready to give the precise statement of what we claim at the beginning of last paragraph. As its proof is quite long, we postpone it to Section \ref{sec:endingpart}.
\begin{prop}\label{prop:Lem11}
\begin{equation}\label{eq:Lem11}
\mP\big({\cal F} \big| {\cal G}^1\cap {\cal G}^2\big) \simeq \sum_{(\overline{\fraka},\overline{\frakb})\in \overline{\Gamma}_i} \sum_{(\overline{\frakc},\overline{\frakd})\in \overline{\Gamma}_j} f(\overline{\fraka},\overline{\frakb},\overline{\frakc},\overline{\frakd}) \overline{\nu}^y(\overline{\fraka},\overline{\frakb})\overline{\nu}^z(\overline{\frakc},\overline{\frakd}). 
\end{equation}
\end{prop}
Finally, we state the main result of this subsection, which now follows as a quick corollary of Proposition \ref{prop:Lem11}.
\begin{prop}\label{prop:boxdecomp}
\begin{equation}\label{eq:H1}
\begin{split}
&P(0\xrightarrow{\gamma} B_i\xrightarrow{\gamma} B_j)\\
\simeq\;&\sum_{y\in\partial B_i}\sum_{z\in\partial B_j}G_{\mathbb{D}}(0,y)G_{\mathbb{D}}(y,z)         
\sum_{(\overline{\fraka},\overline{\frakb})\in \overline{\Gamma}_i} \sum_{(\overline{\frakc},\overline{\frakd})\in \overline{\Gamma}_j} f(\overline{\fraka},\overline{\frakb},\overline{\frakc},\overline{\frakd}) \overline{\nu}^y(\overline{\fraka},\overline{\frakb})\overline{\nu}^z(\overline{\frakc},\overline{\frakd}) 
\mP({\cal G}^1)\mP({\cal G}^2)
\end{split}
\end{equation}
\end{prop}
\begin{proof}
By Proposition \ref{prop:decomp}, \eqref{eq:G12localization} and Proposition \ref{prop:Lem11},  we have 
\begin{equation*}
\begin{split}
&\; P(0\xrightarrow{\gamma} B_i\xrightarrow{\gamma} B_j)\overset{\eqref{eq:decomp}}{\simeq}\sum_{y\in\partial B_i}\sum_{z\in\partial B_j} G_{\mathbb{D}}(0,y)G_{\mathbb{D}}(y,z) \mP({\cal F}) \\
\overset{\eqref{eq:G12localization}}{\simeq}\;&\sum_{y\in\partial B_i}\sum_{z\in\partial B_j}G_{\mathbb{D}}(0,y)G_{\mathbb{D}}(y,z)\mP({\cal F}\;|\;{\cal G}^1\cap{\cal G}^2)\mP({\cal G}^1)\mP({\cal G}^2)\\
\overset{\eqref{eq:Lem11}}{\simeq}\;&\sum_{y\in\partial B_i}\sum_{z\in\partial B_j}G_{\mathbb{D}}(0,y)G_{\mathbb{D}}(y,z)         
\sum_{(\overline{\fraka},\overline{\frakb})\in \overline{\Gamma}_i} \sum_{(\overline{\frakc},\overline{\frakd})\in \overline{\Gamma}_j} f(\overline{\fraka},\overline{\frakb},\overline{\frakc},\overline{\frakd}) \overline{\nu}^y(\overline{\fraka},\overline{\frakb})\overline{\nu}^z(\overline{\frakc},\overline{\frakd}) 
\mP({\cal G}^1)\mP({\cal G}^2).
\end{split}
\end{equation*}
This finishes the proof.
\end{proof}

\subsection{The point-crossing probability}\label{sec:pt}

In this subsection, we are going to deal with the probability of the point-crossing event on the LHS of \eqref{eq:MAIN}. As the situation is very similar to that of cube-crossing probabilities, we will skip the proof but only give analogues of Propositions  \ref{prop:decomp} and \ref{prop:boxdecomp} (see Propositions \ref{lem:Lem12} and \ref{prop:pointdecomp} below) which will be directly quoted in the proof of Proposition \ref{prop4.6}.

As in previous subsections, we still assume that $(i,j)$ is good in the sense of \eqref{goodbad} throughout this subsection and fix $v\in B_i$ and $w\in B_j$.
Let
\begin{itemize}
\item $Y^1$ be the random walk started from $v$, conditioned to hit $0$ (and stopped there) before leaving $\mathbb{D}$;
\item $Y^2$ be the random walk started from $w$, conditioned to hit $v$ (and stopped there) before leaving $\mathbb{D}$;
\item $Y^3$ be the simple random walk started from $w$ stopped upon leaving $\mathbb{D}$.
\end{itemize}
We write $\widetilde{\mP}$ for their joint distribution.
Then a classical last-exit decomposition argument (see e.g.\ Proposition 8.1 of \cite{S}) shows that 

\begin{equation}\label{eq:pointdecomp}
P\big(0\xrightarrow{\gamma}v\xrightarrow{\gamma}w\big)=G_{\mathbb{D}}(0,v)G_{\mathbb{D}}(v,w) \widetilde{\mP}\Big(\big\{  \LE(Y^1)\cap \big(Y^2[1,\Diamond)\cup Y^3\big)=\emptyset,\; \LE(Y^2)\cap Y^3[1,\Diamond]=\emptyset   \big\} \Big).
\end{equation}
Similar to the definition of event ${\cal F}$ in \eqref{eq:Hdef}, we define the following event:
\begin{equation}\label{widetildeH}
\widetilde{\cal F}:=\Big\{\LE(Y^1)\cap \big(Y^2[U_{j,{-k^3+k^2}},\Diamond)\cup Y^3\big)=\emptyset,\; \LE\big(Y^2[0,\widehat{\sigma}]\big)\cap Y^3[1,\Diamond)=\emptyset {\big)} \Big\}.
\end{equation}
where (compare this to \eqref{eq:sigmadef})
$\widehat{\sigma}:=\max\big\{t\mid Y^2(t)\in\partial B(x_i,2^{-k^3})\big\}$.

We now claim that  without incurring a big error, we can replace the event on the RHS of  \eqref{eq:pointdecomp} by $\widetilde{\cal F}$ defined in \eqref{widetildeH}. 
\begin{prop}\label{lem:Lem12}For any good $(i,j)$ and any $v\in B_i$ and $w\in B_j$, 
\begin{equation}
P(0\xrightarrow{\gamma}v\xrightarrow{\gamma}w)\simeq G_{\mathbb{D}}(0,v)G_{\mathbb{D}}(v,w)\widetilde\mP(\widetilde{\cal F}).
\end{equation}
\end{prop}
We omit the proof as it can be proved in a similar way as Proposition \ref{prop:decomp}.

\medskip

We now define (as analogues of \eqref{eq:G1def} and \eqref{eq:G2def}) $\widetilde{\cal G}^1=\widetilde{\cal G}^1(v)$ and $\widetilde{\cal G}^2=\widetilde{\cal G}^2(w)$ where
\begin{equation}\label{eq:Gt12infdef}
\widetilde{\cal G}^1:=\Big\{\LE(Y^1)[0,U_{i,-k^3}]\cap (Y^2)^{\cal R}[1,U_{i,-k^3}]=\emptyset\Big\}\mbox{; }
\widetilde{\cal G}^2:=\Big\{\LE(Y^2[0,\widehat{\sigma}])[0,U_{j,-k^3}]\cap Y^3[1,U_{j,-k^3}]=\emptyset\Big\}.
\end{equation}

Recall the definition of  $\pi$,
$\overline{\Gamma}_i$ and $\overline{\Gamma}_j$ in \eqref{eq:pidef} and \eqref{eq:Gammabardef} respectively.  We let $\overline{\mu}^v(\cdot,\cdot)$ be the probability measure on $\overline{\Gamma}_i$ induced by 
$
\pi\big(\gamma_{\infty}^{v}[0,U_{i,-k^3}], S^v[0,U_{i,-k^3}]\big)
$
governed by $\widetilde\mP^v_\infty(\;\cdot \mid \widetilde{\cal G}^v_\infty)$, 
and define $\overline{\mu}^w(\cdot,\cdot)$ similarly. 
As an analogue of Proposition \ref{prop:boxdecomp}, we have the following decomposition of the point-crossing probability. We omit the proof.
\begin{prop}\label{prop:pointdecomp} For any good $(i,j)$ and any $v\in B_i$ and $w\in B_j$, 
\begin{equation}\label{eq:Ht1}
P(0\xrightarrow{\gamma} v\xrightarrow{\gamma} w)\simeq G_{\mathbb{D}}(0,v)G_{\mathbb{D}}(v,w)\widetilde\mP(\widetilde{\cal G}^1(v))\widetilde\mP(\widetilde{\cal G}^2(w))\sum_{(\overline{\fraka},\overline{\frakb})\in \overline{\Gamma}_i} \sum_{(\overline{\frakc},\overline{\frakd})\in \overline{\Gamma}_j} f(\overline{\fraka},\overline{\frakb},\overline{\frakc},\overline{\frakd}) \overline{\mu}^v(\overline{\fraka},\overline{\frakb})\overline{\mu}^w(\overline{\frakc},\overline{\frakd}).
\end{equation}
Here $f$ is defined as in \eqref{eq:fdef}.
\end{prop}

We now define, as analogue of \eqref{eq:G1infdef} and \eqref{eq:G2infdef},
\begin{equation}\label{eq:Ginfdef}
\widetilde{\cal G}_\infty:=\Big\{\gamma_\infty[0,U_{-k^3}]\cap S[1,U_{-k^3}]=\emptyset\Big\},
\end{equation}
where  $\gamma_\infty$ and $S$ are an infinite LERW and SRW, independent from each other, both started at the origin whose joint law we denote by $\widetilde\mP_\infty$ and for each path $U_{-k^3}$ is the first time the path exits $B(0,2^{-k^3})$, defined in a similar fashion as in \eqref{eq:Udef}. We use superscripts such as $(\gamma^x_\infty,S^x)$ and $\widetilde\mP^x_\infty$ for $x\in 2^{-n}\mathbb{Z}^3$ to represent a translate of what we have just defined. By applying \eqref{ES-infty},
we see that 
\begin{equation}\label{eq:Gtinfasymp}
\widetilde{\mP}_\infty (\widetilde{\cal G}_\infty) \asymp 2^{(2-\beta)(k^3-n)}.
\end{equation}
By an decomposition which is similar to Proposition \ref{prop:pointdecomp} but essentially simpler, the estimate \eqref{TRANSINVA}, and \eqref{eq:LERWvsILERW} to transfer between the beginning parts of $\LE(Y^1)$, $\LE(Y^2)$ and $\gamma_\infty$, one has for any $v\in B_i$ and $w\in B_j$
\begin{equation}\label{eq:Gt12infasymp}
\widetilde\mP(\widetilde{\cal G}^1(v))\simeq \widetilde\mP(\widetilde{\cal G}^1(x_i))\simeq \widetilde{\mP}_\infty (\widetilde{\cal G}_\infty) \simeq  \widetilde\mP(\widetilde{\cal G}^2(x_j)) \simeq  \widetilde\mP(\widetilde{\cal G}^2(w)). 
\end{equation}

\subsection{Proof of Proposition \ref{prop:goodpair}, the main estimate}\label{sec:4.6}
With Propositions \ref{prop:boxdecomp} and \ref{prop:pointdecomp} in hand, we are now ready to give the proof of the main result of this section.

\begin{proof}[Proof of Proposition \ref{prop:goodpair}] We only prove \eqref{eq:MAIN} as the proof for \eqref{eq:MAINprime} is very similar. Recall the definitions of $\overline{\nu}^\bullet$ and $\overline{\mu}^\bullet$ at \eqref{eq:Gammabardef} and \eqref{eq:Ginfdef} respectively.
Comparing Propositions \ref{prop:boxdecomp} and \ref{prop:pointdecomp}, we see that in order to prove \eqref{eq:MAIN}, it suffices to compare $\overline{\mu}^v$ with $\overline{\nu}^y$ for $v\in B_i$, $y\in\partial B_i$, and $\overline{\mu}^w$ with $\overline{\nu}^z$ for $w\in B_j$, $z\in\partial B_j$. Recall that $x_i$ and $x_j$ are the centers of $B_i$ and $B_j$ respectively. 
We observe that by \eqref{coupling2}, we have
\begin{equation}\label{eq:allcoup1}
||\overline{\mu}^v -\overline{\mu}^{x_i}|| \leq c\cdot 2^{-ck^3};
\;
||\overline{\mu}^w -\overline{\mu}^{x_j}|| \leq c\cdot 2^{-ck^3},
\end{equation}
and by \eqref{coupling2}, we have
\begin{equation}\label{eq:allcoup2}
||\overline{\nu}^y -\overline{\mu}^{x_i}||\leq c\cdot 2^{-ck^3};
\;
||\overline{\nu}^z -\overline{\mu}^{x_j}||\leq c\cdot 2^{-ck^3},
\end{equation}
(note that in this case, in the notation of Section \ref{sec:decoupling}, esp.\  \eqref{eq:G1infdef} and \eqref{eq:G2infdef}, the quadruple $(\gamma,\lambda,\Theta_\gamma, \Theta_\lambda)$ in \eqref{LAMBDAbar-k} and \eqref{coupling2} should be taken as 
\begin{equation}\label{eq:couplinginit}
(\gamma^y_\infty[V_1, U_{i, -k^4+k^3}], S^y_2[0, U_{i,-k^4+k^3}],\gamma^y_\infty[0,v_1],B_i) \mbox{ and }(\gamma^z_\infty[V_2, U_{j, -k^4+k^3}], S^z_3[0, U_{j,-k^4+k^3}],\gamma^z_\infty[0,V_2],B_j)
\end{equation}
respectively).
By \eqref{eq:allcoup1} and the observation that $\overline{\mu}^\bullet$ are probability measures, we have
\begin{equation}\label{eq:errorest}
\begin{split}
&\bigg|\sum_{(\overline{\fraka},\overline{\frakb})\in \overline{\Gamma}_i} \sum_{(\overline{\frakc},\overline{\frakd})\in \overline{\Gamma}_j} f(\overline{\fraka},\overline{\frakb},\overline{\frakc},\overline{\frakd})  \Big(\overline{\mu}^{v}(\overline{\fraka},\overline{\frakb})\overline{\mu}^{w}(\overline{\frakc},\overline{\frakd})-  \overline{\mu}^{x_i}(\overline{\fraka},\overline{\frakb})\overline{\mu}^{x_j}(\overline{\frakc},\overline{\frakd})\Big)\bigg|\\
\leq \;&C 2^{-ck^3} \max_{(\overline{\fraka},\overline{\frakb})\in \overline{\Gamma}_i,\; (\overline{\frakc},\overline{\frakd})\in \overline{\Gamma}_j} f(\overline{\fraka},\overline{\frakb},\overline{\frakc},\overline{\frakd})\overset{(*)}{\leq} C' 2^{-c'k^3} \sum_{(\overline{\fraka},\overline{\frakb})\in \overline{\Gamma}_i} \sum_{(\overline{\frakc},\overline{\frakd})\in \overline{\Gamma}_j} f(\overline{\fraka},\overline{\frakb},\overline{\frakc},\overline{\frakd})  \overline{\mu}^{x_i}(\overline{\fraka},\overline{\frakb})\overline{\mu}^{x_j}(\overline{\frakc},\overline{\frakd}).
\end{split}
\end{equation}
where in $(*)$ we apply the following up-to-constants asymptotics:
\begin{lem}
Recall the definition of $\epsilon$ and $r$ in \eqref{eq:rDEF} and let $l=|x_i-x_j|$.
We have
\begin{equation}\label{eq:uptocst1}
\max_{(\overline{\fraka},\overline{\frakb})\in \overline{\Gamma}_i,\; (\overline{\frakc},\overline{\frakd})\in \overline{\Gamma}_j} f(\overline{\fraka},\overline{\frakb},\overline{\frakc},\overline{\frakd}) \asymp \Big( \frac{2^{-k^{3}}}{l}\Big)^{2- \beta} \cdot \Big(\frac{2^{-k^{3}}}{r} \Big)^{2-\beta},
\end{equation}
and 
\begin{equation}\label{eq:uptocst2}
\sum_{(\overline{\fraka},\overline{\frakb})\in \overline{\Gamma}_i} \sum_{(\overline{\frakc},\overline{\frakd})\in \overline{\Gamma}_j} f(\overline{\fraka},\overline{\frakb},\overline{\frakc},\overline{\frakd})  \overline{\mu}^{x_i}(\overline{\fraka},\overline{\frakb})\overline{\mu}^{x_j}(\overline{\frakc},\overline{\frakd}) \asymp \Big( \frac{2^{-k^{3}}}{l}\Big)^{2- \beta} \cdot \Big(\frac{2^{-k^{3}}}{r} \Big)^{2-\beta}.
\end{equation}
\end{lem}
\begin{proof} 
The asymptotics \eqref{eq:uptocst1} follows by first discarding all requirements in the event defining $f$ (see the RHS of \eqref{eq:fdef}) on random walks (i.e., $R^{\fraka,0}$,$R^{\frakc,\frakb}$ and $R^\frakd$) not intersecting $\fraka$, $\frakb$, $\frakc$ and $\frakd$ respectively and then apply  an argument similar to the proof of Proposition \ref{prop:bound-box-hit} to obtain an up-to-constants asymptotics.

The asymptotics \eqref{eq:uptocst2} follows e.g., from observing that the LHS of \eqref{eq:uptocst2} is exactly the last term on the RHS of \eqref{eq:Ht1}, and then applying \eqref{M7prime}, \eqref{green-1} and \eqref{eq:Gt12infasymp} for up-to-constants bounds on other terms.
\end{proof}

Coming back to the proof of Proposition \ref{prop:goodpair}, by \eqref{eq:errorest}, we have
\begin{equation}\label{eq:redu1}
\sum_{(\overline{\fraka},\overline{\frakb})\in \overline{\Gamma}_i} \sum_{(\overline{\frakc},\overline{\frakd})\in \overline{\Gamma}_j} f(\overline{\fraka},\overline{\frakb},\overline{\frakc},\overline{\frakd})  \overline{\mu}^{v}(\overline{\fraka},\overline{\frakb})\overline{\mu}^{w}(\overline{\frakc},\overline{\frakd})
\simeq\sum_{(\overline{\fraka},\overline{\frakb})\in \overline{\Gamma}_i} \sum_{(\overline{\frakc},\overline{\frakd})\in \overline{\Gamma}_j} f(\overline{\fraka},\overline{\frakb},\overline{\frakc},\overline{\frakd})  \overline{\mu}^{x_i}(\overline{\fraka},\overline{\frakb})\overline{\mu}^{x_j}(\overline{\frakc},\overline{\frakd}).
\end{equation}
Similarly, by \eqref{eq:allcoup2} we have 
\begin{equation}\label{eq:redu2}
\sum_{(\overline{\fraka},\overline{\frakb})\in \overline{\Gamma}_i} \sum_{(\overline{\frakc},\overline{\frakd})\in \overline{\Gamma}_j} f(\overline{\fraka},\overline{\frakb},\overline{\frakc},\overline{\frakd})  \overline{\nu}^{y}(\overline{\fraka},\overline{\frakb})\overline{\nu}^{z}(\overline{\frakc},\overline{\frakd})
\simeq\sum_{(\overline{\fraka},\overline{\frakb})\in \overline{\Gamma}_i} \sum_{(\overline{\frakc},\overline{\frakd})\in \overline{\Gamma}_j} f(\overline{\fraka},\overline{\frakb},\overline{\frakc},\overline{\frakd})  \overline{\mu}^{x_i}(\overline{\fraka},\overline{\frakb})\overline{\mu}^{x_j}(\overline{\frakc},\overline{\frakd}).
\end{equation}
By \eqref{green-1},
\begin{equation}\label{eq:redu3}
G_{\mathbb{D}}(0,v)G_{\mathbb{D}}(v,w)\simeq G_{\mathbb{D}}(0,y) G_{\mathbb{D}}(y,z).
\end{equation}
Combining \eqref{eq:redu1}, \eqref{eq:redu2} and \eqref{eq:redu3}, we have 
\begin{equation*}
\begin{split}
\frac{P \Big( 0 \xrightarrow{\gamma} v \xrightarrow{\gamma} w  \Big)}{P \Big( 0 \xrightarrow{\gamma} B_{i} \xrightarrow{\gamma} B_{j} \Big)}
\overset{\eqref{eq:H1}}{\underset{\eqref{eq:Ht1}}{\simeq}} \;&\frac{\widetilde\mP(\widetilde{\cal G}^1)\widetilde\mP(\widetilde{\cal G}^2)\sum_{(\overline{\fraka},\overline{\frakb})\in \overline{\Gamma}_i} \sum_{(\overline{\frakc},\overline{\frakd})\in \overline{\Gamma}_j} f(\overline{\fraka},\overline{\frakb},\overline{\frakc},\overline{\frakd}) \overline{\mu}^{x_i}(\overline{\fraka},\overline{\frakb})\overline{\mu}^{x_j}(\overline{\frakc},\overline{\frakd})}{\sum_{y\in\partial B_i} \sum_{z\in\partial B_j }\mP({\cal G}^1)\mP({\cal G}^2)\sum_{(\overline{\fraka},\overline{\frakb})\in \overline{\Gamma}_i} \sum_{(\overline{\frakc},\overline{\frakd})\in \overline{\Gamma}_j} f(\overline{\fraka},\overline{\frakb},\overline{\frakc},\overline{\frakd}) \overline{\mu}^{x_i}(\overline{\fraka},\overline{\frakb})\overline{\mu}^{x_j}(\overline{\frakc},\overline{\frakd})}
\\
\simeq\;\;\;&\frac{\widetilde\mP\big(\widetilde{\cal G}^1(v)\big)}{\sum_{y\in\partial B_i}\mP\big({\cal G}^1(y)\big)}\frac{\widetilde\mP\big(\widetilde{\cal G}^2(w)\big)}{\sum_{z\in\partial B_j }\mP\big({\cal G}^2(z)\big)}.
\end{split}
\end{equation*}
Recall the definition of 
${\cal G}^1_\infty$ and ${\cal G}^2_\infty$ in \eqref{eq:G1infdef}, \eqref{eq:G2infdef} and of $\widetilde{\cal G}_\infty$ in \eqref{eq:Ginfdef}. By \eqref{eq:G12asymp}, \eqref{eq:Gtinfasymp} and \eqref{eq:Gt12infasymp}, we have
\begin{equation}\label{eq:finaldecouple}
\frac{\widetilde\mP\big(\widetilde{\cal G}^1(v)\big)}{\sum_{y\in\partial B_i}\mP\big({\cal G}^1(y)\big)}\frac{\widetilde\mP\big(\widetilde{\cal G}^2(w)\big)}{\sum_{z\in\partial B_j }\mP\big({\cal G}^2(z)\big)}\simeq\frac{\widetilde\mP_\infty\big(\widetilde{\cal G}_\infty\big)^2}{\sum_{y\in\partial B_i}\mP_\infty\big({\cal G}^1_\infty(y)\big)\sum_{z\in\partial B_j }\mP_\infty({\cal G}^2_\infty(z)\big)}.
\end{equation}
Recall the definition of the reference cube $B_0$, its center $x_0$ at the beginning of Section \ref{sec:sketch}. By an argument which is similar to what we have done in this section to derive \eqref{eq:finaldecouple} but essentially simpler as well as \eqref{TRANSINVA},
we have that 
\begin{equation}\label{eq:translation}
\frac{\widetilde\mP_\infty(\widetilde{\cal G}_\infty)}{\sum_{y\in\partial B_i}\mP_\infty({\cal G}^1_\infty)}\simeq \frac{P(0\xrightarrow{\gamma} x_0)}{P(0\xrightarrow{\gamma} B_0)}\simeq \alpha_0.
\end{equation}
The case for $B_j$ can be estimated similarly.
Summing over $v\in B_i$ and $w\in B_j$ and exchange the role of $i$ and $j$, we arrive at 
\begin{equation*}
\sum_{v\in B_i}\sum_{w\in B_j}P \Big( v,w\in\gamma \Big) \simeq \alpha_0^2 P( \gamma \cap B_i \neq \emptyset,\; \gamma \cap B_j\neq \emptyset)
\end{equation*}
which implies that
$$
E\big[X_i X_j \;\big|\; Y_i=Y_j=1\big] \simeq \alpha_0^2 
$$
as desired.
\end{proof}

\begin{rem} Note that, although not needed in this work, an adaptation of the argument for Proposition \ref{prop:goodpair} also gives the following estimates: for all $v\in B_i$ and $w\in B_j$, 
and for all $v,v' \in B_{i}$ and $w,w'\in B_{j}$,
\begin{equation}\label{EQ2}
P \Big( v, w \in \gamma  \Big) \simeq  P \Big( v', w' \in \gamma  \Big)
\end{equation}
for good $(i, j)$ in the sense of \eqref{goodbad}.
\end{rem}

\subsection{Proof of Proposition \ref{prop:decomp}}\label{sec:4.7}

The goal of this subsection is to give a proof of Proposition \ref{prop:decomp}, where (as we have introduced at the beginning of Section \ref{sec:decomp}) we decompose and rewrite (in an asymptotically equivalent sense) the cube-crossing event $\{ 0 \xrightarrow{\gamma} B_i \xrightarrow{\gamma} B_j \}$ into a form which is easier to analyze.To this end, in this subsection we are going to define a sequence of events (namely ${\cal A}$ through ${\cal E}$ below), each slightly modified from the last one, and show that the probability  will not change drastically (in fact we are going to show that they are all asymptotically equivalent to each other) while undergoing such a metamorphosis.

Recall the notation (e.g. the SRW $S$, and its loop-erasure $\gamma$) at the beginning of Section \ref{sec:3}.

Let $T^{i}$ (resp.\ $T^{j}$) be the {\bf last} time (up to $\Diamond$) that the simple random walk $S$ visits $B_i$ (resp.\ $B_j$).
We now write
$$
\gamma^i=\LE\big(S[0,T^i]\big)\quad\mbox{ and }\quad \gamma^j=\LE\big(S[0,T^j]\big)$$
for the loop-erasure of various segments of $S$ and define
\begin{equation}\label{eq:uidef}
u_i=\min\big\{t \geq 0 \ \big| \ \gamma^i (t) \in B_i\big\}\quad\mbox{ and}\quad u_j=\min\big\{t \geq 0 \ \big| \ \gamma^j (t) \in B_j\big\}.
\end{equation}
We also mention that as a convention in this subsection, we will use lower-case letters to denote random times with respect to the loop-erasures, and upper-case letters for those of the random walks.

 Before diving into the proof, which will be postponed till the end of this subsection, we first point out an observation. By definition of the loop-erasing procedure, it follows that 
\begin{equation}
\Big\{\gamma \cap B_i \neq \emptyset\Big\}
\Longleftrightarrow 
\Big\{\gamma^i [0, u_{i} ] \cap S[T^{i}+1, T] = \emptyset\Big\}. \label{M6}
\end{equation}

We now begin with the metamorphosis of the event $\{ 0 \xrightarrow{\gamma} B_i \xrightarrow{\gamma} B_j \}$. Note although that this event only fixes the order that the LERW $\gamma$ hits cubes $B_i$ and $B_j$ (while does not indicate the order that the SRW $S$ hits these cubes), it is very unlikely that $\gamma$ hits $B_i$ before $B_j$ while $T^j<T^i$ in the meantime, as we are going to show now in the next proposition. 
\begin{prop}\label{lem:3}
\begin{equation}\label{M9}
P\Big({0\overset{\gamma}{\rightarrow}B_i \overset{\gamma}{\rightarrow} B_j}\Big)\simeq P\big({\cal A}\big),
\end{equation}
where
\begin{equation}\label{eq:Odef}
{\cal A}:=\big\{0\overset{\gamma}{\rightarrow}B_i \overset{\gamma}{\rightarrow} B_j\big\}\cap \big\{T^i< T^j\big\}.
\end{equation}
\end{prop}

\begin{proof} 
On the event $\{0\overset{\gamma}{\rightarrow}B_i \overset{\gamma}{\rightarrow} B_j\}$, we claim that 
\begin{equation}\label{eq:hitclaim}
\gamma^j[0,u_j]\cap B_i\neq\emptyset.
\end{equation}
To see this, we argue as follows. Since $\gamma$ hits $B_j$, by \eqref{M6}, we have
$
\gamma^j[0,u_j]\cap S[T^j+1,\Diamond]=\emptyset,
$
which in turn implies that
$$
\gamma[0,m_j]=\gamma^j[0,u_j],
$$
where we write
$$
m_j=\min\big\{t \geq 0 \ \big| \ \gamma (t) \in B_j\big\}.
$$
By definition of the event $\{0\overset{\gamma}{\rightarrow}B_i \overset{\gamma}{\rightarrow} B_j\}$,  $\gamma$ must first hit $B_i$ before hitting $B_j$ at time $m_j$. So we have \eqref{eq:hitclaim} as desired.

Now assume that on the contrary $\{0\overset{\gamma}{\rightarrow}B_i \overset{\gamma}{\rightarrow} B_j\}$ and $\{T^j< T^i\}$ happen simultaneously.  Since in this case
$
S[T^j,\Diamond]\cap B_i\neq\emptyset,
$
we  have 
\begin{equation}\label{eq:H}
P\Big({0\overset{\gamma}{\rightarrow}B_i \overset{\gamma}{\rightarrow} B_j}\mbox{ and $T^j< T^i$}\big)\leq P\left( {\cal A}^{\dagger}
\right),
\end{equation}
where
$$
{\cal A}^{\dagger}:=\Big\{
\gamma^j[0,u_j]\cap S[T^j+1,\Diamond]=\emptyset;\;
\gamma^j[0,u_j]\cap B_i\neq\emptyset;\;
S[T^j,\Diamond]\cap B_i\neq \emptyset
\Big\}.
$$

We now show that 
\begin{equation}\label{eq:Odagger}
P({\cal A}^\dagger)=O(2^{-k^4+k^2}) P\big({0\overset{\gamma}{\rightarrow}B_i \overset{\gamma}{\rightarrow} B_j}\big).
\end{equation}
Writing
$$
T^-=\min\{ t\geq T^j \ \big| \ S(t) \in \partial B(x_j, l/4)\}\quad\mbox{and}\quad T^+=\max\{ t\geq T^j \ \big| \ S(t) \in \partial B(x_j, 4l)\},
$$
it is easy to see that
$$
{\cal A}^\dagger \subset \Big\{
\gamma^j[0,u_j]\cap \big(S[T^j+1,T^-]\cup S[T^+,\Diamond]\big)=\emptyset;\;
\gamma^j[0,u_j]\cap B_i\neq\emptyset\Big\}\cap \big\{ S[T^-,T^+]\cap B_i\neq \emptyset\big\},
$$
and the second event is independent from the first one conditioned on the location of $S[T^-]$ and $S[T^+]$.
WLOG assume $B\subset \frac{2}{3}\mathbb{D}$. By an decomposition similar to that in the proof of Proposition \ref{prop:bound-box-hit}, invoking the asymptotic independence of LERW,  we can show that the probability of the first event is bounded above by  $O\left( \left(\frac{\epsilon}{r}\right)^{{3-\beta}}  \left(\frac{\epsilon}{l}\right)^{{3-\beta}} \right)$ (which is of the same order as $P\big({0\overset{\gamma}{\rightarrow}B_i \overset{\gamma}{\rightarrow} B_j}\big)$ by \eqref{M7}), while that of the second event is bounded above by $O(\epsilon/l)=O(2^{-k^4+k^2})$, uniformly for all possible locations of $S[T^-]$ and $S[T^+]$.  The same argument works for the case $B \subset \mathbb{D}\setminus\frac{1}{3}\mathbb{D}$. This implies \eqref{eq:Odagger} as well as \eqref{M9}.
\end{proof}

With this proposition in mind, under the event $\cal A$,  we can chop $S$ into three segments $S[0,T^i]$, $S[T^i,T^j]$ and $S[T^j,\Diamond]$, and use the non-intersection event of these segments and their loop-erasures to represent the event ${\cal A}$. We then perform various ``surgeries'' to make this event fit the definition of $\cal F$ and show that these surgeries do not drastically change the probability.

Our next goal is to ``decouple'' $S[0,T^i]$ and $S[T^i,T^j]$. 
In other words, writing
\begin{equation}\label{eq:gammaij}
\quad\gamma^{i,j}=\LE\big(S[T^i,T^j]\big),
\end{equation}
we would like to show that in the events we consider, roughly speaking, $\gamma^j=\LE(S[0,T^j])$ can be replaced by the concatenation of $\gamma^i$ and $\gamma^{i,j}$. We write
\begin{equation}\label{eq:uijdef}
u_{i,j}=\min\big\{t \geq 0 \ \big| \ \gamma^{i,j} (t) \in B_j\big\},
\end{equation}
and define
\begin{equation}\label{eq:Jdef}
{\cal B}:=\Big\{T^i< T^j,\; \gamma^i[0,u_i]\cap S[T^i+1,T]=\emptyset,\ \gamma^{i,j}[0,u_{i,j}]\cap S[T^j+1,T]=\emptyset\Big\}.
\end{equation}
\begin{prop}\label{prop:lem5}
\begin{equation}\label{eq:lem5}
P\big({\cal A}\big) 
\simeq P\big({\cal B}\big).
\end{equation}
\end{prop}
Before proceeding to the proof, we need a technical estimate. The next lemma shows that conditioned on $\cal A$, 
with high probability, the LERW $\gamma^i(u_i,\Diamond]$ after the first hitting to $B_i$ has a small diameter, thus cannot form a large ``quasi-loop'' introduced in Section \ref{qloops}. We omit the proof for 
an essentially similar estimate has been proved in \cite{S2} for a two-cube-crossing event (see the arguments leading to (3.70) ibid.).
\begin{lem}\label{lem:4}
Let
\begin{equation}
{\cal A}':= {\cal A} \cap 
\Big\{ \gamma^i[u_i,\Diamond]\subset B\big(x_i, 2^{-k^4+k^2}\big)\Big\}.
\end{equation}
Then, there exists $c>0$ such that
\begin{equation}\label{eq:AAprime}
P\big({\cal A}\setminus{\cal A}'\big)=O\big(2^{-ck^2}\big)P\big({\cal A}\big).
\end{equation}
\end{lem}

\begin{proof}[Proof of Proposition \ref{prop:lem5}]
We only prove the $\leq$ direction as the opposition direction can be proved in a completely similar manner.
In order to prove the $\leq$ direction of \eqref{eq:lem5}, by \eqref{eq:AAprime} it suffices to show 
\begin{equation}\label{eq:J_0J}
P\big({\cal A}'\big) \leq  P\big({\cal B}\big)\big(1+O(2^{-ck^2})\big).
\end{equation}

We now define two random times for $S$:
$$
 \Psi^-=\max\big\{ t\geq T^i \ \big| \ S(t)\in B(x_i, 2^{-k^4+k^2})\big\}
\mbox{ 
and 
 }
\Psi^+=\max\big\{ t\geq T^i\ \big| \ S(t)\in B(x_i, 2^{-k^4+2k^2})\big\},
$$
and write
$$
{\cal A}'':={\cal A}' \cap \Big\{S[T^i, \Psi^-]\subset B\big(x_i, 2^{-k^4+2k^2}\big)\Big\} \cup \Big\{S[T^i, \Psi^+]\subset B\big(x_i, 2^{-k^4+3k^2}\big)\Big\}.
$$
In a similar fashion as Lemma \ref{lem:4} we can show that there exists $c>0$ such that 
\begin{equation}\label{eq:ApApp}
P\big({\cal A}'\setminus {\cal A}''  \big) = O(2^{-ck^2}) P\big({\cal A}'\big).
\end{equation}
If ${\cal A}''$ occurs, then it follows that 
$$
\gamma^{i,j} [0,u_{i,j} ] \cap B\big(x_i, 2^{-k^4+3k^2}\big)^c \subset \gamma^j[0,u_j].
$$
Therefore, we have
\begin{equation}\label{eq:AppBp}
{\cal A}'' \subset {\cal B}':=\left\{ T^i<T^j;\; \gamma^i [0,u_i]\cap S[T^i+1,T]=\emptyset;\;
\gamma^{i,j} [0,u_{i,j}]\cap B\big(x_i, 2^{-k^4+3k^2}\big)^c \cap S[T^j+1 ,\Diamond]=\emptyset\right\}.
\end{equation}
Note that 
$$
{\cal B}'\setminus {\cal B} \subset {\cal B}' \cap \Big\{S [T^j, \Diamond] \cap B \big(x_i, 2^{-k^4+3k^2}\big) \neq \emptyset\Big\}.
$$
By an argument similar to the proof of Proposition \ref{lem:3}, 
 we have
\begin{equation}\label{eq:BpB}
P\big({\cal B}'\big) \leq P\big({\cal B}\big) \big(1+O\big(2^{-k^4+4k^2}\big)\big).
\end{equation} 

Combining \eqref{eq:ApApp}, \eqref{eq:AppBp} and \eqref{eq:BpB}, we obtain \eqref{eq:J_0J} as desired.
\end{proof}

We now perform the next surgery towards $\cal F$. Recall the definition of $u_i$ and $u_{i,j}$ in \eqref{eq:uidef} and \eqref{eq:uijdef}. We define\begin{equation}\label{eq:tidef}
s_i:=\max\big\{t \ \big| \ \gamma^i(t)\in\partial B(x_i, 2^{-k^4+k^3})\big\}\mbox{ and }
u_i^+:= \min \big\{ t\geq s_i \ \big| \; \gamma^i(t)\in B_i\big\},
\end{equation}
and define $s_{i,j}$, ${u}_{i,j}^+$ for $\gamma^{i,j}$ similarly. These random times correspond to $U_{i,-k^4+k^3}(\gamma^1)$, $V_1$, $U_{j,-k^4+k^3}(\gamma^2)$, $V_2$ (see \eqref{eq:Udef} and \eqref{eq:V1V2def} for their definition; note that $\gamma^i$ and $\gamma^1$ travel in ``opposite'' directions and same for  $\gamma^{i,j}$ and $\gamma^2$) respectively in the definition of $\cal F$.  With the definition of ${\cal B}$ (see \eqref{eq:Jdef}) in mind, we now define
\begin{equation}\label{eq:kdef}
{\cal C}:=\left\{ T^i<T^j;\;\gamma^i [0,u_i^+]\cap S[T^i+1,\Diamond]=\emptyset;\;
\gamma^{i,j} [0,{u}_{i,j}^+] \cap S[T^j+1,\Diamond]=\emptyset\right\}.
\end{equation}
with $u_i^+$ and ${u}_{i,j}^+$ replacing $u_i$ and $u_{i,j}$ in $\cal B$. Note that $u_i\leq u_i^+$, and $u_i<u_i^+$ if and only if the loop-erasure $\gamma^i[u_i,\Diamond] \cap \partial B \big(x_i, 2^{-k^4+k^3}\big) \neq \emptyset$. The same observation remains true for $u_{i,j}$ and ${u}_{i,j}^+$ in the same fashion.
By an argument similar to Lemma \ref{lem:4}, we can show actually
$u_i=u_i^+$ and $u_{i,j}={u}_{i,j}^+$ with high probability. More precisely, we have the following lemma. 
\begin{lem}\label{lem:6}
\begin{equation}\label{eq:BtoC}
P\big({\cal B}\big)\simeq P\big({\cal C}\big).
\end{equation}
\end{lem}
We omit the proof as it is almost the same as Lemma \ref{lem:4}.

Out next goal is to cut out the last segment of $S[T^i,T^j]$ in the definition of $\cal C$ (see \eqref{eq:kdef}). More precisely,  
let 
\begin{equation}\label{eq:Sigmadef}
 \Xi: = T^i + U_{j,-k^3+k^2}(S[T^i+\cdot])
\end{equation}
be the first time after $T^i$ the walk $S$ hits $\partial B\big(x_j,2^{-k^3+k^2}\big)$
(note that this corresponds to $U_{j,-k^3+k^2}$ of $(X^2)^{\cal R}$ in the definition of ${\cal F}$ in \eqref{eq:Hdef}).
By replacing $S[T^i+1,\Diamond]$ by $S[T^i+1, \Xi ]\cup S[T^j,\Diamond]$ in the definition of $\cal C$ in \eqref{eq:kdef}, we now define
\begin{equation}\label{eq:ldef}
{\cal D}:=\left\{
T^i<T^j;\;\gamma^i [0,u_i^+]\cap \big(S[T^i+1, \Xi ]\cup S[T^j,\Diamond]\big)=\emptyset;\;
\gamma^{i,j} [0,{u}_{i,j}^+] \cap S[T^j+1,\Diamond]=\emptyset
\right\}.
\end{equation}
We now show that this surgery will not change the probability very much.

\begin{lem} \label{lem:7prime}

\begin{equation}\label{eq:7prime}
 P({\cal C})\simeq P \left({\cal D}\right).
\end{equation}
\end{lem}

\begin{proof}
Note that by definition, ${\cal C}\subset {\cal D}$. Hence to prove \eqref{eq:7prime},
it suffices to obtain an upper bound for the probability of 
$$
{\cal D}^\dagger:= {\cal D}\setminus {\cal C}= {\cal D}  \cap  \Big\{
\gamma^i[0,u_i^+]\cap S[ \Xi ,T^j]\neq\emptyset \Big\}.
$$
We now decompose ${\cal D}^\dagger$ according to at which scale $\gamma^i[0,u_i^+]$ and $S[ \Xi ,T^j]\neq\emptyset$ meet.  

Recall that $l=|x_i-x_j|$. We first deal with 
$${\cal D}^\dagger_1:={\cal D} \cap \Big\{ \gamma^i [0,u_i^+]\cap B(x_j, 2^{-k^2}l)\neq \emptyset\Big\}.$$
WLOG assume $B\subset \frac{2}{3}\mathbb{D}$. A similar argument as in the proof of Proposition \ref{lem:3} gives that:
\begin{equation}
\begin{split}
P[{\cal D}^\dagger_1]&\leq C \frac{2^{-k^2}l}{r} \Es_n (2^{-k^2}l, l)\cdot  \frac{\epsilon}{l}\cdot\Es_n(\epsilon,l) \cdot \frac{\epsilon}{l}\cdot\Es_n(\epsilon,l)\Es_n(l,r)\\
& \leq C' \left(\frac{\epsilon}{r}\right)^{3-\beta}\left(\frac{\epsilon}{l}\right)^{3-\beta} 2^{(\beta-3)k^2}\leq C''2^{(\beta-3)k^2}P[{\cal C}].
\end{split}
\end{equation}
Therefore, we only need to estimate the probability of 
$${\cal D}^\dagger_2:={\cal D}^\dagger \cap \big\{ {\rm dist}\big(\gamma^i[0,u_i^+],x_j\big)\geq 2^{-k^2} l\big\}.$$
With this in mind, for $q=0,1,\ldots,k^2$,  let 
$$
{\cal M}^{q}:= 
\left\{
 2^{-k^2+q}\cdot l  \; < {\rm dist}\Big(\gamma^i[0, u_i^+],x_j\Big)\leq 2^{-k^2+q+1}\cdot l
 \right\}.
$$
Observe that 
$$
\big({\cal D}^\dagger\cap {\cal M}^q\big)  \subset {\cal L}^q := {\cal D} \cap \left\{ \begin{split}
&\gamma^i [0, u_i^+] \cap B(x_j ,2^{-k^2}\cdot l \cdot 2^{q+1})\neq \emptyset \\
& S[ \Xi,T^j] \cap  B(x_j, 2^{-k^2}\cdot l \cdot 2^{q}) \neq \emptyset.
\end{split}\right\}
$$
As in the proof of Proposition \ref{lem:3} (note that in this case we need to truncate according to ${\rm dist}\big(x_j,S[\Xi, Y]\big)$, where $Y$ is the first time after $\Xi$ that $S$ hits $\partial B(x_j, 2^{-k^2+q}\cdot l)$), we have
\begin{equation*}
P({\cal L}^q)\leq C\cdot \frac{2^{-k^2+q}\cdot l}{r} \cdot \Es_n(2^{-k^2+q}\cdot l,l) \cdot \frac{\epsilon}{l}\cdot \Es_n(\epsilon,l) \cdot\frac{\epsilon}{l}\cdot \frac{2^{-k^3+k^2}}{2^{-k^2+q}\cdot l}\cdot \Es_n(\epsilon,2^{-k^3+k^2})\Es_n(l,r).
\end{equation*}
After simplification, we know that for some $c,C,C'>0$,
$$
P({\cal L}^q)\leq C\cdot \Big(\frac{\epsilon}{r}\Big)^{2-\beta} \Big(\frac{\epsilon}{l}\Big)^{2-\beta} 2^{-c(k^3+q)}\leq C'   2^{-c(k^3+q)} P[{\cal C}].
$$
Hence,
\begin{equation}
P[{\cal D}^\dagger]\leq P[{\cal D}^\dagger_1]+P[{\cal D}^\dagger_2]\leq C 2^{(\beta-3)k^2}P[{\cal C}] +\sum_{q=0}^{k^2} C'2^{-c(k^3+q)} P[{\cal C}]\leq C'' 2^{-c'k^2} P[{\cal C}].
\end{equation}
This finishes the proof of \eqref{eq:7prime}.
\end{proof}

Before proving Proposition \ref{prop:decomp}, we need to perform one final surgery. Let 
\begin{equation}
\Sigma:=T^i+U_{i,-k^3}(S[T^i+\cdot])
\end{equation} be the first time after $T^i$ the walk $S$ hits $\partial B(x_i, 2^{-k^3})$ (note that this corresponds to the definition of $\sigma$ in \eqref{eq:sigmadef} in setting up for the event $\cal F$). Define 
$$
\overline{\gamma}^{i,j}=\LE \big(S[\Sigma,T^j]\big)
$$
and define $\overline{u}_{i,j}^+$ in a similar fashion as ${u}_{i,j}^+$ for $\gamma^{i,j}$ in \eqref{eq:tidef}.
We now replace $\gamma^{i,j}$ by $\overline{\gamma}^{i,j}$ and ${u}_{i,j}^+$ by $\overline{u}_{i,j}^+$ in the definition of $\cal D$ (see \eqref{eq:ldef}) and define the following event 
\begin{equation}
{\cal E}:=\left\{
\begin{split}
&T^i<T^j;\quad \gamma^i [0,u^+_i]\cap \big(S[T^i+1, \Xi ]\cup S[T^j,T]\big)=\emptyset;\\
&\overline{\gamma}^{i,j} [0,\overline{u}_{i,j}^+]\cap  S[T^j+1,T]=\emptyset
\end{split}
\right\}.
\end{equation}
By an argument very similar to the proof of Lemma \ref{lem:7prime},  we can show that
\begin{lem}
\begin{equation}\label{eq:LN}
P[{\cal E}]\simeq P[{\cal D}].
\end{equation}
\end{lem}
\begin{proof}[Sketch of Proof]
On ${\cal D} \setminus {\cal E}$, $\overline{\gamma}^{i,j}$ must intersect $S[T^j+1,\Diamond]$ before it merges with $\gamma^{i,j}$. Similar to the proof of Lemma \ref{lem:7prime}, one can further decompose  ${\cal D} \setminus {\cal E}$ according to the scale at which they meet and show that $P[{\cal D} \setminus {\cal E}]=O(2^{-ck^2})P({\cal D})$. The same argument works for ${\cal E} \setminus {\cal D}$.
\end{proof}
Finally, we arrive at Proposition \ref{prop:decomp}. We remark that the last-exit decomposition argument below bears some similarity to (3.19) of \cite{S2}.
\begin{proof}[Proof of Proposition \ref{prop:decomp}]
By \eqref{M9}, \eqref{eq:lem5}, \eqref{eq:BtoC}, \eqref{eq:7prime} and \eqref{eq:LN},
$$
P\Big({0\overset{\gamma}{\rightarrow}B_i \overset{\gamma}{\rightarrow} B_j}\Big)\simeq P({\cal A})\simeq P({\cal B})\simeq P({\cal C})\simeq P({\cal D}) \simeq P({\cal E}).
$$
Decompositions at $T^i$ and sub-sequently at $T_j$ for the random walk $S$ restricted on the event ${\cal E}$ allow us to treat $S[0,T^i]$ $S[T^i,T^j]$ and $S[T^j,\Diamond]$ as three independent walks conditioned on the location of $S[T^i]$ and $S[T^j]$. This along with the reversibility of LERW gives that (in the notation of Proposition \ref{prop:decomp})
$$
P({\cal E})=\sum_{y\in\partial_i B_i} \sum_{z\in\partial_i B_j} G_{\mathbb{D}} (0,y) G_{\mathbb{D}} (y,z) \mP^{y,z,z}\big({\cal F}(y,z)\big)
$$
as desired. \end{proof}

\subsection{Proof of Proposition 
\ref{prop:Lem11} (and of Lemma \ref{lem:FFocompare})}\label{sec:endingpart}

In this subsection, we will give the proof of Proposition \ref{prop:Lem11} and sketch the proof of Lemma \ref{lem:FFocompare}. 

We now briefly discuss the strategy for proving Proposition \ref{prop:Lem11}. Recall various notation and conventions from Section \ref{sec:decoupling} (in particular the collection of pairs of paths $\Delta_i$, $\Delta_j$, the truncation mapping $\pi(\cdot,\cdot,\cdot,\cdot)$, the random walks $R^{\fraka,0}$, $R^{\frakc,\frakb}$, $R^\frakd$ and their joint law $\mathbb{Q}$, etc.). We first define another function $h$ to measure the ``typicality'' of the quadruple of paths $(\fraka,\frakb,\frakc,\frakd)$, which operates in a similar fashion as $f$, but much easier to be decoupled into products of two functions measuring the typicality of $(\fraka,\frakb)$ and $(\frakc,\frakd)$ separately (see e.g.\ Lemma \ref{lem:hdecomp} and the definition of typical paths in \eqref{eq:Cijdef}). Then, as a crucial step in the proof of this proposition, we will show in Proposition \ref{prop:Lem9} that for ``typical'' quadruples, (see \eqref{eq:Cijdef} for a precise definition),  $f(\fraka,\frakb,\frakc,\frakd)$ almost  depends on the end parts $\pi(\fraka,\frakb,\frakc,\frakd)$ only. This follows from an analysis in the proof to preclude rare ``backtrackings'' of random walk.

We start with the definition of $h$. For $(\fraka,\frakb) \in \Delta_i$ and $(\frakc,\frakd)\in \Delta_j$, recall the definition of $V_1$ and $V_2$ from \eqref{eq:V1V2def} and write
\begin{equation}\label{eq:hdef}
h(\fraka,\frakb,\frakc,\frakd)=\mathbb{Q}\left(
\begin{split}
& \big(\fraka[V_1,\Diamond]\cup \LE(R^{\fraka,0})[0,U_{i,-k^3+1}]\big)\bigcap \big(R^{\frakb,\frakc}[0,U_{i,-k^3+1}]\cup \frakb\cup R^\frakd\cup \frakd\big)=\emptyset, \\
& \big(\LE(R^{\frakc,\frakb})[0,U_{j,{-k^3+1}}] \cup \frakc[V_2,\Diamond]\big)\bigcap \big(R^\frakd[0,U_{j,{-k^3+1}}]\cup \frakd\big)=\emptyset 
\end{split}
\right).
\end{equation}
Note that in contrast to the definition of $f$ in \eqref{eq:fdef}, here we only require non-disconnection up to scale $2^{-k^3+1}$.
One may compare our definition here with the definition of $h$ for the one-point estimate case in in Section 3.3 of \cite{Escape}. See the discussion in Remark 3.6, ibid., for the significance of such functions.

For $s,t>0$, write \begin{equation}
K(k,s,t)=\left(2^{-k^3}/s \times 2^{-k^3}/t\right)^{2-\beta}
\end{equation}
for short. 
Recall the definition of $r$ and $l$ in \eqref{eq:rDEF} and \eqref{goodbad} respectively. We now claim that for the same pairs of paths, the ratio between $f$ and $h$ is up-to-constants equivalent to the explicit quantity $[\Es_n(2^{-k^3},l)]^2\Es_n(l,r)$ (which is also equivalent to $K(k,r,l)$ in the notation above), explaining heuristically the effect of extra non-disconnection in the definition of $f$ compared to $h$. 
\begin{lem} \label{lem:Lem8}
For $(\fraka,\frakb) \in \Delta_i$ and $(\frakc,\frakd)\in \Delta_j$, it follows that 
$$
f(\fraka,\frakb,\frakc,\frakd)\asymp K(k,r,l) h(\fraka,\frakb,\frakc,\frakd).
$$
\end{lem}
This lemma is essentially a two-point version of Proposition 3.5 of \cite{Escape}, where we compare asymptotics of similar functions on the non-intersection probability of one pair of paths only, which is a rather standard application of asymptotic independence of LERW (i.e., quasi-multiplicativity of non-intersection probabilities; see Section \ref{sec:basicproperties}), Harnack principle and the separation lemma (for the lower bound). 


 Recall that $l=|x_i-x_j|$ as defined in \eqref{goodbad}. To decouple the analysis of the behavior of paths up to $l/4$ near $B_i$ and $B_j$, we define for this subsection only
\begin{equation}\label{eq:truncatedl3}
\gamma^1:= \LE(R^{\fraka,0})[0,U_{i,p(l)}],\quad \lambda^1:=R^{\frakb,\frakc}[0,U_{i,p(l)}];
\quad
\gamma^2:=\LE(R^{\frakc,\frakb})[0,U_{j,p(l)}];\quad \lambda^2=R^\frakd[0,U_{j,p(l)}]
\end{equation}
where we write $p(l):= \log(l/4)$ for short. Then, we introduce the following events ${\cal L}^1={\cal L}^1(\fraka,\frakb,\frakc)$ and ${\cal L}^2={\cal L}^2(\frakb,\frakc,\frakd)$ to describe the non-intersection near $B_i$ and $B_j$ up to scale $l/4$:
\begin{equation}\label{eq:L1L2def}
{\cal L}^1:=\left\{
 \big(\fraka[V_1,\Diamond]\cup \gamma^1\big)\bigcap \big(\lambda^1\cup \frakb\big)=\emptyset 
\right\}\quad\mbox{and}
\quad
{\cal L}^2:=\left\{
 \big(\gamma^2 \cup \frakc[V_2,\Diamond]\big)\bigcap (\lambda^2\cup \frakd\big)=\emptyset 
\right\}.
\end{equation}

\begin{lem} \label{lem:hdecomp} For all $(\fraka,\frakb)\in\Delta_i$ and $(\frakc,\frakd)\in\Delta_j$,
\begin{equation}\label{eq:hdecomp}
  P \left({\cal L}^1\cap {\cal L}^2 \right)
  \asymp P \left( {\cal L}^1 \right) \times P\left( {\cal L}^2
\right) \asymp h(\fraka,\frakb,\frakc,\frakd) K(k,l,l).
\end{equation}
\end{lem}
We omit the proof for the same reason for omitting the proof of Lemma \ref{lem:Lem8}.

To define the ``typicality'' of pairs of paths, we now define some reference quantities: let $R^{\frakc,x_i}$ be conditioned random walks defined as $R^{\frakc,\frakb}$ in Definition \ref{defnf}, but with $\frakb$ replaced by $x_i$. We also define $R^{\frakb,x_j}$ in a similar fashion. With slight abuse of notation we still denote by $P$ the corresponding probability measure.
We then define ${\cal L}^1_0={\cal L}^1_0(\fraka,\frakb)$ and ${\cal L}^2_0={\cal L}^2_0(\frakc,\frakd)$ with $R^{\frakb,x_j}$ and $R^{\frakc,x_i}$ replacing the roles of $R^{\frakb,\frakc}$ and $R^{\frakc,\frakb}$ respectively. 
Note that by \eqref{eq:asymp-variant} and Harnack principle, for all $(\fraka,\frakb)\in \Delta_i$ and $(\frakc,\frakd)\in \Delta_j$,
$$
P({\cal L}^1(\fraka,\frakb,\frakc))\simeq P({\cal L}^1_0(\fraka,\frakb))\quad\mbox{ and }\quad P({\cal L}^2(\frakb,\frakc,\frakd))\simeq P({\cal L}^2_0(\frakc,\frakd)).
$$
For pairs of paths $(\fraka,\frakb)\in \Delta_i$ and $(\frakc,\frakd)\in \Delta_j$, 
we write 
\begin{equation}\label{eq:h1h2def}
h^1(\fraka,\frakb)=P({\cal L}^1_0)\mbox{ and }h^2(\frakc,\frakd)=P({\cal L}^2_0).
\end{equation}
Recall that $\pi$ is the truncation operator defined in and below \eqref{eq:pidef}. We now set
\begin{equation}\label{eq:Cijdef}
{\cal C}_i=\big\{(\fraka,\frakb)\in\Delta_i\big| h^1(\pi(\fraka,\frakb))\geq 2^{-k^3}\big\}
\mbox{ and }
{\cal C}_j=\big\{(\frakc,\frakd)\in\Delta_j\big| h^2(\pi(\frakc,\frakd))\geq 2^{-k^3}\big\}.
\end{equation}
for the collection of pairs of ``typical'' paths. As discussed at the beginning of this subsection, we claim that for typical paths, the value of $f(\fraka,\frakb,\frakc,\frakd)$ is very close to $f(\pi(\fraka,\frakb,\frakc,\frakd))$.
\begin{prop}\label{prop:Lem9}
For $(\fraka,\frakb)\in {\cal C}_i$ and $(\frakc,\frakd)\in{\cal C}_j$, 
 it follows that 
\begin{equation}\label{eq:ffpi}
\big| f(\fraka,\frakb,\frakc,\frakd)-f(\pi(\fraka,\frakb,\frakc,\frakd))\big| \leq C2^{-k^3} f(\pi(\fraka,\frakb,\frakc,\frakd)).
\end{equation}
\end{prop}
Note that the function $f$ is defined as a conditional probability. In order to analyze the difference of $f(\cdot)$ and $f(\pi(\cdot))$, we need to rewrite it into an unconditioned form. To this end, for $(\fraka,\frakb)\in\Delta_i$ and $(\frakc,\frakd)\in\Delta_j$, we define three new independent RW's $R^1$, $R^2$, $R^3$ such that
\begin{itemize}
\item $R^1$ starts from $x_\fraka$, is conditioned to hit $0$ before leaving $\mathbb{D}$ and stops upon hitting $0$.
\item $R^2$ starts from $x_\frakc$, is conditioned to hit $x_\frakb$ before leaving $\mathbb{D}$, and is truncated at the last time it visits $x_\frakb$. 
\item $R^3$: starts from $x_\frakd$ stopped upon leaving $\mathbb{D}$.
\end{itemize}
(Recall that $x_\bullet$ stands for the ending point of $\bullet$ for $\bullet=\fraka,\frakb,\frakc,\frakd$.) We write $\mP_0$ for their joint law. Note that the definitions of $R^3$ and $R^\frakd$ are identical while the difference between $R^1$, $R^2$ and $R^{\fraka,0}$ and $R^{\frakc,\frakb}$ is that here we do not require $R^1$ and $R^2$  to avoid $\fraka$ and $\frakb$ respectively in the way $R^{\fraka,0}$ and $R^{\frakc,\frakb}$ are defined.
\begin{proof}
Define 
$$
H_1=H_1(\fraka,\frakb,\frakc):=\left\{R^1[1,\Diamond]\cap \fraka=\emptyset, R^2[1,\Diamond]\cap \frakc=\emptyset\right\},
$$
and  
\begin{equation*}
H_2=H_2(\fraka,\frakb,\frakc,\frakd):=H_1(\fraka,\frakb,\frakc)\bigcap\left\{
\begin{split}
& \big(\fraka[V_1,\Diamond]\cup \LE(R^1)\big)\bigcap \big((R^{2})^{\cal R}[0,U_{j,-k^3+k^2}]\cup \frakb\cup R^3\cup \frakd\big)=\emptyset \\
& \big(\LE(R^{2}) \cup \frakc[V_2,\Diamond]\big)\bigcap (R^3\cup \frakd\big)=\emptyset 
\end{split}
\right\}.
\end{equation*}
Note that the second event on the RHS above is just the event on the RHS of \eqref{eq:fdef} where $f(\fraka,\frakb,\frakc,\frakd)$ is defined with $R^{\fraka,0}$, $R^{\frakc,\frakb}$ and $R^\frakd$ replaced by $R^1$, $R^2$ and $R^3$ respectively, hence
\begin{equation}\label{eq:ftoH}
f(\fraka,\frakb,\frakc,\frakd)=\mP_0(H_2|H_1).
\end{equation}
We then define $H'_1$ and $H'_2$ through replacing $(\fraka,\frakb,\frakc,\frakd)$ by $\pi(\fraka,\frakb,\frakc,\frakd)$ in the definition of $H_1$ and $H_2$ respectively. In the same fashion, we have
\begin{equation}\label{eq:fptoHp}
f(\pi(\fraka,\frakb,\frakc,\frakd))=\mP_0(H'_2|H'_1).
\end{equation}
By definition, it follows that
$$
H_1\subset H'_1\mbox{ and }H_2\subset H'_2.
$$

We will first show that $\mP_0(H_1)$ is close to $\mP_0(H'_1)$. Note that 
\begin{equation}\label{eq:Hprime1Hprime}
\mP_0(H'_1)-\mP_0(H_1)\leq \mP_0\left(H'_1, R^1[1,\Diamond]\cap \fraka\setminus\pi(\fraka) \neq \emptyset\right)+\mP_0\left(H'_1, R^2[1,\Diamond]\cap \frakc\setminus\pi(\frakc)] \neq \emptyset\right).
\end{equation}
We now claim that both terms on the RHS above is bounded by $ O(2^{-9k^3}) \mP_0(H'_1)$.
We first look at the first term and focus on the behavior of $R^1$. Note that the event $\{R[1,\Diamond]\cap \pi(\fraka) =\emptyset; R[1,\Diamond]\cap \fraka\setminus\pi(\fraka) \neq \emptyset\}$ implies the intersection of the following two events $\cal S$ and $\cal T$, with ${\cal T}$ independent from $\cal S$ conditioned on the location of first exit from $B(x_\fraka,2^{-k^3-1})$ by $R^1$:
\begin{enumerate}
\item[1)]  ${\cal S}:=\left\{\mbox{$R^1$ avoids $\pi(\fraka)$ up to first exit from $B(x_\fraka,2^{-k^3-1})$}\right\}$;
\item[2)]  ${\cal T}:=\left\{\mbox{after exiting $B(x_\fraka, 2^{-k^3-1})$, $R^1$ hits $\fraka[0,U_{i,-10k^3})=\fraka\setminus\pi(\fraka)$}\right\}$.
\end{enumerate}  
We now claim that $\mP_0[{\cal S}]\asymp \mP_0[R^1[1,\Diamond]\cap \pi(\fraka) =\emptyset]$. It suffices to prove the ``$\leq$'' side. By Proposition 6.1 of \cite{S}, conditioning on $\cal S$, there is a uniformly positive probability (regardless of $\fraka$) that the exit location of $R^1$ from $B(x_\fraka,2^{-k^3-1})$ is of order $2^{-k^3-1}$ away from $\pi(\fraka)$, hence there is also a uniformly positive probability to avoid $\pi(\fraka)$ for the rest of its duration.
In the meanwhile, uniformly for $\fraka$ and the exit location, we have $\mP_0[{\cal T}]=O(2^{-9k^3})$. Taking $R^2$ into consideration, we see that the first term on the RHS of \eqref{eq:Hprime1Hprime} is bounded by $O(2^{-9k^3})\mP_0[H'_1]$.
In a similar fashion, the same bound applies for the second term. Thus,
\begin{equation}\label{eq:H1H1p}
\mP_0(H_1)\simeq \mP_0(H'_1).
\end{equation}

Next we will compare $\mP_0(H_2)$ and $\mP_0(H'_2)$ under the assumptions that $(\fraka,\frakb)\in {\cal C}_i$, $(\frakc,\frakd)\in{\cal C}_j$. By an argument similar to the analysis of \eqref{eq:Hprime1Hprime}, again with an application of Proposition 6.1 of \cite{S}, we have
\begin{equation}\label{eq:Hprime2Hprime}
 \mP_0(H'_2)-\mP_0(H_2)
\leq  \mP_0\Big(H'_2, \Big(R^1\cup R^2\cup R^3\Big)\bigcap B(\{x_i,x_j\}, 2^{-10k^3})\neq \emptyset\Big)\leq c\cdot 2^{-9k^3} \mP_0(H'_1).
\end{equation}
By applying Lemma \ref{lem:Lem8} to $\pi(\fraka,\frakb,\frakc,\frakd)$ and then invoking Lemma \ref{lem:hdecomp}, we have
$$
f(\pi(\fraka,\frakb,\frakc,\frakd))\asymp h^1(\pi(\fraka,\frakb)) \cdot h^2(\pi(\frakc,\frakd)) \cdot K(k,r,l).
$$
When $(\fraka,\frakb)\in {\cal C}_i$ and $(\frakc,\frakd)\in {\cal C}_j$, it follows that 
$$
f(\pi(\fraka,\frakb,\frakc,\frakd))\geq c\cdot 2^{-2k^3}\cdot K(k,r,l) \geq c2^{-4k^3}.
$$
Therefore, in this case
$$
\mP_0(H'_2)-\mP_0(H_2)\leq c\cdot 2^{-5k^3} \cdot P(H'_2),
$$
implying that
\begin{equation} \label{eq:H2H2p}
\mP_0(H_2)\simeq\mP_0(H'_2).
\end{equation}
Finally,  by \eqref{eq:ftoH}, \eqref{eq:fptoHp}, \eqref{eq:H1H1p} and \eqref{eq:H2H2p}, for $(\fraka,\frakb)\in {\cal C}_i$ and $(\frakc,\frakd)\in {\cal C}_j$,
\begin{equation}
f(\fraka,\frakb,\frakc,\frakd)=\mP_0(H_2|H_1)\simeq \mP_0(H'_2|H'_1)= f(\pi(\fraka,\frakb,\frakc,\frakd)),
\end{equation}
which completes the proof.
\end{proof}

Proposition \ref{prop:Lem11} now follows as an application of Lemmas \ref{lem:Lem8} and \ref{lem:hdecomp}, Proposition \ref{prop:Lem9} and the separation lemma.
\begin{proof}[Proof of Proposition \ref{prop:Lem11}]
By the separation lemma (\ref{SEPARATION-LEMMA-1}), it follows that (recall the definition of well-separatedness below \eqref{SEPARATION-LEMMA-1})
$$
\nu^y\big(\{(\fraka,\frakb)\in\Delta_i\;|\;(\fraka,\frakb)\mbox{ is well-separated}\}\big)\geq c.
$$
and 
$$
\nu^z\big(\{(\frakc,\frakd)\in\Delta_j\;|\;(\frakc,\frakd)\mbox{ is well-separated}\}\big)\geq c.
$$
If $(\fraka,\frakb)$ and $(\frakc,\frakd)$ are well-separated, it is easy to see that 
$$
h^1(\fraka,\frakb)\geq c\mbox{ and } h^2(\frakc,\frakd)\geq c.
$$
Thus, by Lemmas \ref{lem:Lem8} and \ref{lem:hdecomp} along with \eqref{eq:h1h2def} which altogether relate $f$ to $h^1$ and $h^2$,
$$
f(\fraka,\frakb,\frakc,\frakd)\geq c \cdot K(k,r,l) \mbox{,  if $(\fraka,\frakb)$ and $(\frakc,\frakd)$ are well-separated.}
$$
Therefore we have
\begin{equation}\label{eq:Hdecomp3}
\sum_{(\fraka,\frakb)\in \Delta_i} \sum_{(\frakc,\frakd)\in \Delta_j} f(\fraka,\frakb,\frakc,\frakd) \nu^y(\fraka,\frakb)\nu^z(\frakc,\frakd) \asymp K(k,r,l).
\end{equation}
Recalling the definition of ${\cal  C}_i$ and ${\cal  C}_j$ in \eqref{eq:Cijdef} and noting that by definition $h^1(\fraka,\frakb)\leq h^1(\pi(\fraka,\frakb))$ and $h^2(\frakc,\frakd)\leq h^1(\pi(\frakc,\frakd))$, again by Lemma \ref{lem:Lem8}, \eqref{eq:hdecomp} and \eqref{eq:h1h2def}, we have if $(\fraka,\frakb)\notin {\cal C}_i$ or $(\frakc,\frakd)\notin {\cal C}_j$,
$$
f(\fraka,\frakb,\frakc,\frakd)\leq c 2^{-k^3}\cdot K(k,r,l) \mbox{,  if $(\fraka,\frakb)\notin {\cal C}_i$ or $(\frakc,\frakd)\notin {\cal C}_j$}.
$$
Hence,
\begin{equation}\label{eq:exceptional}
\sum_{(\fraka,\frakb)\notin {\cal C}_i\mbox{ or } (\frakc,\frakd)\notin {\cal C}_j} f(\fraka,\frakb,\frakc,\frakd) \nu^y(\fraka,\frakb)\nu^z(\frakc,\frakd)\leq C2^{-k^3} K(k,r,l),
\end{equation}
and by \eqref{eq:Hdecomp3},
\begin{equation}\label{eq:temp0}
\sum_{(\fraka,\frakb)\in \Delta_i} \sum_{(\frakc,\frakd)\in \Delta_j} f(\fraka,\frakb,\frakc,\frakd) \nu^y(\fraka,\frakb)\nu^z(\frakc,\frakd) \simeq \sum_{(\fraka,\frakb)\in {\cal C}_i}\sum_{(\frakc,\frakd)\in {\cal C}_j} f(\fraka,\frakb,\frakc,\frakd) \nu^y(\fraka,\frakb)\nu^z(\frakc,\frakd).
\end{equation}
Finally,
\begin{equation*}
\begin{split}
\mP({\cal F} \big| {\cal G}^1\cap {\cal G}^2)\overset{\eqref{eq:Hdecomp2}}{\simeq}\,&\sum_{(\fraka,\frakb)\in \Delta_i} \sum_{(\frakc,\frakd)\in \Delta_j} f(\fraka,\frakb,\frakc,\frakd) \nu^y(\fraka,\frakb)\nu^z(\frakc,\frakd)  \\
\overset{\eqref{eq:temp0}}{\simeq}& \sum_{(\fraka,\frakb)\in {\cal C}_i} \sum_{(\frakc,\frakd)\in {\cal C}_j} f(\fraka,\frakb,\frakc,\frakd) \nu^y(\fraka,\frakb)\nu^z(\frakc,\frakd)\\
\overset{\eqref{eq:ffpi}}{\simeq}& \sum_{(\fraka,\frakb)\in \Delta_i} \sum_{(\frakc,\frakd)\in \Delta_j} f(\pi(\fraka,\frakb,\frakc,\frakd)) \nu^y(\fraka,\frakb)\nu^z(\frakc,\frakd)\\
\overset{*}{=}\;\;\;& \sum_{(\overline{\fraka},\overline{\frakb})\in \overline{\Gamma}_i} \sum_{(\overline{\frakc},\overline{\frakd})\in \overline{\Gamma}_j} f(\overline{\fraka},\overline{\frakb},\overline{\frakc},\overline{\frakd}) \overline{\nu}^y(\overline{\fraka},\overline{\frakb})\overline{\nu}^z(\overline{\frakc},\overline{\frakd}). 
\end{split}
\end{equation*}
where the last line ($*$) follows from the definition of $\overline{\nu}^\bullet$ above the statement of Proposition \ref{prop:Lem11}. This finishes the proof.
\end{proof}

Finally, we briefly discuss the proof of Lemma \ref{lem:FFocompare}.
\begin{proof}[Sketch of proof for Lemma \ref{lem:FFocompare}]
Recall the comment on the difference between $\cal F$ and ${\cal F}^o$ below \eqref{eq:Fodef}. It is obvious that ${\cal F}\subset{\cal F}^o$ and the analysis of $\mP({\cal F}^o\setminus {\cal F} \mid {\cal G}_{\fraka,\frakb})$ is quite similar to that of \eqref{eq:Hprime2Hprime} in the proof of Proposition \ref{prop:Lem9}, and again an application of Proposition 6.1 of \cite{S}, incurring an upper bound of $O(2^{-k^4+k^3}) \mP( {\cal F} \mid {\cal G}_{\fraka,\frakb})$.
\end{proof}

\section{$L^{2}$-approximation of the occupation measure $\mu_{n} $}\label{sec:5}

In this section, we are going finalize our preparatory works for the proof of Theorem \ref{1st} by giving in Sections  \ref{sec:geometric} and  \ref{sec:L2approx} respectively a geometric property of LERW and the final form of the key $L^2$-estimate.

We now briefly review notation and the setup. Recall that we write $\gamma_{n} $ for the LERW on $2^{-n}\mathbb{Z}^3$ in $\mathbb{D}$ from the origin. We define the renormalized occupation measure $\mu_{n}$ by
\begin{equation}\label{measure}
\mu_{n} := f_n^{-1} \sum_{x \in \gamma_{n} \cap 2^{-n} \mathbb{Z}^{3} } \delta_{x},
\end{equation} 
where $\delta_{x}$ stands for the Dirac measure at $x$ and $\beta$ is the growth exponent defined in \eqref{growth}. With this in mind, take a cube $B \subset \mathbb{D}$ as introduced in Definition \ref{def:Cubes} satisfying \eqref{eq:Breq} and  divide it into smaller cubes $B_{1}, B_{2}, \cdots, B_{M_k}$ as in Definition \ref{def:Cubes}. Recall that $X_{i}$ stands for the number of lattice points in $B_{i}$ hit by the LERW, $Y_{i}$ stands for the indicator function of the event that the LERW hits $B_i$ and $X,\ Y$ stand for their respective sums as in \eqref{X} and \eqref{Y}.

Note that by definition $\mu_n(B)$ is very close to $f_n^{-1} X$, and the error only comes from the possible overlaps or cracks between neighboring $B_i$'s in counting lattice points. In Section \ref{sec:4}, we have already proved that $X$ can be $L^2$-approximated by $\alpha_{0} Y$ with an appropriate choice of $\alpha_{0}=\alpha_{0}(n,k)$ as in \eqref{Beta}. What is left here is to show that asymptotically $\alpha_{0}/f_n$ does not really depend on $n$. To this end, in Section \ref{sec:geometric} we prove a geometric property for LERW which roughly states that when LERW hits some cube, with high probability it hits the interior of that cube. This allows us to replace the indicator function that some cube is hit by $\gamma_n$ by that by ${\cal K}$ with little extra probability cost and to conclude the section with Proposition \ref{prop3}, the final form of the key $L^2$-estimate.

\subsection{A geometric property of LERW}
\label{sec:geometric}

We begin by setting some notation. 
\begin{dfn}\label{0710-1}
Let $x \in \mathbb{D}$ and $r > 0$. Set $B = B_\infty(x,r)$. 
Suppose that 
\begin{equation}\label{0711}
\text{$r \in (0, 10^{-2})$,\ \ $ 0 < \delta < e^{- \frac{1}{r}}$, \ \ $B \subset \mathbb{D}$ \ \  and \ \ $b:= {\rm dist} (B, \{ 0 \} \cup \partial \mathbb{D} ) \in (10 r,  \frac{1}{2}) $.}
\end{equation}
 We set $\widehat{B}:=B_\infty(x,r-\delta)$ and $\check{B}:=B_\infty(x,r+\delta)$ respectively. 
\end{dfn}

\begin{prop}\label{PROP2}
Let $B$,  $\widehat{B}$, $\check{B}$, $b$ and $\delta$ be as defined in Definition \ref{0710-1}. Then there exist  universal constants $c, \, C \in (0, \infty)$  such that 
\begin{equation}\label{3eq3}
P \Big( \gamma_{n} \cap B \neq \emptyset \text{ and } \gamma_{n} \cap \widehat{B} = \emptyset \Big) \le C 
 \frac{ \log b}{\log \delta},
\end{equation}
for all $n \ge n_{\delta}$ where the constant $n_{\delta} \ge 1$ depends only on $\delta$.
Similarly, we have
\begin{equation}\label{large}
P \Big( \gamma_{n} \cap \check{B} \neq \emptyset \text{ and } \gamma_{n} \cap B = \emptyset \Big) \le C  \frac{ \log b}{\log \delta},
\end{equation}
for all $n \ge n_{\delta}$.
\end{prop} 

\begin{proof}
As two claims are very similar, we will only prove \eqref{3eq3}. 
We define $n_{\delta} \ge 1$ as the smallest integer $n$ satisfying $2^{n} \delta^{100} > 2^{100}$. 
From now on, we assume $n \ge n_{\delta}$.

Let $\tau$ be the first time that $\gamma= \gamma_{n}$ hits $B$. We will first prove that $\gamma ( \tau ) $ is not too close to an edge of $B$. If we write $x = (x_{1}, x_{2}, x_{3})$ and set $a_{i} = x_{i} -r$, the cube $B$ is described in the form of $B = \prod_{i=1}^{3} [a_{i}, a_{i} + 2r ]$. We set $e_{1}, e_{2}, \cdots , e_{12}$ for the twelve edges of $B$ and write $\widetilde{e} = \bigcup_{i=1}^{12} e_{i}$ for its union. Let 
\begin{equation*}
F = \big\{ y  \ \big| \ {\rm dist} (y , \widetilde{e} ) \le 2 \delta \big\}
\end{equation*}
be the set of points within a distance $2 \delta$ from $\widetilde{e}$. Suppose that $\gamma \cap F \neq \emptyset$. Then $S[0, T]$ must hit $F$. Now consider a line $\ell$ which satisfies that  $b \le \text{dist} (0, \ell )  \le 2b$ and that $\ell$ is parallel to either $x$-axis, $y$-axis or $z$-axis. Recall that $r$,  $\delta$ and $b$ satisfy \eqref{0711}. It then follows from \cite[Proposition 6.4.1]{LawLim} that the probability that $S [0, T]$ contains a point within a distance $2 \delta$ of the line $\ell$ is bounded above by $C  \frac{ \log b}{\log \delta}$. Using this, we have 
\begin{equation}\label{0708-1}
P \Big( \gamma \cap F \neq \emptyset \Big) \le C_{0}  \frac{ \log b}{\log \delta},
\end{equation}
for some universal constant $C_{0} < \infty$.
Therefore, if we write $F_{1}, F_{2}, \cdots , F_{6}$ for the faces of $B$, then we have
\begin{align*}
&P \Big( \gamma \cap B \neq \emptyset \text{ and } \gamma \cap \widehat{B} = \emptyset \Big) \\
&\le C_{0}  \frac{ \log b}{\log \delta}  + \sum_{i=1}^{6} P \Big( \gamma \cap B \neq \emptyset, \  \gamma \cap \widehat{B} = \emptyset, \ \gamma (\tau ) \in F_{i}, \ \gamma \cap F = \emptyset \Big).
\end{align*} 
So it suffices to show that 
\begin{equation}\label{EQEQ}
P \Big( \gamma \cap B \neq \emptyset, \  \gamma \cap \widehat{B} = \emptyset, \ \gamma (\tau ) \in F_{i}, \ \gamma \cap F = \emptyset\Big) \le C  \frac{ \log b}{\log \delta}
\end{equation}
for each $i = 1, 2, \cdots , 6$. We will only prove the inequality above for $i=1$. Without loss of generality, we can assume $F_{1} = [a_{1}, a_{1} + 2r] \times [a_{2}, a_{2} + 2r] \times \{ a_{3} \}$.

Let $\widetilde{B} = (0,0,\delta) +B$. 
We also define $\widetilde{\gamma} (t) := \gamma (t) + (0, 0, \delta)$ as the simple path obtained by translating $\gamma$  with distance $\delta$ in the direction of $z$-axis. Let $\widetilde{\tau}$ be the first time that $\widetilde{\gamma}$ hits $B$.

Suppose that $ \gamma \cap B \neq \emptyset $, $\gamma (\tau ) \in F_{1} $ and $\gamma \cap F = \emptyset $. Then clearly $\widetilde{\gamma}$ lies in $\widehat{B}$ at time $\tau$.
We write $\widetilde{t}$ for the first time that $\widetilde{\gamma}$ hits $\widetilde{B}$. By definition, we have $\widetilde{t} = \tau$. We claim that $\widetilde{\gamma} (\widetilde{\tau})$ the location of $\widetilde{\gamma}$ at its first hitting of $B$ is lying on $F_{1}$. To see it, we first note that 
\begin{equation*}
{\rm dist} \big( \widetilde{\gamma} (\widetilde{\tau}), F_{1} \big) \le \delta
\end{equation*}
since otherwise $\widetilde{\gamma}$ would have another first hitting point of $\widetilde{B}$ which is different from $\widetilde{\gamma} ( \widetilde{t} )$. However, if $\widetilde{\gamma} (\widetilde{\tau}) \notin F_{1}$, then $\gamma \cap F \neq \emptyset$ which leads a contradiction, see Figure \ref{0709-f-1} for this. Therefore, $\widetilde{\gamma} (\widetilde{\tau}) \in F_{1}$. Consequently, we have 
\begin{equation*}
\gamma \cap B \neq \emptyset ,  \  \gamma (\tau ) \in F_{1}  \text{ and }  \gamma \cap F = \emptyset \Longrightarrow  \widetilde{\gamma}  \cap \widehat{B} \neq \emptyset \text{ and } \widetilde{\gamma} (\widetilde{\tau}) \in F_{1}.
\end{equation*}
Therefore, we have
\begin{equation}\label{3EQ4}
P \Big(  \widetilde{\gamma}  \cap \widehat{B} \neq \emptyset, \  \widetilde{\gamma} (\widetilde{\tau}) \in F_{1} \Big) \ge P \Big( \gamma \cap B \neq \emptyset ,  \  \gamma (\tau ) \in F_{1} \Big) - C_{0}  \frac{ \log b}{\log \delta},
\end{equation}
where we used \eqref{0708-1} to bound on the probability that $\gamma$ hits $F$.

\begin{figure}[h]\label{supp-fig-1}
\begin{center}
\includegraphics[scale=0.55]{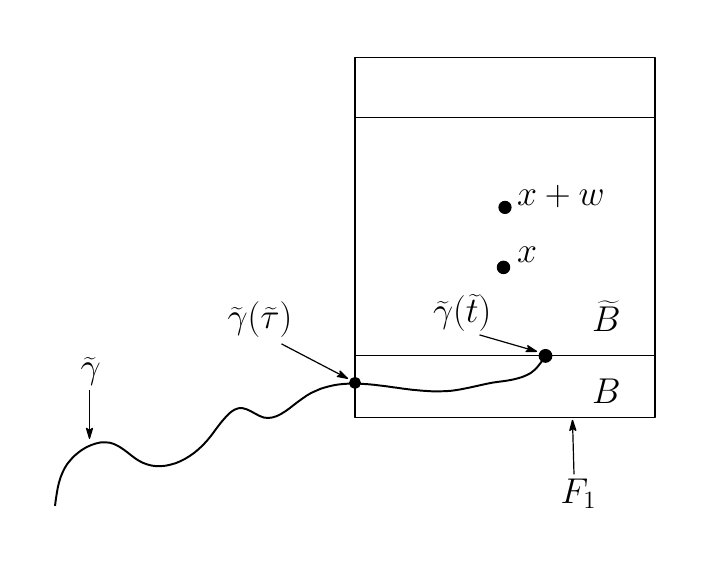}
\caption{Illustration for $B$, $\widetilde{B}$ and $\widetilde{\gamma}$. The point $x$ stands for the center of the cube $B$ while the point $x +w = x+ (0,0, \delta)$ represents the center of $\widetilde{B}$. If $\widetilde{\gamma} (\widetilde{\tau}) \notin F_{1}$ as drawn in the picture, $\widetilde{\gamma}$ must hit $F$ a neighborhood of the set of twelve edges of $B$, which is unlikely to happen.}\label{0709-f-1}
\end{center}
\end{figure}

We will show that the difference between 
\begin{equation*}
P \Big(  \widetilde{\gamma}  \cap \widehat{B} \neq \emptyset, \  \widetilde{\gamma} (\widetilde{\tau}) \in F_{1} \Big)\quad \mbox{ and } \quad
P \Big( \gamma \cap \widehat{B} \neq \emptyset ,  \  \gamma (\tau ) \in F_{1} \Big)
\end{equation*}
is accurately small by using a simple coupling argument. To do it, let $w = (0,0, \delta)$. We write $\widetilde{S}$ for the simple random walk on $2^{-n} \mathbb{Z}^{3}$ started at $w$ and set $\widetilde{T}$ for the first time that $\widetilde{S}$ exits from $\mathbb{D} + w := \{ y + w \ | \ y \in \mathbb{D} \}$. Then, by the translation invariance, $\widetilde{\gamma}$ has the same distribution as that of $\gamma^{\ast} :=\text{LE} ( \widetilde{S} [0, \widetilde{T} ] )$. We write $T$ for the first time that $\widetilde{S}$ exits from $\mathbb{D}$ as well. Let $\gamma' := \text{LE} ( \widetilde{S} [0, T] )$. We first compare 
\begin{equation*}
P \Big(  \widetilde{\gamma}  \cap \widehat{B} \neq \emptyset, \  \widetilde{\gamma} (\widetilde{\tau}) \in F_{1} \Big) = P \Big(  \gamma^{\ast} \cap \widehat{B} \neq \emptyset, \  \gamma^{\ast} (\tau^{\ast}) \in F_{1} \Big)
\end{equation*}
with
\begin{equation*}
P \Big(  \gamma'  \cap \widehat{B} \neq \emptyset, \  \gamma' (\tau') \in F_{1} \Big),
\end{equation*}
where $\tau^{\ast}$ (resp. $\tau'$) stands for the first time that $\gamma^{\ast}$ (resp. $\gamma'$) hits $B$, see Table \ref{symbols-1} for a list of symbols. To do this, let 
$$
G_{1} = \Big\{   \gamma^{\ast} \cap \widehat{B} \neq \emptyset, \  \gamma^{\ast} (\tau^{\ast}) \in F_{1}  \Big\}\quad\mbox{ and }\quad
G_{2} = \Big\{ \gamma'  \cap \widehat{B} \neq \emptyset, \  \gamma' (\tau') \in F_{1}  \Big\}.
$$
It follows from the gambler's ruin estimate (see \cite[Proposition 5.1.6]{LawLim} for this) that 
\begin{equation*}
P \Big( \text{diam} \big( \widetilde{S} [ \widetilde{T} \wedge T, \widetilde{T} \vee T ] \big) \ge \delta^{\frac{1}{2}} \Big) \le C  \delta^{\frac{1}{2}}.
\end{equation*}
Suppose that the event  $G_{3} : = \Big\{ \text{diam} \big( \widetilde{S} [ \widetilde{T} \wedge T, \widetilde{T} \vee T ] \big) \le \delta^{\frac{1}{2}} \Big\}$ occurs. Write $u$ (resp. $u'$) for the first time that $\gamma^{\ast}$ (resp. $\gamma'$) exits from $\big\{ y \in \mathbb{R}^{3} \ \big| \ |y| < 1 - 2 \sqrt{\delta} \big\}$. Then we have 
\begin{equation}\label{0709-b-1}
\gamma^{\ast} [0, u] = \gamma'  [0, u']. 
\end{equation}
The reason for this is as follows. First, suppose that $T < \widetilde{T}$. Since the diameter of $\widetilde{S} [ T, \widetilde{T} ]$ is smaller than $\delta^{\frac{1}{2}}$ by $G_{3}$, this means that $\gamma'  [0, u']$ is not destroyed by $\widetilde{S} [ T, \widetilde{T} ]$ when we make  $\gamma^{\ast}$. Therefore, it follows that $\gamma^{\ast} [0, u] = \gamma'  [0, u']$ in this case. For the case that $T > \widetilde{T}$, we show \eqref{0709-b-1} similarly.

\begin{table}[hbtp]
  \centering
  \begin{tabular}{|c|c||c|c|}
     \hline 
    Symbol  & Meaning  & Symbol  & Meaning \\
    \hline \hline 
    $S$  & SRW on $2^{-n} \mathbb{Z}^{3}$ stated at the origin  & $\gamma' $ & $ \text{LE} \big( \widetilde{S} [0, T] \big)$  \\ \hline
    $\widetilde{S}$ & SRW on $2^{-n} \mathbb{Z}^{3}$ stated at $w := (0,0, \delta )$ & $\tau$ & First time that $\gamma$ hits $B$ \\  \hline
    $T$ & First time that $S$ exits from $\mathbb{D}$ & $\widetilde{\tau}$ & First time that $\widetilde{\gamma}$ hits $B$  \\ \hline  
    $\widetilde{T}$ & First time that $\widetilde{S}$ exits from $\mathbb{D} + w$ & $\widetilde{t}$ & First time that $\widetilde{\gamma}$ hits $\widetilde{B}$ \\ \hline 
    $\widetilde{\gamma}$ & $\gamma + w$ & $\tau^{\ast}$  & First time that $\gamma^{\ast}$ hits $B$ \\ \hline 
   $\gamma^{\ast} $ & $ \text{LE} \big( \widetilde{S} [0, \widetilde{T} ] \big)$  & $\tau'$ & First time that $\gamma'$ hits $B$ \\ 
    \hline 
  \end{tabular}
  \caption{List of symbols used in the proof of Proposition \ref{PROP2}}
  \label{symbols-1}
\end{table}

Thus, if $G_{1} \cap G_{2}^{c} \cap G_{3}$ or $G_{1}^{c} \cap G_{2} \cap G_{3}$ occur, then $\widetilde{S} [v, \widetilde{T} \vee T ] \cap B \neq \emptyset$ where $v$ stands for the first time that $\widetilde{S}$ exits from $\big\{ y \in \mathbb{R}^{3} \ \big| \ |y| < 1 - 2 \sqrt{\delta} \big\}$. However, since $b:= {\rm dist} \big( B, \{0 \} \cup \partial \mathbb{D} \big)$, by our choice of $\delta$ and $b$ as in \eqref{0711}, we have $b > 2 \sqrt{\delta}$. Thus, it follows from the gambler's ruin estimate (see \cite[Proposition 5.1.6]{LawLim} for this) again that 
\begin{equation*}
P \Big( \widetilde{S} [v, \widetilde{T} \vee T ] \cap B \neq \emptyset \Big) \le  \frac{ C \, \delta^{\frac{1}{2}}}{b}.
\end{equation*}
Consequently, we see that 
\begin{equation}\label{3EQ1}
\Big| P \Big(  \widetilde{\gamma}  \cap \widehat{B} \neq \emptyset, \  \widetilde{\gamma} (\widetilde{\tau}) \in F_{1} \Big) - P \Big(  \gamma'  \cap \widehat{B} \neq \emptyset, \  \gamma' (\tau') \in F_{1} \Big) \Big| \le  \frac{ C \, \delta^{\frac{1}{2}}}{b}.
\end{equation}

We will next compare
\begin{equation*}
P \Big(  \gamma'  \cap \widehat{B} \neq \emptyset, \  \gamma' (\tau') \in F_{1} \Big)
\quad \mbox{
with
 }\quad
P \Big( \gamma \cap \widehat{B} \neq \emptyset ,  \  \gamma (\tau ) \in F_{1} \Big).
\end{equation*}
Recall that $\gamma$ is the loop-erasure of $S[0, T]$ where $S$ starts from the origin while $\gamma'$ is the loop-erasure of $\widetilde{S} [0,T]$ where $\widetilde{S}$ starts from $w = (0,0, \delta)$ (we write $T$ for the first exit time from $\mathbb{D}$ for both $S$ and $\widetilde{S}$). We also recall that $\tau$ (resp. $\tau'$) stands for the first time that $\gamma$ (resp. $\gamma'$) hits $B$. The starting point of $S$ is different from that of $\widetilde{S}$. However, considering  the path of $S$ until  the first time that it hits a plane $\Pi = \{ (p, q, \frac{\delta}{2}) \ | \ (p, q) \in \mathbb{R}^{2} \}$ and then taking a reflection of the path with respect to the plane $\Pi$, we can define $S$ and $\widetilde{S}$ on the same probability space such that 
\begin{equation}\label{0709-1}
P (H) \ge 1 - C \delta^{\frac{1}{2}},
\end{equation}
where 
\begin{align}\label{0709-2}
&\text{$t_{1}$ (resp. $t_{2}$) stands for the first time that $S$ (resp. $\widetilde{S}$) exits from $\delta^{\frac{1}{2}} \mathbb{D} $, and} \notag \\
&H = \Big\{ S (k + t_{1} ) = \widetilde{S} (k + t_{2} ) \text{ for all } k \ge 0 \Big\}.
\end{align}
Here we used the gambler's ruin estimate (see \cite[Proposition 5.1.6]{LawLim} for this) again to show that the probability that $S [0, t_{1} ]$ does not hit the plane $\Pi$ is bounded above by $C \, \delta^{\frac{1}{2}}$.

From now on, we assume the coupling of $S$ and $\widetilde{S}$ so that the inequality \eqref{0709-1} holds. Then they have a cut time with high probability in the following sense. Let $t_{1}'$ (resp. $t_{2}'$) stand for the first time that $S$ (resp. $\widetilde{S}$) exits from $\delta^{\frac{1}{3}} \mathbb{D} $.  We call $k$ a cut time for $S$ (resp. $\widetilde{S}$) if the following three conditions hold (see Figure \ref{0709-f-2} for the setup):
\begin{itemize}
\item $t_{1} \le k \le t_{1}'$ (resp. $t_{2} \le k \le t_{2}'$),

\item $S[k +1, T] \cap \big( S[t_{1}, k] \cup \delta^{\frac{1}{2}} \mathbb{D} \big) = \emptyset$ (resp. $\widetilde{S} [k +1, T] \cap \big( \widetilde{S} [t_{2}, k] \cup \delta^{\frac{1}{2}} \mathbb{D} \big) = \emptyset$),

\item $S(k) \in \delta^{\frac{1}{3}} \mathbb{D} \setminus \delta^{\frac{1}{2}} \mathbb{D}$ (resp. $\widetilde{S}(k) \in \delta^{\frac{1}{3}} \mathbb{D} \setminus \delta^{\frac{1}{2}} \mathbb{D}$).

\end{itemize}
Here we recall that $t_{1}$ and $t_{2}$ are as defined in \eqref{0709-2} and $T$ stands for the first exit time from $\mathbb{D}$ for both $S$ and $\widetilde{S}$.
Let $H'$ be the event that $S$ has such a cut time. It is then shown in the inequality  \cite[(16)]{Lawcut} that 
\begin{equation}\label{0709-3}
P (H') \ge 1 - C \delta^{c},
\end{equation}
for some universal constants $c, C \in (0, \infty)$.

Now suppose that $H \cap H'$ occurs. We write $k + t_{1}$ (with $k \ge 0$) for a cut time of $S$. Then the event $H$ ensures that $k + t_{2}$ is a cut time for $\widetilde{S}$. Furthermore, by definition of the loop-erasing procedure, one can see that 
\begin{align*}
&\gamma = \text{LE} (S[0,T]) = \text{LE} \big( S[0, k+ t_{1} ] \big) \oplus \text{LE} \big( S[k+ t_{1} + 1, T ] \big), \\
&\gamma' = \text{LE} (\widetilde{S} [0, T] ) = \text{LE} \big(\widetilde{S} [0, k+ t_{2} ] \big) \oplus \text{LE} \big( \widetilde{S} [k+ t_{2} + 1, T ] \big). 
\end{align*}
Note that the event $H$ guarantees that $\text{LE} \big( S[k+ t_{1} + 1, T ] \big) = \text{LE} \big( \widetilde{S} [k+ t_{2} + 1, T ] \big)$. Also we mention that $\text{LE} \big( S[0, k+ t_{1} ] \big), \text{LE} \big(\widetilde{S} [0, k+ t_{2} ] \big) \subset \delta^{\frac{1}{3}} \mathbb{D}$ by our definition of the cut time as above. In particular, on the event $H \cap H'$, it follows that the event
\begin{equation*}
H'' := \Big\{ \gamma (k + v_{1} ) = \gamma' (k + v_{2} ) \text{ for all } k \ge 0 \Big\}
\end{equation*}
occurs. Here $v_{1}$ (resp. $v_{2}$) stands for the first time that $\gamma$ (resp. $\gamma'$) exits from $\delta^{\frac{1}{3}} \mathbb{D}$.

\begin{figure}[h]\label{supp-fig-2}
\begin{center}
\includegraphics[scale=0.55]{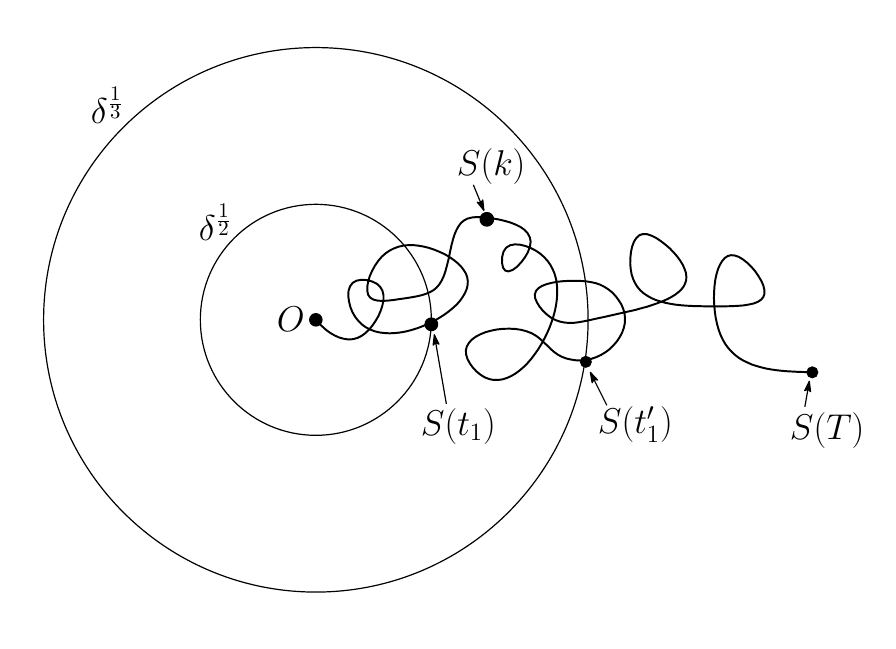}
\caption{Cut time for $S$.}\label{0709-f-2}
\end{center}
\end{figure}

Since $b:= {\rm dist} \big( B, \{0 \} \cup \partial \mathbb{D} \big)$,  our choice of $\delta$ and $b$ as in \eqref{0711} ensures that    $B \cap \delta^{\frac{1}{3}} \mathbb{D} = \emptyset$ and       
\begin{equation*}
H'' \cap \Big\{ \gamma \cap \widehat{B} \neq \emptyset ,  \  \gamma (\tau ) \in F_{1} \Big\} \Leftrightarrow H'' \cap \Big\{ \gamma' \cap \widehat{B} \neq \emptyset ,  \  \gamma' (\tau' ) \in F_{1} \Big\}. 
\end{equation*}
Thus, using \eqref{0709-1}, \eqref{0709-3} and the fact that $H \cap H' \subset H''$, we have 
\begin{align}
&P \Big( \gamma \cap \widehat{B} \neq \emptyset ,  \  \gamma (\tau ) \in F_{1} \Big) \notag \\
&= P \Big( \gamma \cap \widehat{B} \neq \emptyset ,  \  \gamma (\tau ) \in F_{1}, \ H'' \Big) + P \Big( \gamma \cap \widehat{B} \neq \emptyset ,  \  \gamma (\tau ) \in F_{1}, \ (H'')^{c} \Big) \notag \\
&=P \Big( \gamma' \cap \widehat{B} \neq \emptyset ,  \  \gamma' (\tau' ) \in F_{1}, \ H'' \Big) + O \big( \delta^{c} \big) \notag \\
&= P \Big( \gamma' \cap \widehat{B} \neq \emptyset ,  \  \gamma' (\tau' ) \in F_{1} \Big) +O \big( \delta^{c} \big), \label{3EQ2}
\end{align}
for some universal constant $c > 0$. 

Combining \eqref{3EQ1} with \eqref{3EQ2}, we see that 
\begin{equation}\label{3EQ3}
\Big| P \Big(  \widetilde{\gamma}  \cap \widehat{B} \neq \emptyset, \  \widetilde{\gamma} (\widetilde{\tau}) \in F_{1} \Big) - P \Big( \gamma \cap \widehat{B} \neq \emptyset ,  \  \gamma (\tau ) \in F_{1} \Big) \Big| \le \frac{C \, \delta^{c}}{b},
\end{equation}
for some universal constants $c, C > 0$. On the other hand, by \eqref{3EQ4}, we have 
\begin{align*}
&P \Big( \gamma \cap \widehat{B} \neq \emptyset ,  \  \gamma (\tau ) \in F_{1} \Big) \\
&= P \Big( \gamma \cap B \neq \emptyset ,  \  \gamma (\tau ) \in F_{1} \Big) - P \Big( \gamma \cap B \neq \emptyset, \  \gamma \cap \widehat{B} = \emptyset, \ \gamma (\tau ) \in F_{1} \Big) \\
&\le P \Big(  \widetilde{\gamma}  \cap \widehat{B} \neq \emptyset, \  \widetilde{\gamma} (\widetilde{\tau}) \in F_{1} \Big) + C_{0} \frac{ \log b}{\log \delta}  - P \Big( \gamma \cap B \neq \emptyset, \  \gamma \cap \widehat{B} = \emptyset, \ \gamma (\tau ) \in F_{1} \Big).
\end{align*}
From this and \eqref{3EQ3}, we finally conclude that 
\begin{align*}
&P \Big( \gamma \cap B \neq \emptyset, \  \gamma \cap \widehat{B} = \emptyset, \ \gamma (\tau ) \in F_{1} \Big) \\
& \le P \Big(  \widetilde{\gamma}  \cap \widehat{B} \neq \emptyset, \  \widetilde{\gamma} (\widetilde{\tau}) \in F_{1} \Big) - P \Big( \gamma \cap \widehat{B} \neq \emptyset ,  \  \gamma (\tau ) \in F_{1} \Big) + C_{0} \frac{ \log b}{\log \delta}  \\
&\le \frac{C \, \delta^{c}}{b} + C_{0} \frac{ \log b}{\log \delta}.
\end{align*}
However, our choice of $\delta$ and $b$ as in \eqref{0711} ensures that $\frac{ \delta^{c}}{b}$ is bounded above by $C \frac{ \log b}{\log \delta}$ for some universal constant $C < \infty$. Therefore, we get \eqref{EQEQ} and finish the proof.
\end{proof}
\begin{rem} A similar coupling of two LERW's with different initial configurations can also be found in  \cite[Section 2.6]{Lawrecent}. 
\end{rem}

We define 
$\{ \gamma_{n} \}_{n \ge 1}$ and ${\cal K}$ on the same probability space so that 
\begin{equation}\label{3eq1}
\lim_{n \to \infty} d_{\text{Haus}} \big( \gamma_{n}, {\cal K} \big) = 0
\end{equation}
almost surely. 
%
%
%
%
%
%
Let $B_n=B\cap 2^{-n}\mathbb{Z}^3$ be the discretization of $B$. 
Let $\widetilde{Y}_{n}$,$\widetilde{Y}_{n}^{-}$, $\widetilde{Y}_{n}^{+}$ and $\widetilde{Z}$ be the indicator function of the event that $\gamma_{n} \cap B_n \neq \emptyset$, $\gamma_{n} \cap B_n^\circ \neq \emptyset$, $\gamma_{n} \cap \overline{B_n} \neq \emptyset$ and  ${\cal K} \cap B \neq \emptyset$ respectively. Then we have the following corollary.

\begin{cor}\label{COR}
With the notation as defined above, it follows that for all $\delta \in (0, e^{- \frac{1}{r}} )$ there exists $N= N_{\delta}$ depending on $\delta$ such that for all $n \ge N_{\delta}$
\begin{equation}
P \Big( \widetilde{Y}_{n}=\widetilde{Y}_{n}^{-}= \widetilde{Y}_{n}^{+} = \widetilde{Z} \Big) \ge 1 - C \, \frac{\log b}{\log \delta}
\end{equation}
for some universal constant $C < \infty$. 
\end{cor}

\begin{proof}
Without loss of generality we only prove the claim for $\widetilde{Y}_{n}$.
We recall that $\{ \gamma_{n} \}_{n \ge 1}$ and ${\cal K}$ are coupled such that \eqref{3eq1} holds. Therefore, it follows that there exists $N= N_{\delta}$ depending on $\delta$ such that for all $n \ge N_{\delta}$,
\begin{equation}\label{0711-1}
P \Big( d_{\text{Haus}} \big( \gamma_{n}, {\cal K} \big) \ge \frac{\delta}{2}  \Big) \le \exp \big\{ - \delta^{-2} \big\}.
\end{equation}
We may assume that $N_{\delta} \ge n_{\delta}$ where $n_{\delta}$ is the constant as in Proposition \ref{PROP2}.

Now suppose that $\widetilde{Y}_{n} \neq \widetilde{Z}$. Then either $A_{1} := \{ \widetilde{Y}_{n} = 1, \widetilde{Z}= 0 \}$ or $A_{2} := \{ \widetilde{Y}_{n}= 0, \widetilde{Z}=1 \}$ occur. Note that by \eqref{0711-1}
\begin{equation*}
P (A_{1}) \le \exp \big\{ - \delta^{-2} \big\} + P \Big( d_{\text{Haus}} \big( \gamma_{n}, {\cal K} \big) < \frac{\delta}{2}, \  \ A_{1} \Big).
\end{equation*}
So we further assume that $d_{\text{Haus}} \big( \gamma_{n}, {\cal K} \big) < \frac{\delta}{2}$ and $A_{1}$ hold. Recall that 
$\widehat{B}$ and $\check{B}$ are as defined in Definition \ref{0710-1}. If $\gamma_{n}$ hits $\widehat{B}$, then the condition that $d_{\text{Haus}} \big( \gamma_{n}, {\cal K} \big) < \frac{\delta}{2} $ implies that ${\cal K} \cap B \neq \emptyset$. Thus, the event $A_{1}$ ensures that $\gamma_{n} \cap \widehat{B} = \emptyset$ and $\gamma_{n} \cap B \neq \emptyset$. Therefore, by \eqref{3eq3}, we have
\begin{equation*}
P (A_{1}) \le C \, \frac{\log b}{\log \delta},
\end{equation*}
where we used $ \exp \big\{ - \delta^{-2} \big\} \le C \frac{\log b}{\log \delta}$ because of our choice of $\delta$ and $b$ as in \eqref{0711}.
Similarly,
\begin{equation*}
P (A_{2}) \le C \, \frac{\log b}{\log \delta},
\end{equation*}
which completes the proof.
\end{proof}

\subsection{$L^2$-approximation}
\label{sec:L2approx}
We continue our work from Proposition \ref{prop4.6}.
In this subsection, we first obtain an accurate asymptotics of $\alpha_0$. Then, we will rewrite Proposition \ref{prop4.6} in the form of Proposition \ref{prop3}, which allows us to show in the next section that any possible scaling limit of the occupation measure must be a measurable function of the scaling limit ${\cal K}$.

We start with the asymptotics of $\alpha_0$. Recall the convention regarding $\simeq$ in \eqref{eq:simeqconv}.
\begin{prop}\label{prop:alpha0}
There exists universal constant $c > 0$, and constants $\zeta(k)$, $N_{k}$ depending only on $k \ge 1$ such that for all $k \ge 1 $ and $n \ge N_{k}$, 
\begin{equation}\label{eq:alphazeta}
\alpha_{0}(n,k) \simeq \zeta(k){f_n}.
\end{equation}
\end{prop}

\begin{proof}
By definition of $\alpha_{0}$ in \eqref{Beta}, we see that
\begin{equation}\label{atuiyo}
 \frac{\alpha_{0}}{f_n} = \frac{ \sum_{x \in B_{0} \cap 2^{-n} \mathbb{Z}^{3}} P \big( x \in \gamma_{n} \big) }{f_n P \big( Y_{0} = 1\big) },
\end{equation}
where $Y_{0}$ denotes the indicator function of the event that $\gamma_{n}$ hits $B_{0}$, as defined above \eqref{Beta}.

By \eqref{TRANSINVA} and \eqref{eq:fndef}, we have 
 $$
 \frac{ \sum_{x \in B_{0} \cap 2^{-n} \mathbb{Z}^{3}} P \big( x \in \gamma_{n} \big) }{f_n P \big( Y_{0} = 1\big) } \simeq \frac{|B_0|P(x_0\in \gamma_n)}{f_n P \big( Y_{0} = 1\big) }=  \frac{c_{1} 2^{-3k^{4}}}{ P \big( Y_{0} = 1 \big)}
 $$
 for some $c_1>0$.
On the other hand, applying Corollary \ref{COR} to the case that $x=x_{0}$, $r = 3 \cdot 2^{-k^{4}} $ and $\delta = 2^{- 2^{k^{10}}}$ (note that $\frac{\log b}{\log \delta} \le C 2^{-k^{10}} $ in this case), we see that 
\begin{equation}\label{0712-o}
\Big| P \big( Y_{0} = 1 \big) - P \big( W_k = 1 \big) \Big| \le c 2^{-k^{10}},
\end{equation}
for all $n \ge n_{k}'$ where the constant $n_{k}'$ depends only on $k$ and $W_{k}$ is the indicator function of the event that ${\cal K}$ hits $B_0$:
\begin{equation}\label{WK}
W_k: = 1\{{\cal K} \cap B_{0} \neq \emptyset\}.
\end{equation}
(Note that when applying Corollary \ref{COR} by picking $B$ there as $B_0$ in the definition of $Y_0$ (see above \eqref{Beta}), one can verify that $1_{\widetilde{Y}_{n}^{-}} \leq 1_{Y_0} \leq 1_ {\widetilde{Y}_{n}^{+}}$.)
It then follows from \eqref{20221118-1} that
\begin{equation}\label{1026-0912}
P \big( Y_0 = 1 \big) \asymp 2^{-(3- \beta ) k^{4}}. 
\end{equation}
Hence,
\begin{equation*}
P \big( Y_{0} = 1 \big) \simeq  P \big( W_k = 1 \big),
\end{equation*}
for all $n \ge n_{k}'$. Thus, writing $N_{k} = n_{k} \vee n_{k}'$ and letting
\begin{equation}\label{eq:zetadef}
\zeta(k):=\frac{c_1 2^{-3 k^{4}}}{P \big( W_k = 1 \big) },
\end{equation}
we conclude that the claim \eqref{eq:alphazeta} holds for all $n \ge N_{k}$. Note that by \eqref{1026-0912} we have
\begin{equation}\label{0726}
\zeta(k) \asymp 2^{- \beta k^{4}}.
\end{equation}
Consequently, we have 
\begin{equation}\label{2022-10-8}
\frac{\alpha_{0}}{f_n} \simeq \zeta (k),
\end{equation}
which finishes the proof. 
\end{proof}

We are ready to show the following proposition, which bears some similarity to   \cite[Proposition 4.11]{GPS} which also gives an $L^2$-approximation of the occupation measure  in the case of critical planar percolation.

\begin{prop}\label{prop3}
With the notation above, it follows that there exist universal constant $c, C, c^{\ast} > 0$ and a constant $N_{k}$ depending only on $k \ge 1$ such that for all $k \ge 1$ and $n \ge N_{k}$
\begin{equation}
E \Big[ \Big(  \mu_{n} (B) -  \zeta(k) Y \Big)^{2} \Big]\le C 2^{-c k^{2}},
\end{equation}
where $\zeta (k)$ is as defined in \eqref{eq:zetadef}.
\end{prop}

\begin{proof}
Note that we have
\begin{align}\label{0711-e-1}
\parallel \mu_{n} (B) -  \zeta(k) Y  \parallel_{2} \;&\le\; \parallel \mu_{n} (B) -  f_n^{-1} X\parallel_2+  f_n^{-1}\parallel  X - \alpha_{0}  Y \parallel_{2} + \parallel  f_n^{-1}\alpha_{0} Y  -  \zeta(k) Y  \parallel_{2} \notag \\
 &:=\; {\rm I}_1+{\rm I}_2+{\rm I}_3.
\end{align}
where $\| \cdot \|_{2}$ stands for the $L^{2}$-norm.

We start with ${\rm I}_1$. By definition we have
$$
|\mu_n(B)-f_n^{-1} X| \leq f_n^{-1} \sum_{i=1}^{M_k} \sum_{x\in \partial_i B_i  \cup \partial B_i}(1_{x\in \gamma_n}).
$$
It then follows from the two-point function estimates \eqref{M7prime} (and the corresponding version for $B\subset \mathbb{D}\setminus \frac{1}{3}\mathbb{D}$) that 
$$
{\rm I}_1 \leq f_n^{-1} E \left[\left(\sum_{i=1}^{M_k} \sum_{x\in \partial_i B_i  \cup \partial B_i}(1_{x\in \gamma_n}) \right)^2\right]^{1/2} \leq C  f_n^{-1}2^{k^4-k} 2^{(2n- 2k)} 2^{(\beta-3)n}  \leq C2^{k^4-3k} 2^{-n} \leq C' 2^{-k^2},
$$
when $n\geq N_k$ for some sufficiently large $k$.
Next,  by Propositions \ref{prop4.6}, \eqref{eq:fnasymp} and \eqref{improve}
$${\rm I}_2=f_n^{-1}E \Big( \big( X - \alpha_{0} Y \big)^{2} \Big)^{\frac{1}{2}}   \leq C 2^{-c k^{2}} f_n^{-1} E (X) \leq  C' 2^{-c k^{2}}.
$$
Finally,  by \eqref{eq:fnasymp} and Proposition \ref{prop:alpha0},
$${\rm I}_3=
 \Big| \alpha_{0} f_n^{-1}- \zeta(k) \Big| \cdot E (Y^{2})^{\frac{1}{2}} 
\le\;C 2^{- c k^{3}} \zeta(k) E (Y^{2})^{\frac{1}{2}}.
$$
It follows from Proposition \ref{2nd-mom-mom} that
$$E (Y^{2} ) \le C 2^{2 \beta (k^{4} - k)}.
$$
Since by \eqref{0726}, $\zeta (k) \asymp 2^{-\beta k^{4} }$, the right hand side of \eqref{0711-e-1} is thus smaller than $C 2^{-c k^{2}}$ and we finish the proof.
\end{proof}

\section{Convergence of the occupation measure}\label{sec:6}
We recall that the rescaled occupation measure $\mu_{n}$ is defined as in \eqref{measure}. The goal of this section is to prove Theorem \ref{1st}, which states that the joint law of $(\gamma_{n}, \mu_{n} )$ converges weakly to some $( {\cal K}, \mu )$ with respect to the topology of ${\cal H} (\overline{\mathbb{D} } ) \times {\cal M} (\overline{\mathbb{D}} ) $ where ${\cal K}$ stands for Kozma's scaling limit of the LERW (see Section \ref{metric} for ${\cal H} (\overline{\mathbb{D} } ) $ and ${\cal M} (\overline{\mathbb{D} } ) $).

To achieve this, 
we will establish the uniqueness of sub-sequential limits  of $(\gamma_{n}, \mu_{n} )$ in Section  \ref{sec:unique}.

\subsection{Uniqueness of the sub-sequential limit}\label{sec:unique}
In Proposition \ref{tight} we proved that $ \{ (\gamma_{n}, \mu_{n}) \}_{n \ge 1}$ is tight. We also know that $\gamma_{n}$ converges weakly to ${\cal K}$ as $n \to \infty$ with respect to the Hausdorff distance (see Section \ref{SCALING}). Therefore, we can find a sub-sequence $\{ n_{l} \}_{l \ge 1}$ and a random measure $\mu^{\ast}$ such that 
\begin{equation}\label{sub}
(\gamma_{n_{l}}, \mu_{n_{l}}) \xrightarrow{ d } ({\cal K}, \mu^{\ast} ) \ \ \text{ as } \  l \to \infty.
\end{equation}
In this subsection, we will show that $\mu^{\ast}$ does not depend on the chosen sub-sequence $\{ n_{l} \}_{l \ge 1}$, i.e., there exists a unique possible choice for the measure $\mu^{\ast}$. In order to prove this, it is enough to show that  $\mu^{\ast}$ is a measurable function of ${\cal K}$. This is summarized in the following theorem, which is the last piece to complete the proof of Theorem \ref{1st}.

\begin{thm}\label{goal}
Any sub-sequential limit $\mu^{*}$ defined from \eqref{sub}, is a measurable function of $\cal K$.
\end{thm}
Although the proof of this theorem is quite long, it is nevertheless very straightforward. To show the measurability, we consider the inner product of the measure $\mu^*$ with one of the orthogonal bases, and show that it can be $L^2$-approximated by various  intermediate quantities $J_1$ through $J_4$ (see Table \ref{symbols-3} for a recap of the definitions), where $J_1$ and $J_2$ are defined through the discrete occupation measures, $J_3$ through the trace of LERW and $J_4$ through ${\cal K}$, the scaling limit. Note that the key $L^2$ approximation Proposition \ref{prop3} is applied here to communicate between $J_2$ and $J_3$. Finally, we show that the $L^2$-convergence in the approximation scheme above can be upgraded to a.s.\ convergence and the inner product is indeed a function of $\cal K$.
\begin{proof}

By Skorokhod's representation theorem, we can define $ \{ (\gamma_{n_{l}}, \mu_{n_{l}}) \}_{l \ge 1}$ and $({\cal K}, \mu^{\ast} )$ on the same probability space such that 
\begin{equation}\label{eq:6.2I}
\mbox{$\gamma_{n_{l}} \to {\cal K}$ almost surely with respect to the topology of $({\cal H} (\overline{\mathbb{D} } ), d_{\text{Haus}} )$;}
\end{equation}
\begin{equation}\label{eq:6.2II}
\mbox{ $\mu_{n_{l}} \to \mu^{\ast}$ almost surely with respect to the topology of the weak convergence on ${\cal M}  (\overline{\mathbb{D} } )$.}
\end{equation}

Let $( \Phi (\overline{\mathbb{D} } ) , \parallel \cdot \parallel )$ be the space of real-valued continuous functions on $\overline{\mathbb{D} } $ endowed with the uniform norm $\parallel \cdot \parallel$. We take a countable basis $\{ \phi_{j} \}_{j \ge 1 }$ of the space $(\Phi (\overline{\mathbb{D} } ) , \parallel \cdot \parallel )$. We may assume that $\parallel \phi_{j} \parallel \le 1$ for all $j \ge 1$.

Fix $j \ge 1$. We will show that $\mu^{\ast} ( \phi_{j} ) := \int_{\overline{\mathbb{D} }} \phi_{j} (x) \, d \mu^{\ast} ( x)$ is a measurable function of ${\cal K}$. Recall that $\parallel \phi_{j} \parallel \le 1$. Since $\phi_{j}$ is a uniformly continuous function on $\overline{\mathbb{D} } $, for each $m \ge 1$ there exists $k_m \in \mathbb{N}$ such that 
\begin{equation}\label{0}
\max \Big\{ | \phi_{j} (x) - \phi_{j} (y) | \ \Big| \ x, y \in \overline{\mathbb{D} }, \ | x- y | \le 4 \cdot 2^{-k_m} \Big\} \le 2^{-m}.
\end{equation}
We may assume that $k_m \ge m$ and $k_m < k_{m+1}$ for all $m \ge 1$. When there is no confusion, we will write $k$ instead of $k_m$ for simplicity.  

We consider a covering of $\overline{\mathbb{D} } $ by a collection of cubes $\{ E_{q} \}_{q=1}^{L_{k}}$ in the form of 
\begin{equation}\label{0724-1}
E_{q}:= B_\infty(x_{q}, 2^{-k -1}) \ \text{ for } \ 1 \le q \le L_{k}
\end{equation}
where $\| \cdot \|_{\infty}$ stands for the $L^{\infty}$-norm. We assume that $\{ E_{q} \}_{q=1}^{L_{k}}$ satisfies 
\begin{align}
&\text{(i) \  $|x_{q} | < 2$ for all $1 \le q \le L_{k}$ and $\| x_{q} - x_{q'} \|_{\infty} \ge 2^{-k}$ for $q \neq q'$,} \notag \\
&\text{(ii) \  $E_{q} \cap E_{q'} \cap 2^{-n} \mathbb{Z}^{3} = \emptyset $ for $q \neq q'$.} \label{0724-2}
\end{align}
Note that $L_{k} \asymp 2^{3 k}$. 
 
 We set 
\begin{equation}\label{0715-a-1}
I_{k} = \Big\{ 1 \le q \le L_{k} \ \Big| \ {\rm dist} \big( E_{q}, \{ 0 \} \cup \partial \mathbb{D}  \big) > 2^{-k} \Big\}
\end{equation}
for the set of indices of ``typical" cubes. Note any typical cube $E_q$ for $q\in I_k$ satisfies the requirement \eqref{eq:Breq} for the cube $B$ considered in Definition \ref{def:Cubes}.
Let 
\begin{equation}\label{0714-d-1}
H_{k} = \bigcup_{q \in I_{k} } E_{q}
\end{equation}
be the union of such typical cubes. We also write $I_{k}^{c} = \big\{ 1 \le q \le L_{k} \ \big| \ q \notin I_{k}  \big\}$
for the complement of $I_{k}$.  

For each $q \in I_{k}$, we further decompose $E_{q}$ into smaller cubes $D^{1}_{q}, D^{2}_{q}, \cdots , D^{M_k}_{q}$ in the form of
\begin{equation}\label{0724}
D_{q}^{i} = E_{q} \cap B_\infty(x_{q}^{i},2^{-k^{4} -1}) \ \text{ for }  \  1 \le i \le M_k,
\end{equation}
where the following conditions are satisfied:
\begin{align}
&\text{(a) \  $x_{q}^{i} \in E_{q}$ for each $1 \le i \le M_k$ \ \ \ and  \ \ \ $E_{q} = \bigcup_{i=1}^{M_k} D^{i}_{q}  $,
} \notag \\
&\text{(b) \  $\| x^{i}_{q} - x^{i'}_{q} \|_{\infty} \ge  2^{- k^{4}} $ \  \ \ and  \ \ \ $D^{i}_{q} \cap D^{i'}_{q} \cap 2^{-n} \mathbb{Z}^{3} = \emptyset $  \ \ for  \ $i \neq i'$.} \label{0726-1}
\end{align}
Note that each $D^{i}_{q}$ is the intersection of $E_{q}$ and a cube of side length  $2^{- k^{4}}$ centered at $x^{i}_{q}$. Note that these cubes ``roughly'' corresponds to the discrete cubes $B_i$'s considered in Definition \ref{def:Cubes}, but the latter, being a discrete object, may leave cracks or overlappings at the boundary. We write $B_i^q$ for the discrete cube corresponding to $D^{i}_{q}$.   We also note that $M_k$, the number of smaller cubes satisfies $M_k \asymp 2^{3 (k^{4} - k )}$ for every $1 \le q \le L_{k}$, and does not depend on $q$ thanks to the translation invariance.

Recall that we assume the coupling of $ \{ (\gamma_{n_{l}}, \mu_{n_{l}}) \}_{l \ge 1}$ and $({\cal K}, \mu^{\ast} )$ such that both \eqref{eq:6.2I} and \eqref{eq:6.2II} hold.
From now on, we write $n = n_{l}$ for simplicity.
Let 
\begin{align}
&Y^{i}_{q} := {\bf 1} \Big\{ \gamma_{n} \cap B_{i}^{q} \neq \emptyset \Big\}, \label{IND1} \\
&\widetilde{W}^{i}_{q} := {\bf 1} \Big\{ {\cal K} \cap D^{i}_{q} \neq \emptyset \Big\}. \label{IND3}
\end{align}
Applying Corollary \ref{COR}  to the case that $r = 3 \cdot 2^{- k^{4}}$, $b > 2^{- k}$ and $\delta = 2^{- 2^{k^{21}}}$, it follows that there exist a universal constant $C < \infty$ and a constant $N_{k} < \infty$ depending only on $k$ such that for all $n \ge N_{k}$, $i \in I_{k}$  and $1 \le q \le M_k$,
\begin{equation}\label{0714-c-1}
P \big( Y^{i}_{q} = \widetilde{W}^{i}_{q} \big) \ge 1 - C 2^{- k^{20}}.
\end{equation}
Since the total number of smaller cubes $D^{i}_{q}$ is bounded above by $L_{k}  M_k \asymp 2^{3 k^{4}}$, by \eqref{0714-c-1}, we have that 
\begin{equation}\label{0714-c-2}
P \Big( Y^{i}_{q} = \widetilde{W}^{i}_{q} \ \text{ for all  $q \in I_{k}$  and $1 \le i \le M_k$}  \Big) \ge 1 - C 2^{- k^{19}},
\end{equation}
for all $n \ge N_{k}$. 

Recall that $H_{k}$ is the union of typical cubes as defined in \eqref{0714-d-1}. We will first compare 
\begin{equation}\label{0715}
\mu_{n} (\phi_{j} ) = \int_{\overline{\mathbb{D} } }  \phi_{j} d \mu_{n}
\quad\mbox{ with }\quad
J_{1} :=  \int_{H_{k}}  \phi_{j} d \mu_{n}.
\end{equation}
Since $\parallel \phi_{j} \parallel \le 1$ for all $j \ge 1$, it holds  that
\begin{equation*}
\Big| \mu_{n} (\phi_{j} )  - J_{1} \Big| \le \int_{\overline{\mathbb{D} } \setminus H_{k} } | \phi_{j}| d \mu_{n} \le \mu_{n} \Big( \overline{\mathbb{D} } \setminus H_{k} \Big).
\end{equation*}
We will prove that 
\begin{equation}\label{0714-d-2}
E \bigg\{  \Big[  \mu_{n} \Big( \overline{\mathbb{D} } \setminus H_{k} \Big) \Big]^{2} \bigg\} \le C 2^{-2k}
\end{equation}
for some universal constant $C < \infty$. To see this, we write 
\begin{itemize}
\item $g_{1} $  for the number of points in $A_{1} := 2^{-n} \mathbb{Z}^{3} \cap \{ x \in \mathbb{R}^{3} \ | \ |x| \le 2^{-k} \}$  hit by $\gamma_{n}$;

\item $g_{2}$ for the number of points in $A_{2} := 2^{-n} \mathbb{Z}^{3} \cap \{ x \in \mathbb{R}^{3} \ | \  1- 2^{-k} \le |x| \le 1 \}$ hit by $\gamma_{n}$;

\item $g_{3}$ for the  number of points on the boundaries of $E_q$ for each $q$ hit by $\gamma_{n}$;
\end{itemize}
By Cauchy-Schwarz inequality and \eqref{eq:fnasymp}, it follows that 
\begin{equation}\label{0714-d-3}
\Big[  \mu_{n} \Big( \overline{\mathbb{D} } \setminus H_{k} \Big) \Big]^{2} \le C \cdot 2^{- 2 \beta n} \, \big( g_{1}^{2} + g_{2}^{2}+g_{3}^{2} \big).
\end{equation}
We will first estimate $E (g_{1}^{2} ) =\sum_{x, y \in A_{1} } P ( x, y \in \gamma_{n} ) $. By
applying \eqref{M7prime} in summation, we have
\begin{equation}\label{0714-f-1-1}
E (g_{1}^{2} ) \le C (2^{-k} \cdot 2^{n} )^{2 \beta},
\end{equation}
for some universal constant $C < \infty$. 
Similar bounds hold for $g_2$ and $g_3$.
We hence obtain \eqref{0714-d-2} as desired, which gives
\begin{equation}\label{3}
\parallel \mu_{n} (\phi_{j} ) - J_{1} \parallel _{2} \le C_{1} 2^{-k}
\end{equation}
for some universal constant $C_{1} < \infty$. 

By definitions of $H_{k}$ and $J_{1}$ as in \eqref{0714-d-1} and \eqref{0715}, and dealing with the boundary we note that  
\begin{equation}\label{eq:J1J1p}
\parallel J_{1} - J'_{1} \parallel_2 \leq C2^{-k^2},
\end{equation}
when $n\geq N_k$ where
\begin{equation}\label{4}
J'_{1} = \sum_{q \in I_{k}} \int_{E_{q}} \phi_{j} d \mu_{n}.
\end{equation}

Recall that $x_{q}$ stands for the center of $E_{q}$. We next compare $J_{1}$ (in fact $J'_{1}$ with the Riemann sum
\begin{equation}\label{5}
J_{2}:= \sum_{q \in I_{k}}  \phi_{j} ( x_{q} ) \mu_{n} (E_{q}).
\end{equation}
 Since $|x- y| \le 3 \cdot 2^{-k_m}$ for all $x, y \in E_{q}$ (note that we write $k = k_{m}$ for simplicity)
and for any $i \in I_{k}$,  using \eqref{0}, we see that 
\begin{equation*}
|J'_{1} - J_{2}| \le 2^{-m} \mu_{n} (\overline{\mathbb{D} }).
\end{equation*}
However, 
as $E ( M_{n}^{2} ) \leq C \big( 2^{2 n}  \Es ( 2^{n} ) \big)^{2}$ by \cite[Theorem 8.4]{S} (or directly from asymptotics of two-point functions such as \eqref{M7prime}), 
by \eqref{eq:fnasymp} we see that \begin{equation*}
\sup_{n} E \bigg\{  \Big[ \mu_{n}  (\overline{\mathbb{D} }) \Big]^{2} \bigg\} < \infty.
\end{equation*}
This joined with \eqref{eq:J1J1p} implies 
\begin{equation}\label{6}
\parallel J'_{1} - J_{2} \parallel _{2}  \le C_{2} 2^{-m}
\end{equation}
for some universal constant $C_{2} < \infty$.

We now  replace $J_{2}$ by some macroscopic quantity $J_{3}$ defined as follows. Take $q  \in I_{k}$. 
Using the indicator function $Y^{i}_{q}$ as defined in \eqref{IND1}, we then set 
\begin{equation}\label{IND2}
U_{q} := \sum_{i=1}^{M_k} Y^{i}_{q}.
\end{equation}

We are ready to replace $\mu_{n} (E_{q} ) $ with some macroscopic quantity. Applying Proposition \ref{prop3} to the case that $B = E_{q}$ and $Y = U_{q}$, we see that for each $i \in I_{k}$
\begin{equation}\label{7}
E \bigg\{  \Big[ \mu_{n} (E_{q} )  - \zeta(k) \, U_{q} \Big]^{2} \bigg\} \le C 2^{- c k^{2}},
\end{equation}
where $c, C$ are universal constants as in Proposition \ref{prop3} and $\zeta(\cdot)$ is defined in \eqref{eq:zetadef}. 
With this in mind, we define 
\begin{equation}\label{8}
J_{3} := \sum_{q \in I_{k}}  \phi_{j} ( x_{q} )\zeta(k) U_{q}.
\end{equation}
Since $\parallel \phi_{j} \parallel \le 1$ for all $j \ge 1$ and $J_{2}$ is defined as in \eqref{5},  it follows from \eqref{7} that 
\begin{equation}\label{9}
\parallel J_{2} - J_{3} \parallel _{2} \le  \sum_{q \in I_{k}} \parallel Z_{q} \parallel_{2}  \le  C 2^{- c k^{2}} L_{k} \le C 2^{- c k^{2}} 2^{3 k} \le C_{3} 2^{- c_{4} k^{2}}
\end{equation}
where 
 we write
\begin{equation*}
Z_{q} := \Big| \mu_{n} (E_{q} )  -\zeta(k) U_{q} \Big|.
\end{equation*}

\begin{table}[hbtp]
  \centering
  \begin{tabular}{|c|c||c|c|}
    \hline \hline
    Symbol  & Meaning  & Symbol  & Meaning     \\
    \hline \hline
    $({\cal K}, \mu^{\ast} )$ & Limit of $ \{ (\gamma_{n_{l}}, \mu_{n_{l}}) \}_{l \ge 1}$  &
    $\{ \phi_{j} \}_{j \ge 1 }$ & Basis of the space $( \Phi (\overline{\mathbb{D} } ) , \parallel \cdot \parallel )$  \\ \hline
    $\{ E_{q} \}_{q=1}^{L_{k}}$ & Covering of $\overline{\mathbb{D} } $ satisfying \eqref{0724-1} and \eqref{0724-2} &
    $I_{k} $ & $ \big\{ 1 \le i \le L_{k} \ \big| \ {\rm dist} \big( E_{q}, \{ 0 \} \cup \partial \mathbb{D}  \big) > 2^{-k} \big\} $ \\ \hline
    $H_{k} $ & $ \bigcup_{q \in I_{k} } E_{q} $ &
    $\{ D^{i}_{q}  \}_{i =1}^{M_k}$ & Covering of $E_{q}$ satisfying \eqref{0724} and \eqref{0726-1} \\ \hline
    $Y^{i}_{q} $ & Indicator function  $ {\bf 1} \big\{ \gamma_{n} \cap D^{i, +}_{q} \neq \emptyset \big\} $ &
    $\widetilde{W}^{i}_{q} $ & Indicator function $ {\bf 1} \big\{ {\cal K} \cap B_{i}^{q} \neq \emptyset \big\} $ \\ \hline
    $J_{1} $ & $\int_{H_{k}}  \phi_{j} d \mu_{n}$ &
    $J_{2} $ & $\sum_{q \in I_{k}}  \phi_{j} ( x_{q} ) \mu_{n} (E_{q})$ \\ \hline
    $J_{3}$ &  $\sum_{q \in I_{k}}  \phi_{j} ( x_{q} )\zeta(k) U_{q}$ &
    $J_{4}$ & $\sum_{q \in I_{k}}  \phi_{j} ( x_{q} ) \zeta(k) \widetilde{W}_{q}$ \\
   \hline \hline
  \end{tabular}
  \caption{List of symbols used in the proof of Theorem \ref{goal}}
  \label{symbols-3}
\end{table}

Finally, we replace $J_{3}$ by $J_{4}$ which is a measurable quantity with respect to the scaling limit ${\cal K}$ as follows. 
For each $q \in I_{k}$, we let 
\begin{equation}\label{IND4}
\widetilde{W}_{q} := \sum_{i=1}^{M_k} \widetilde{W}^{i}_{q}.
\end{equation}
Now we define
\begin{equation}\label{macro}
J_{4} = J_{4}^{j, k} := \sum_{q  \in I_{k}}  \phi_{j} ( x_{q} ) \zeta(k) \widetilde{W}_{q},
\end{equation}
see \eqref{eq:zetadef} for $\zeta(\cdot)$.
Note that $J_{4}$ is a measurable function of ${\cal K}$ depending on $j$ and $k$ (we recall that $j$ is a subscript for $\phi_{j}$ the basis of the space $(\Phi (\overline{\mathbb{D} } ) , \parallel \cdot \parallel )$ and that we write $k = k_{m}$, see \eqref{0} for $k_{m}$). We already know that the inequality  \eqref{0714-c-2} holds for all $n \ge N_{k}$.   With this in mind, we define the event $F$ by 
\begin{equation}\label{11}
F := \Big\{ U_{q} \neq \widetilde{W}_{q} \text{ for some } q \in I_{k} \Big\}.
\end{equation}
Then the inequality  \eqref{0714-c-2}  guarantees that 
\begin{equation}\label{0715-b-1}
P (F) \le C 2^{- k^{19}},
\end{equation}
for all $n \ge N_{k}$.

Thus, by H\"older's inequality, we have 
\begin{equation}\label{12}
\parallel J_{3} - J_{4} \parallel _{2} \leq \big( \parallel J_{3}  \parallel _{3} + \parallel J_{4}  \parallel _{3} \big)^{\frac{1}{2}} \cdot P (F)^{\frac{1}{3}} \le  C 2^{-\frac{1}{3} \cdot  k^{19}} \big( \parallel J_{3}  \parallel _{3} + \parallel J_{4}  \parallel _{3} \big)^{\frac{1}{2}}.
\end{equation}
Here $\parallel \cdot \parallel _{3}$ stands for the $L^{3}$-norm.
However, we have  that 
$$
|J_{3} | = \Big|  \sum_{q \in I_{k}}  \phi_{j} ( x_{q} ) \zeta(k) U_{q} \Big|  \le C L_{k} M_k \zeta(k)  \le C 2^{3 k} 2^{3 (k^{4} -k ) } 2^{- \beta k^{4} } \le C 2^{2k^{4}}.
$$
 Similarly, we have
\begin{equation*}
|J_{4}| \le C 2^{2k^{4}}.
\end{equation*}
Thus, we have 
\begin{equation*}
\parallel J_{3}  \parallel _{3} + \parallel J_{4}  \parallel _{3}  \le C 2^{2k^{4}}.
\end{equation*}
Combining this with \eqref{12}, it follows that 
\begin{equation}\label{13}
\parallel J_{3} - J_{4} \parallel _{2} \le C_{5} 2^{-c_{6}  k^{18}}.
\end{equation}

Changing $k = k_{m}$ if necessary, we may assume that 
\begin{equation}\label{0715-e-1}
C_{1} 2^{-k} + C_{3} 2^{- c_{4} k^{2}} + C_{5} 2^{- c_{6} k^{18}} \le 2^{-m},
\end{equation}
where $C_{1}$, $C_{3}$, $c_{4}$, $C_{5}$ and $c_{6}$ are universal constants given in \eqref{3}, \eqref{9} and \eqref{13}. 
With this in mind, using \eqref{3}, \eqref{6}, \eqref{9} and \eqref{13},
it follows that for all $m \ge 1$ and $n_{l} \ge N_{k_m}$
\begin{equation}\label{14}
\parallel \mu_{n_{l}} (\phi_{j}) - J_{4} \parallel _{2} \le C 2^{-m},
\end{equation}
where we recall that $J_{4} = J_{4}^{j, k}$ is a measurable function of ${\cal K}$ which depends on $j$ and $m$ but does not depend on $n_{l}$, see \eqref{macro} for $J_{4}$. Therefore, the inequality of \eqref{14} implies that the sequence of ${\cal K}$-measurable variables $\{ J_{4}^{j, k_m} \}_{m \ge 1}$ is a Cauchy sequence in $L^{2}$. So it has a limit in $L^{2}$. We write 
\begin{equation}\label{15}
\mu (\phi_{j} ) := \lim_{m \to \infty} J_{4}^{j, k_m} \text{ in } L^{2}.
\end{equation}

Note that the limit $\mu (\phi_{j} )$ is of course a measurable function of ${\cal K}$. The inequality \eqref{14} also gives that 
\begin{equation}\label{16}
\parallel \mu (\phi_{j} )  - J_{4}^{j, k_m} \parallel _{2} \le C 2^{-m}.
\end{equation}
Thus, the triangle inequality tells that 
\begin{equation}\label{17}
\parallel \mu (\phi_{j} )  -  \mu_{n_{l}} (\phi_{j})  \parallel _{2} \le C 2^{-m}.
\end{equation}
for all $n_{l} \ge N_{k_m}$. From this, by choosing a sub-sequence $\{ n_{l_{p}} \}_{p \ge 1}$ of $\{ n_{l} \}_{l \ge 1}$ appropriately, we can conclude that $ \mu_{n_{l_{p}}} (\phi_{j}) $ converges to $\mu (\phi_{j} )$  almost surely as $p \to \infty$. However, since  $\mu_{n_{l}}$ converges weakly to $\mu^{\ast}$ almost surely by \eqref{eq:6.2II}, one has $\mu^{\ast} (\phi_{j} ) = \mu (\phi_{j} )$. In particular,  $\mu^{\ast} (\phi_{j} ) $ is a measurable function of ${\cal K}$  for each $j$. This implies that $\mu^{\ast}$ is also measurable with respect to ${\cal K}$. 
\end{proof}

Finally, we are ready to prove our first main theorem.
\begin{proof}[Proof of Theorem \ref{1st}]
By Theorem \ref{goal}, the measure $\mu^{\ast}$ can be characterized uniquely by $\cal K$, which means it does not depend on the chosen sub-sequence $\{ n_{l} \}_{l \ge 1}$. Letting $\mu=\mu^{\ast}$, we obtain that $(\gamma_n,\mu_n)$ converges weakly to $({\cal K},\mu)$, which finishes the proof.
\end{proof}

\section{Convergence in the natural parametrization}\label{sec:7}
In this section, we finally prove Theorem \ref{3rd}, the weak convergence of the naturally rescaled $\{ \eta_{n} \}_{n \ge 1}$ with respect to the supremum distance. 

In Section \ref{sec:guidepara}, we will first give a sketch of the proof.  In Section \ref{sec:regularity} we will establish some regularity properties of LERW. Then in Section \ref{property}, we establish some basic properties of the limiting occupation measure $\mu$ which will be used to parametrize the scaling limit. We will finish the proof of Theorem \ref{3rd} in Section \ref{Natural} by showing that any sub-sequential limit of $\{ \eta_{n} \}_{n \ge 1}$ coincides with the curve parametrized by $\mu$. 

\subsection{Guideline for parametrization of ${\cal K}$ and uniform convergence of $\eta_{n}$}\label{sec:guidepara}

In this subsection, we will give a guideline for the parametrization works we will be carrying out in Sections \ref{property} and \ref{Natural}.

We first make a convention here: to distinguish elements between $ {\cal C}  (\overline{\mathbb{D}})$ and  ${\cal H}  (\overline{\mathbb{D}})$, we will write $\widehat{\lambda} = \{ \lambda (t) \ | \ 0 \le t \le t_{\lambda} \}$ for the range of $\lambda$ if $\lambda \in {\cal C}  (\overline{\mathbb{D}})$ where $t_{\lambda}$ stands for the time duration of $\lambda$. Note that naturally $\widehat{\lambda} \in {\cal H}  (\overline{\mathbb{D}})$. 

By the tightness of $\eta_n$ and the occupation measure $\mu_n$ (see Corollary \ref{0626} and Proposition \ref{tight} resp.), we can find a sub-sequence $\{ n_{l} \}_{l \ge 1 }$, a random element $\zeta$ of $ {\cal C} (\overline{\mathbb{D}} )$ and a random element of $\nu$ of $ {\cal M} (\overline{\mathbb{D}} )$ so that 
\begin{equation}\label{AA}
(\eta_{n_{l}}, \mu_{n_{l}} ) \xrightarrow{d} (\zeta, \nu) \ \  (\text{as } l \to \infty),
\end{equation}
with respect to the product topology of $ \big(  {\cal C} (\overline{\mathbb{D}} ), \rho \big)$ and the topology of the weak convergence on $ {\cal M} (\overline{\mathbb{D}} )$. 
By Theorem \ref{goal}, 
\begin{equation}\label{B}
(\widehat{\eta}_{n_{l}}, \mu_{n_{l}} ) \xrightarrow{d} ({\cal K}, \mu) \ \  (\text{as } l \to \infty),
\end{equation}
with respect to the topology of $ {\cal H} (\overline{\mathbb{D}} ) \times {\cal M} (\overline{\mathbb{D}} )$.  
This gives that 
\begin{equation}\label{C}
(\widehat{\zeta}, \nu)  \overset{d}{=} ({\cal K}, \mu)
\end{equation}
as a random element of the product space ${\cal H} (\overline{\mathbb{D}} )  \times   {\cal M} (\overline{\mathbb{D}} )$. From now on we no longer distinguish $(\widehat{\zeta}, \nu) $ and $({\cal K}, \mu)$  and use $(\widehat{\zeta},\mu)$ throughout this section.
By Skorokhod's representation theorem, we can couple $\{ (\eta_{n_{l}}, \mu_{n_{l}} ) \}_{l \ge 1}$ and $(\zeta, \nu)$ in the same probability space (which we also denote by $P$) such that 
\begin{align}\label{eq:subcoup}
&\mbox{ $\displaystyle \lim_{l \to \infty } \rho (\eta_{n_{l}} , \zeta ) = 0$ almost surely, and} \notag \\
&\mbox{ $\mu_{n_{l}}\to \mu$ with respect to the topology of weak convergence on $ {\cal M} (\overline{\mathbb{D}} )$ almost surely.}
\end{align}

In the next subsection, we will show the following two basic properties of $\mu$ from \eqref{AA}: 
\begin{itemize}
\item[(a)] The support of $\mu$ satisfies
\begin{equation}\label{D}
\text{supp} (\mu) = \widehat{\zeta}
\end{equation}
almost surely. In particular, almost surely 
\begin{equation}\label{E}
\mu \big( \widehat{\zeta}_{x} \big) < \mu \big( \widehat{\zeta}_{y} \big) 
\end{equation}
as long as two points $x, y \in \widehat{\zeta}$ satisfy $\widehat{\zeta}_{x} \subsetneq \widehat{\zeta}_{y}$. Here $\widehat{\zeta}_{x}$ stands for the range of the curve between the origin and $x$ lying in $\widehat{\zeta}$ for $x \in \widehat{\zeta}$ (thus, we regard  $\widehat{\zeta}_{x}$ as an element of ${\cal H} (\overline{\mathbb{D}} )$). Note that by the fact that ${\cal K}$ is a.s.\ a simple curve, for each $x \in \widehat{\zeta}$, $\widehat{\zeta}_{x} \subset \widehat{\zeta}$ is (the range of) the simple curve whose starting point is the origin and endpoint is $x$.

\item[(b)] Almost surely for each $x \in \widehat{\zeta}$,
\begin{equation}\label{F}
\lim_{y \in \widehat{\zeta} \atop  |x-y| \to 0} \mu \big( \widehat{\zeta}_{y} \big) = \mu \big( \widehat{\zeta}_{x} \big).
\end{equation}

\end{itemize}

Once we obtain the properties of $\mu$ as above, we can parametrize $\widehat{\zeta}$ using the measure $\mu$. Note that for each $t \in [0, \mu ( \widehat{\zeta} ) ]$, there is a unique point $x_{t} \in \widehat{\zeta} $ satisfying that $ \mu \big( \widehat{\zeta}_{x_{t}} \big) = t$. Therefore, we may define
\begin{equation}\label{G}
\eta^{\ast} (t) := x_{t}   \text{ for } t \in  [0, \mu ( \widehat{\zeta} ) ].
\end{equation}
Then the properties (a) and (b) as described above ensure that $\eta^{\ast}$ is a random element of ${\cal C} (\overline{\mathbb{D}}) $ which is a measurable function of $(\widehat{\zeta}, \mu) $.

In Section \ref{Natural}, we will show that 
\begin{equation}\label{K}
\zeta{=}\eta^{\ast}\quad\mbox{ a.s., and moreover }\quad(\eta_{n_{l}}, \mu_{n_{l}} ) \xrightarrow{d} (\eta^{\ast}, \mu) \ \  (\text{as } l \to \infty),
\end{equation}
with respect to the product topology of $ \big(  {\cal C} (\overline{\mathbb{D}} ), \rho \big)$ and the topology of the weak convergence on $ {\cal M} (\overline{\mathbb{D}} )$. Recall that $\eta^{\ast}$ is a measurable function of $\mu$ and does not depend on the choice of the sub-sequence. This gives Theorem \ref{3rd}. 

We will prove (a) and (b) (i.e., \eqref{D} and \eqref{F}) in  Section \ref{property} (see Proposition \ref{basic2} and Proposition \ref{basic3}) and then show \eqref{K} in Section \ref{Natural}. 

From now on, in this section we always assume the coupling of  $\{ (\eta_{n_{l}}, \mu_{n_{l}} ) \}_{l \ge 1}$ and $(\zeta, \mu)$ on which \eqref{eq:subcoup} holds.

\subsection{Regularity controls for LERW}\label{sec:regularity}
In this section, we give a few regularity controls for LERW crucial for the analysis of the (sub-sequential) scaling limit in Sections \ref{property} and \ref{Natural}.

We start with a partitioning of the space.

\begin{dfn}\label{def:minimal}
Given $\epsilon\in (0,1)$, we decompose $\mathbb{R}^{3}$ into a collection of closed cubes $\{ B_{i} \}_{i=1}^{\infty} = \{ B_{i}^{\epsilon} \}_{i=1}^{\infty}$ with side length $\epsilon$ which satisfies $x_{i} \in \frac{\epsilon}{2}\mathbb{Z}^3$ and $x_{i} \neq x_{j} $ if $i \neq j$, where $x_{i}$ stands for the center of $B_{i}$.
We assume $x_{1} = 0$. We further assume that $\{ B_{i} \}_{i=1}^{M} = \{ B_{i}^{\epsilon} \}_{i=1}^{M}$ is a covering of $\overline{\mathbb{D}}$ such that $\overline{\mathbb{D}} \subset \bigcup_{i=1}^{M} B_{i}$, $\overline{\mathbb{D}} \cap B_{i} \neq \emptyset$ for all $1 \le i \le M$ and $B_{i} \cap B_{i'} \cap 2^{-n} \mathbb{Z}^{3} = \emptyset $ for $i \neq i'$. Note that  $M = M_{\epsilon} \asymp \epsilon^{-3}$. 
\end{dfn}
Let $X^{n, \epsilon}_{i}$ be the number of points in $2^{-n} \mathbb{Z}^{3} \cap B_{i}$ hit by $\eta_{n}$.
Also, we define $$I = \big\{ 1 \le i \le M \ \big| \  \text{dist} \big( B_{i}, \partial \mathbb{D} \big) \ge 10 \epsilon \big\}.$$ 
We define the event $V_{n} = V_{n,\epsilon}$ 
as follows:
\begin{equation}\label{0705-3-1}
V_{n} = V_{n,\epsilon} := \Big\{ X^{n, \epsilon}_{i} \le \sqrt{\epsilon} 2^{\beta n} \text{ for all } i \in I \Big\}.
\end{equation}
By an argument essentially the same as that of \cite[Theorem 8.4]{S} (see also \cite[Lemma 4.6]{Holder}), there exists a universal constant $C_{0} < \infty$ such that for all $k \ge 1$
\begin{equation}\label{0702-a-1}
E \Big\{ \big( X^{n, \epsilon}_{i} \big)^{k} \Big\} \le C_{0}^{k} k !  \big( \epsilon 2^{n} \big)^{k \beta},
\end{equation}
for $i\in I$.
Letting $c_{1} = \frac{1}{2 C_{0}}$, it follows from \eqref{0702-a-1} that for all $i\in I$,
\begin{equation}\label{0702-a-2}
E \Big\{ \exp \Big( \frac{ c_{1} X^{n, \epsilon}_{i}  }{ \epsilon^{\beta } 2^{\beta n } }\Big) \Big\} < \infty.
\end{equation}
Using this, Markov's inequality implies that for all $\kappa \ge 1$ and $i\in I$.
\begin{equation}\label{exp-tail-up}
P \Big( X^{n, \epsilon}_{i} > \kappa \big( \epsilon 2^{n} \big)^{\beta} \Big) \le C e^{- c \kappa }
\end{equation}
for some universal constants $0 < c, C < \infty $. 
It then follows that 
\begin{equation}\label{card}
P \big( V_{n, \epsilon} \big) \ge 1- C e^{- \frac{c}{\sqrt{\epsilon}} },
\end{equation}
for some universal constants $0 < c, C< \infty$.

We  also want to control the time that $\eta_{n}$ stays in 
\begin{equation}\label{eq:Hbdef}
 H_{b}: = \mathbb{D}_{1+b}\setminus \mathbb{D}_{1-b}
\end{equation}
for some $b\in(0,1/10)$. By first moment bound from Section \ref{sec:ONE}, there exists $C_3<\infty$ that
\begin{equation}\label{0629-ver-2}
P \big( \mu_{n} ( H_{b} ) \ge r \big) \le \frac{C_{3} b^{\beta }}{r}
\end{equation}
for all $r > 0$. Let $u (b)  = u_{n} (b) = \inf \big\{ j \ge 0 \ \big| \ \eta_{n} (j) \in H_{b} \big\}$ be the first time that $\eta_{n}$ hits $H_{b}$. It follows from \eqref{srwbound} that there exists $C_4<\infty$ that
\begin{equation}\label{0703-a-1}
P \big( \eta_{n}  [u(b^{2} ), t_{\eta_n} ] \not\subset H_{b} \big) \le C_{4} b.
\end{equation}
Thus, if we define 
\begin{equation}\label{0703-b-2}
 U_{n, b, r} := \Big\{ \mu_{n} ( H_{b} ) \le r, \ \eta_{n}  [u(b^{2} ), t_{\eta_n} ] \subset H_{b} \Big\},
\end{equation}
we have that 
\begin{equation}\label{0703-b-3}
P \big( U_{n, b, r} \big) \ge 1 -  \frac{C_{3} b^{\beta }}{r} - C_{4} b.
\end{equation}


We next recall the event $K_{n,\delta,\epsilon}$ from Proposition \ref{HC-upper} controlling the ``modulus of continuity'' of $\eta_n$.
To simplify notation, we now fix a certain $\beta'\in(\beta,2)$ and write
\begin{equation}\label{eq:Wndef}
W_{n,a}:=K_{n,a,1/\beta-1/\beta'}
\end{equation}
for $a\in (0,1/10)$.
We also recall the definition of quasi-loops from  Section \ref{qloops}.  Write
\begin{equation}\label{0705-5}
 F_{n} =F_{n,\epsilon} := \Big\{ {\rm QL} \big( \epsilon^{L}, \epsilon ; \eta_{n} \big) = \emptyset \Big\}.
\end{equation} 
where $L$ is the quantity defined in Lemma \ref{lem:noquasiloop} for the event that there is no $(\epsilon ,\epsilon^{L})$-quasi-loop for $\eta_n$.

We now summarize all the regularity properties we need in the following proposition:
\begin{prop}\label{prop:regularity}
For any $0<c<1$, there exist $a_0,\epsilon_0,b_0,r_0\in(0,\infty)$ such that for any $0<a<a_0$, $0<\epsilon<\epsilon_0$, $0<b<b_0$, $0<r<r_0$ and any sub-sequence $(n_l)_{l\geq 1}$ 
\begin{equation}\label{eq:regularity}
P \Big( W_{n_l,a} \cap V_{n_l,\epsilon} \cap F_{n_l,\epsilon} \cap U_{n_l, b, r} \text{ occurs for infinitely many } l \Big) \ge 1 - c.
\end{equation}
\end{prop}
\begin{proof}
The claim follows from  \eqref{eq:noquasiloop}, \eqref{card}, \eqref{eq:HC-upper}, \eqref{0703-b-3} and the reverse Fatou lemma. 
\end{proof}

Before finishing this subsection, we will show that any sub-sequential scaling limit $\zeta : [0, t_{\zeta}] \to \widehat{\zeta}$ defined in \eqref{AA} is a bijection. Here we recall that $\widehat{\zeta} = \{ \zeta (t) \ | \ 0 \le t \le t_{\zeta} \}$ stands for the range of $\zeta$.
\begin{prop}\label{basic4}
The mapping $\zeta : [0, t_{\zeta}] \to \widehat{\zeta}$ is a bijection almost surely.
Furthermore, 
$\zeta (0) = 0$, $\zeta ( t_{\zeta} ) \in \partial \mathbb{D}$, $\zeta (t) \in \mathbb{D}$ for all $t \in [0, t_{\zeta} )$ and $t_{\zeta} \in (0, \infty)$ almost surely.
\end{prop}
\begin{proof}

We start with the first claim. It suffices to show that $\zeta : [0, t_{\zeta}] \to \widehat{\zeta}$ is an injection almost surely since it is a trivially surjection. We will prove it by contradiction. To do it, define
\begin{equation*}
R_{0} := \sup \Big\{ |t-t'| \ \Big| \ 0 \le t, t' \le t_{\zeta}, \ \zeta (t) = \zeta (t') \Big\}.
\end{equation*}
We have that $\zeta : [0, t_{\zeta}] \to \widehat{\zeta}$ is an injection if and only if $R_{0} = 0$. So suppose that 
\begin{equation*}
P (R_{0} > 0) \ge c_{1} > 0
\end{equation*}
for some positive constant $c_{1}$. Then we can find $\delta > 0$ such that 
\begin{equation}\label{as1}
P (R_{0} \ge \delta) \ge \frac{c_{1}}{2}.
\end{equation}

By Proposition \ref{prop:regularity}, one can choose $0<\epsilon,b,r<1/10$ such that 
$\epsilon < 
\min(\delta^2/10^5, b^2/ 10)$, $r<\delta/40$ and the event 
$$
 \{R_{0} \ge \delta\}\cap \{F_{n_{l_j},\epsilon} \cap V_{n_{l_{j}},\epsilon} \cap U_{n_{l_{j}}, b, r}\mbox{ occurs for a sub-sub-sequence } n_{l_{j}},\;j = 1, 2, \cdots \}   
 $$
happens with probability at least $\frac{c_{1}}{4}$. In this proof, from this point onwards we always assume this event.

By the first assertion  $R_{0} \ge \delta$, we can find $0 \le t < t' \le t_{\zeta}$ with $ t' - t > \frac{\delta}{2}$ such that $\zeta (t) = \zeta (t')$. 
Taking $j$ sufficiently large such that $\rho (\eta_{n_{l_{j}}} , \zeta ) < 10^{-4} \epsilon^{L} \cdot \min \{ 1, t_{\zeta} \}$ (recall that the constant $L$ is defined in Lemma \ref{lem:noquasiloop}), we can find $0 \le s < s' \le t_{\eta_{n_{l_{j}}}}$ such that the following conditions hold:
\begin{itemize}
\item[(i)] $s' - s > \frac{\delta}{4}$;

\item[(ii)] $\big| \eta_{n_{l_{j}}} (s) - \zeta (t) \big| \le \frac{\epsilon^{L}}{100}$ \ \  and \ \  $\big| \eta_{n_{l_{j}}} (s') - \zeta (t') \big| \le \frac{\epsilon^{L}}{100}$,
\end{itemize}
Because $\zeta (t) = \zeta (t')$, the second condition (ii) ensures that 
\begin{equation}\label{as3}
\big| \eta_{n_{l_{j}}} (s) - \eta_{n_{l_{j}}} (s') \big| \le \frac{\epsilon^{L}}{50}.
\end{equation}

\begin{figure}[h]\label{supp-fig-6}
\begin{center}
\includegraphics[scale=0.65]{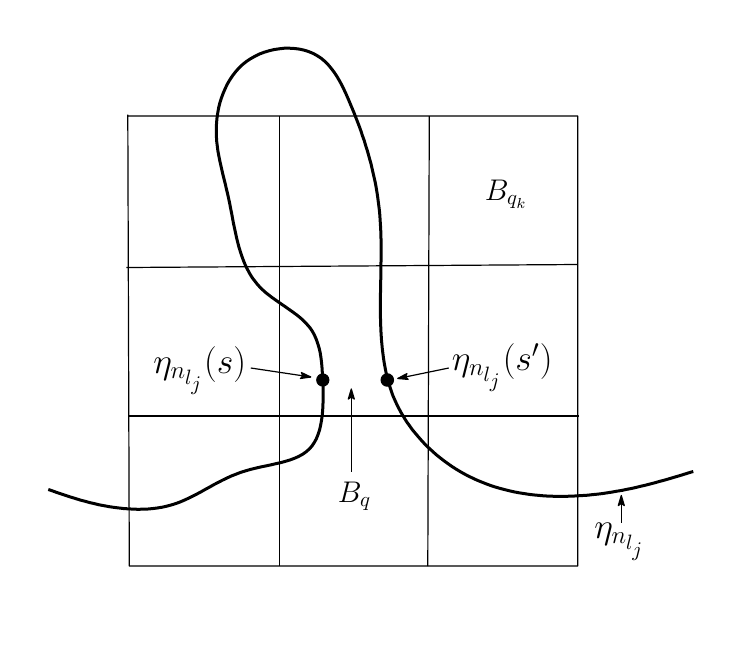}
\caption{Illustration for a quasi-loop at $\eta_{n_{l_{j}}} (s)$. The point $\eta_{n_{l_{j}}} (s)$ is contained in $B_{q}$. The cube $B_{q}$ is surrounded by $\{ B_{q_{k}} \}_{k=2}^{27}$.}\label{ql-q-2}
\end{center}
\end{figure}

Recall the covering $\{ B_{i}^{\epsilon}\}_{i=1}^{M}$
from Definition \ref{def:minimal}. 
Let $B_{q}$ be the cube which contains $\eta_{n_{l_{j}}} (s)$. Recall the definition of  $H_{b}$ 
in \eqref{eq:Hbdef}. We first claim that 
\begin{equation}\label{0704}
B_{q} \cap H_{b_{1}^{2}/3} = \emptyset.
\end{equation}  
Otherwise, by our choice of $\epsilon$, $\eta_{n_{l_{j}}} (s) \in H_{b^{2}/2} \subset H_{b}$, which contradicts with the event $U_{n_{l_{j}}, b, r}$ because $\mu_{n_{l_{j}}} (H_{b} ) \ge s' - s > \delta/ 4 > 10 r$ by our choice of $r$.

Therefore, we get \eqref{0704} and the index $q$ is an element of $I$ in Definition \ref{def:minimal}.
by our choice of $\epsilon$. 

We take $B_{q} = B_{q_{1}}, B_{q_{2}}, \cdots , B_{q_{27}}$ from our covering $\{ B_{i} \}_{i=1}^{M} $
which satisfies that 
$B_{q_{k}} \cap B_{q} \neq \emptyset$ for each $k =2, \cdots , 27$. Note that by our choice of $\epsilon$ 
we have  $q_{k} \in I$ for all $1 \le k \le 27$.

Since the event $V_{n_{l_{j}},\epsilon}$ occurs, it follows that 
\begin{equation}\label{as4}
\eta_{n_{l_{j}}} [s, s'] \not\subset \bigcup_{k=1}^{27} B_{q_{k}}.
\end{equation}
Otherwise, one has
\begin{equation*}
X^{n_{l_{j}}, \epsilon}_{q_{1}} + X^{n_{l_{j}}, \epsilon}_{q_{2}} + \cdots + X^{n_{l_{j}}, \epsilon}_{q_{27}}  \ge (s'-s) 2^{\beta n_{l_{j}} }  \ge \frac{\delta}{4} \cdot 2^{\beta  n_{l_{j}} },
\end{equation*}
Combining this with our choice of $\epsilon$ 
we see that there exists $k \in \{ 1,2, \cdots , 27 \}$ such that 
\begin{equation*}
X^{n_{l_{j}}, \epsilon}_{q_{k}} \ge \frac{\delta}{108} \cdot 2^{\beta n_{l_{j}} } > \sqrt{\epsilon} 2^{\beta n_{l_{j}} },
\end{equation*}
which contradicts the event $V_{n_{l_{j}},\epsilon}$. Thus, \eqref{as4} holds. However, using \eqref{as3} and \eqref{as4}, we see that $ \eta_{n_{l_{j}}} $ has a $(\frac{\epsilon^{L}}{50}, \epsilon )$-quasi-loop at $\eta_{n_{l_{j}}} (s)$; see Figure \ref{ql-q-2}. This contradict with the event $F_{n_{l_{j}},\epsilon}$.
So, we see that $\zeta$ is a bijection almost surely.

The second claim is easy to prove. 
Since $\eta_{n_{l}}$ converges to $\zeta$ with respect to the metric $\rho$ almost surely, we have $\zeta (0) = \lim_{l \to \infty} \eta_{n_{l}} (0) = 0$, $\zeta (t_{\zeta} ) = \lim_{l \to \infty} \eta_{n_{l}} \big( t_{n_{l}} \big) \in \partial \mathbb{D}$ and $\zeta [0, t_{\zeta} ] \subset \overline{\mathbb{D}}$ almost surely. We already know that $\zeta$ is a bijection, hence, if there exists $t \in [0, t_{\zeta} )$ such that $\zeta (t) \notin \mathbb{D}$, it holds that $\zeta (t) \in \partial \mathbb{D}$ with $\zeta (t) \neq \zeta (t_{\zeta} )$. But this contradicts \cite[Theorem 1.2]{SS} which ensures that $ \widehat{\zeta} \cap \partial \mathbb{D}$  is a singleton almost surely. It is clear that $t_{\zeta} > 0$ because  $\zeta (0) = 0$ while $\zeta ( t_{\zeta} ) \in \partial \mathbb{D}$. Finally, by \eqref{eq:fnasymp} and \eqref{exp-tail}, it follows 
that there exist universal constants $c >0 $ and $C < \infty$ that for all $\kappa \ge 1$
\begin{equation*}
\sup_{n \ge 1}  P \big( t_{\eta_n} \ge  \kappa \big) \le C e^{ - c \kappa },
\end{equation*}
Combining this with the fact that $t_{\eta_{n_{l}}} \to t_{\zeta}$ with probability one, we have $t_{\zeta} < \infty$ almost surely.
Thus, we obtain the second claim and finish the proof.
\end{proof}

\subsection{Properties of $\mu$}\label{property}
In this subsection we will show several basic properties of $\mu$ of Theorem \ref{1st}. 

\begin{prop}\label{basic1}
Take a cube $B \subset \mathbb{D}$ with  $0 \notin \partial B$. It follows that $\mu ( \partial B) = 0$ almost surely (see \eqref{AA} for the measure $\mu$). Furthermore, we have $\mu (    \partial \mathbb{D} ) = 0$ with probability one.
\end{prop}
The proof is an easy application of Markov inequality with the upper bound of the one-point function (see Section \ref{sec:ONE} for details). Hence we omit it.

\begin{rem}
As a quick corollary of Proposition \ref{basic1}, for any cube $B \subset \mathbb{D}$ with $0 \notin \partial B$, 
\begin{equation}\label{darui-5}
\lim_{l \to \infty } \mu_{n_{l}} (B)  = \mu (B)
\end{equation}
almost surely. From this, if we consider a (at most countable)  collection of cubes $B_{1}, B_{2} , \cdots $ with $0 \notin \partial B_{i}$ for each $i \ge 1$, writing $B' = B_{1} \cup B_{2} \cup \cdots $, it follows that 
\begin{equation}\label{darui-6}
\lim_{l \to \infty } \mu_{n_{l}} (B')  = \mu (B')\quad\mbox{ almost surely.}
\end{equation}
 Furthermore, we have 
\begin{equation}\label{darui-7}
\lim_{l \to \infty } \mu_{n_{l}} ( \mathbb{D})  = \mu (\mathbb{D})\quad\mbox{ almost surely.}
\end{equation}

\end{rem}

Next we will deal with the support of the measure $\mu$. 

\begin{prop}\label{basic2}
It follows that 
\begin{equation}\label{darui-8}
{\rm supp} (\mu) = \widehat{\zeta}\quad \mbox{almost surely},
\end{equation}
see \eqref{AA} for $\mu$ and $\zeta$. Also, recall that $\widehat{\zeta} = \{ \zeta (t) \ | \ 0 \le t \le t_{\zeta} \}$ stands for the range of $\zeta$.
\end{prop}

\begin{proof}
%
We will first prove that $\text{supp} (\mu) \subset \widehat{\zeta} $ almost surely. It suffices to prove that with probability one $\mu ( F) = 0$ for each closed set $F \subset \overline{\mathbb{D}} \setminus \widehat{\zeta}$. So take a closed subset $F \subset \overline{\mathbb{D}} \setminus \widehat{\zeta}$. Taking $\epsilon > 0$ sufficiently small,  we see that $F_{\epsilon} :=  \{ x \in \mathbb{R}^{3} \ | \ \text{dist} (x, F )  < \epsilon \}$ also has a positive distance from $\widehat{\zeta}$. Since $ \lim_{l \to \infty }  \rho (\eta_{n_{l}} , \zeta ) = 0$ with probability one, it follows that $F_{\epsilon} \cap \widehat{ \eta}_{n_{l}} = \emptyset$ for sufficiently large $l$, which implies that $\mu_{n_{l}} (F_{\epsilon}) = 0$. Since $F_{\epsilon}$ is open, it follows that $\mu ( F_{\epsilon}) = 0$, which gives $\mu (F) =0$ with probability one. Thus, we see that  $\text{supp} (\mu) \subset \widehat{\zeta} $ almost surely.

We will next prove that $\text{supp} (\mu) = \widehat{\zeta} $ almost surely. For a point $x \in \widehat{\zeta}$, we write $\widehat{\zeta}_{x}$ for the range of the sub-path of $\widehat{\zeta}$ between the origin and $x$. Note that $\widehat{\zeta}_{x}$ is a random element of ${\cal H} (\overline{\mathbb{D}} )$ which is a simple path for each $x \in \widehat{\zeta}$ (see Section \ref{SCALING} for this fact). In order to prove  $\text{supp} (\mu) = \widehat{\zeta} $, it suffices to show that for each $x, y \in \widehat{\zeta}$ with $\widehat{\zeta}_{x} \subsetneq \widehat{\zeta}_{y}$ we have $\mu \big( \widehat{\zeta}_{y} \setminus \widehat{\zeta}_{x} \big) > 0$. 

We will show this by contradiction. Let 
\begin{equation*}
D_{0} = \sup \Big\{ |x-y| \ \Big| \ x, y \in \widehat{\zeta} \text{ with } \widehat{\zeta}_{x} \subseteq \widehat{\zeta}_{y}, \ \mu \big( \widehat{\zeta}_{y} \setminus \widehat{\zeta}_{x} \big)  = 0 \Big\}.
\end{equation*}
Suppose that there exist  $x, y \in \widehat{\zeta}$ with $\widehat{\zeta}_{x} \subsetneq \widehat{\zeta}_{y}$ such that $\mu \big( \widehat{\zeta}_{y} \setminus \widehat{\zeta}_{x} \big) = 0$. Then we have $D_{0} > 0$. So suppose that 
\begin{equation}\label{darui-9}
P ( D_{0} > 0 )  \ge c_{0} > 0.
\end{equation}
Taking $q > 0$ sufficiently small, we have 
\begin{equation}\label{darui-10}
P ( D_{0} \ge q  )  \ge \frac{c_{0}}{2}.
\end{equation}

By \eqref{darui-10} and Proposition \ref{prop:regularity}, 
we can pick some small $0<a < (q/2)^{\beta'}$ (recall the choice of $\beta'$ above \eqref{eq:Wndef}), such that the event
$$
 \{ D_{0} \ge q \} \cap \{W_{n_{l_{j}}, a } \mbox{ occurs for a sub-sub-sequence } n_{l_{j}},\;j = 1, 2, \cdots \}   
$$
happens with probability at least $c_{0}/4$.
We work on this event from now on. 
In this case, there exist $x, y  \in \widehat{\zeta} \text{ with } \widehat{\zeta}_{x} \subseteq \widehat{\zeta}_{y}$ such that $\mu \big( \widehat{\zeta}_{x, y} \big)  = 0$ and $|x-y| \ge \frac{ 2 q}{3}$. Here $\widehat{\zeta}_{x, y}$ stands for the range of the sub-path of $\zeta$ between $x$ and $y$ which is a closed set, i.e., $\widehat{\zeta}_{x, y } = \big( \widehat{\zeta}_{y} \setminus \widehat{\zeta}_{x} \big) \cup \{ x \}$. We fix such two points $x, y \in  \widehat{\zeta}$. 

Take $\epsilon \in (0,1)$. Recall that $\{ B_{i} \}_{i = 1}^{M} = \{ B_{i}^{\epsilon} \}_{i=1}^{M_{\epsilon}}$ is a covering of $\overline{\mathbb{D}}$ as defined in Definition \ref{def:minimal}. We set 
\begin{equation*}
I_{\epsilon}^{x,y} := \Big\{ 1 \le i \le M \ \Big| \ B_{i} \cap \widehat{\zeta}_{x,y} \neq \emptyset \Big\}  \ \ 
\mbox{ and }
 \ \ J_{\epsilon}^{x,y} := \Big\{ 1 \le i \le M \ \Big| \ B_{i} \cap B_{j} \neq \emptyset \text{ for some } j \in I_{\epsilon}^{x,y} \Big\}.
\end{equation*}
Note that
\begin{equation*}
{\cal B}_{\epsilon}^{x,y} :=  \displaystyle \bigcup_{i \in J_{\epsilon}^{x,y}} B_{i}
\end{equation*}
is a covering of $\widehat{\zeta}_{x,y} $ satisfying that ${\cal B}_{\epsilon}^{x,y} \downarrow \widehat{\zeta}_{x,y}$ as $\epsilon \downarrow 0$. We also mention that
the set of points within distance $\frac{\epsilon}{2}$ from $\widehat{\zeta}_{x,y}$ is contained in  ${\cal B}_{\epsilon}^{x,y}$.

Since  $\mu \big( \widehat{\zeta}_{x, y} \big)  = 0$ and  ${\cal B}_{\epsilon}^{x,y} \downarrow \widehat{\zeta}_{x,y}$ as $\epsilon \downarrow 0$, using the monotone convergence theorem, we have $\mu \big( {\cal B}_{\epsilon}^{x,y}  \big) \downarrow 0$ as $\epsilon \downarrow 0$. Thus, we can take $\epsilon_{0} \in  (0, \frac{q}{100})$ sufficiently small so that $\mu \big( {\cal B}_{\epsilon_{0}}^{x,y}  \big) \le  \frac{a}{10}$. 
Furthermore, since $\mu_{n_{l_{j}}}$ converges weakly to $\mu$  and $\mu \big( \partial {\cal B}_{\epsilon_{0}}^{x,y}  \big) = 0$ almost surely by Proposition \ref{basic1}, taking $j$ sufficiently large, it follows that
\begin{equation}\label{atama}
\mu_{n_{l_{j}}} \big(  {\cal B}_{\epsilon_{0}}^{x,y}  \big)  \le \frac{a}{5}.
\end{equation}

Since $\displaystyle \lim_{j \to \infty } \rho (\eta_{n_{l_{j}}} , \zeta ) = 0$, taking $j$ sufficiently large, it follows that there exist $s < s'$ such that the following two conditions hold (see Figure \ref{supp-fig-5} for this):
\begin{figure}[h]
\begin{center}
\includegraphics[scale=0.6]{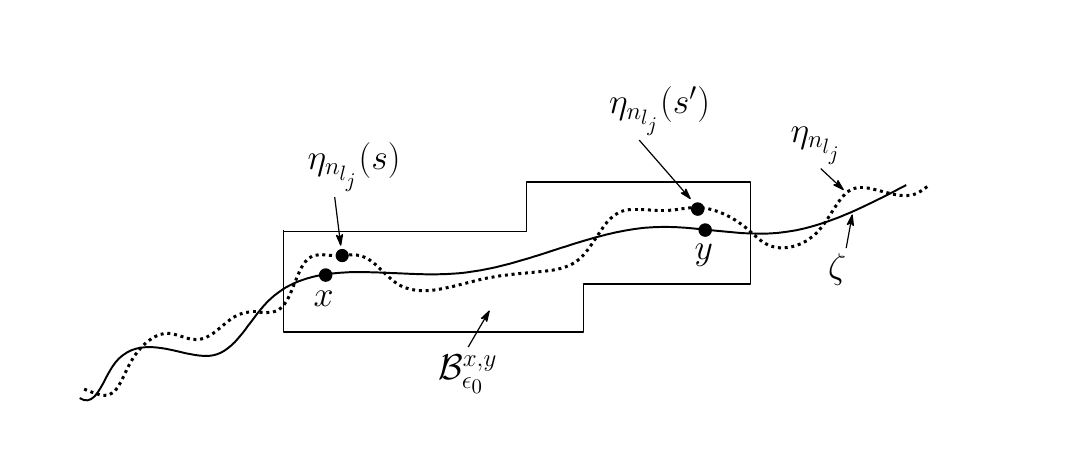}
\caption{Illustration for $\eta_{n_{l_{j}}}$ and $\zeta$. The solid curve stands for $\zeta$ while the dotted curve represents $\eta_{n_{l_{j}}}$.} \label{supp-fig-5}
\end{center}
\end{figure}
\begin{itemize}
\item $\eta_{n_{l_{j}}}[ s, s'] \subset  {\cal B}_{\epsilon_{0}}^{x,y} $;

\item $|\eta_{n_{l_{j}}} (s) -x | \le \frac{\epsilon_{0}}{8}, \ |\eta_{n_{l_{j}}} (s') -y | \le \frac{\epsilon_{0}}{8}$.
\end{itemize}

Since $|x-y| \ge \frac{2 q}{3}$ and $\epsilon_{0} < \frac{q}{100}$, using the triangle inequality
\begin{equation*}
|\eta_{n_{l_{j}}} (s) - \eta_{n_{l_{j}}} (s')  | \ge |x-y| - \frac{\epsilon_{0}}{4} \ge \frac{q}{2} > a^{1/\beta'}.
\end{equation*}
However, the event $W^{n_{l_{j}},a}$ ensures that $s' - s > a$.
Combining this with the fact that $\eta_{n_{l_{j}}}[ s, s']  \subset  {\cal B}_{\epsilon_{0}}^{x,y}$, we have
\begin{equation}\label{atama-1}
\mu_{n_{l_{j}}} \big(  {\cal B}_{\epsilon_{0}}^{x,y}  \big) \ge   \mu_{n_{l_{j}}} \big( \eta_{n_{l_{j}}}[ s, s'] \big) \ge s' - s  > a.
\end{equation}
which contradicts \eqref{atama}. Thus, we conclude that $\text{supp} (\mu) = \widehat{\zeta}$ almost surely and finish the proof.
\end{proof}

Recall that we write $\widehat{\zeta}_{x}$ for the range of the sub-path of $\eta$ between the origin and $x \in \widehat{\zeta}$. The next proposition shows that $\mu ( \widehat{\zeta}_{x} )$ is a continuous function in $x \in \widehat{\zeta}$.

\begin{prop}\label{basic3}
Almost surely,
\begin{equation}\label{singleton}
    \mu(\{x\})=0,\quad \forall x \in \widehat{\zeta}.
\end{equation}
In particular, for each $x \in \widehat{\zeta}$, 
\begin{equation}\label{conti-meas}
\lim_{y \in \widehat{\zeta} \atop  |x-y| \to 0} \mu \big( \widehat{\zeta}_{y} \big) = \mu \big( \widehat{\zeta}_{x} \big).
\end{equation}
\end{prop}

\begin{proof}
%
We will show \eqref{singleton} by contradiction. To do it, let 
\begin{equation*}
L_{0} := \sup \{ \mu ( \{ x \} )  \ | \ x \in \widehat{\zeta} \}
\end{equation*}
so that \eqref{singleton} is equivalent to the equation $L_{0} = 0$. So suppose that 
\begin{equation}
P ( L_{0} > 0 ) \ge c_{0} > 0
\end{equation}
for some positive constant $c_{0}$. Then by taking $\delta > 0$ sufficiently small, we have 
\begin{equation}\label{hairi}
P ( L_{0} >  \delta ) \ge \frac{c_{0}}{2}.
\end{equation}
Proposition  \ref{basic1} guarantees that if we take $a > 0$ sufficiently small we have that 
\begin{equation}\label{0702-b-1}
P \big( R_{\delta,a} \big) \ge 1- \frac{c_{0}}{4},
\end{equation}
where
\begin{equation}
R_{\delta,a} :=\bigg\{    \mu \Big( \big\{ y \in \mathbb{R}^{3} \ \big| \ 1-a \le |y| \le 1 \big\} \Big) \le \frac{\delta}{100}\bigg\},
\end{equation}
By Proposition \ref{prop:regularity}, one can choose $0<\epsilon< \min(\delta^2/100,a/100)$ such that the event 
$$
 \{L_{0} > \delta\} \cap R_{\delta,a} \cap \{ V_{n_{l_{j}},\epsilon} \mbox{ occurs for a sub-sub-sequence } n_{l_{j}},\;j = 1, 2, \cdots \}   
 $$
happens with probability at least $c_0/8$. From now on we always assume this event. Note that in this case, 
%
%
there exists $x \in \widehat{\zeta}$ such that $|x| < 1-a$, $\mu ( \{ x \} ) \ge \frac{\delta}{2}$.  Thus, taking a cube $B_{i}$ with $x \in B_{i}$, we have that $i \in I$ since $a > 100 \epsilon$. Also, it holds that $\mu ( B_{i} ) \ge \mu ( \{ x \} ) \ge \frac{\delta}{2}$.
Since $\mu_{n_{l_{j}}}$ converges weakly to $\mu$ and $\mu ( \partial B_{i} ) = 0$ by Proposition \ref{basic1}, it follows that 
\begin{equation}\label{qq1}
\mu_{n_{l_{j}}} (B_{i} ) \ge \frac{\delta}{4}
\end{equation}
for sufficiently large $j$. However, this contradicts with the event $V_{n_{l_{j}},\epsilon}$ which 
 guarantees that
\begin{equation}\label{qq2}
\mu_{n_{l_{j}}} (B_{i} ) \le \sqrt{\epsilon_{0}} < \frac{\delta}{10},
\end{equation}
which contradicts \eqref{qq1}. This finishes the proof of \eqref{singleton}.
\end{proof}

\subsection{Uniqueness of the sub-sequential limit}\label{Natural}
Let us summarize our current standing point. Recall that $\zeta$ stands for the limit of $\eta_{n}$ along some sub-sequence $n_l$, $l\geq 1$. By Propositions \ref{basic2} and  \ref{basic3}, 
we may consider the reparametrized curve $\eta^{\ast}\in {\cal C} (\overline{\mathbb{D}})$ via the measure $\mu$ as in \eqref{G}.  We mention that  $\eta^{\ast} : [0, \mu (\widehat{\zeta} ) ] \to \widehat{\zeta} $ is a bijection which is a measurable function of $(\widehat{\zeta}, \mu) $. 
We now show that 
\begin{prop}\label{end}
\begin{equation}\label{conv-unif-re}
\zeta = \eta^{\ast}\mbox{ almost surely.}
\end{equation}
\end{prop}
This gives \eqref{K} as well as Theorem \ref{3rd} as discussed in Section \ref{sec:guidepara}. 
\begin{proof}

We need to prove that with probability one 
\begin{equation}\label{coincide2}
t_{\zeta} = t_{\eta^{\ast}} \  \ \text{ and }  \ \     \zeta (t) = \eta^{\ast} (t) \  \text{ for all }  \ t \in [0, t_{\zeta}].
\end{equation}

We will first prove that $t_{\zeta} = t_{\eta^{\ast}}$ almost surely. Note that since $\displaystyle \lim_{l \to \infty } \rho (\eta_{n_{l}} , \zeta ) = 0$, it follows that $t_{\eta_{n_{l}}} \to t_{\zeta}$ as $l \to \infty$ almost surely. But, by definition, we have
\begin{equation*}
t_{\eta_{n_{l}}}  =\mu_{n_{l}} (\overline{\mathbb{D}} ).
\end{equation*}
This implies that $\mu_{n_{l}} (\overline{\mathbb{D}} ) \to t_{\zeta}$ as $l \to \infty$ with probability one. On the other hand, we know that $\mu_{n_{l}}$ converges to $\mu$ with respect to the topology of weak convergence on $ {\cal M} (\overline{\mathbb{D}} )$. Combining this with Proposition \ref{basic1}, it follows that $\mu_{n_{l}} (\overline{\mathbb{D}} ) \to \mu  (\overline{\mathbb{D}} ) $. Thus, we have with probability one
\begin{equation*}
t_{\zeta} = \mu  (\overline{\mathbb{D}} ) = \mu (\widehat{\zeta} ) =   \mu ( \widehat{ \eta^{\ast}} ) = t_{ \eta^{\ast}},
\end{equation*}
where in the second equation we used  the fact that the support of $\mu$ coincides with $\widehat{\zeta}$ which is proved in  Proposition \ref{basic2}. Also, we note that the third equality follows from $\widehat{\zeta} = \widehat{\eta^{\ast}}$. 

Throughout this proof, we write
\begin{equation}\label{0726-u}
U:=t_{\zeta} = t_{\eta^{\ast}}.
\end{equation}
We will next prove that 
\begin{equation}\label{essen}
\zeta (t) = \eta^{\ast} (t) \text{ for all } t \in [0, U]
\end{equation}
by contradiction. To do it, we define
\begin{equation}\label{CONTRA}
Z_{0} := \sup \Big\{ |t - t' | \ \Big|  \  t, t' \in [0,U], \    \zeta (t) = \eta^{\ast} (t') \Big\}.
\end{equation}
By Proposition \ref{basic4}, $\zeta: [0, U] \to \widehat{\zeta} $ is a continuous and bijective curve. Propositions \ref{basic2} and \ref{basic3} show that the same is true for $\eta^{\ast}: [0, U] \to  \widehat{\zeta} $. 
Therefore, we have  
\begin{equation}\label{equiva}
Z_{0} = 0 \Leftrightarrow \text{\eqref{essen}}.
\end{equation}

With this in mind, suppose that $P ( Z_{0} > 0 ) \ge c_{0} > 0$ for some positive constant $c_{0}$. We can then find $\delta > 0$ such that 
\begin{equation}\label{CONTRA2}
P ( Z_{0} > \delta ) \ge \frac{c_{0}}{2}.
\end{equation}
Suppose that $Z_{0} > \delta $. Then we can find $t, t' \in [0, U]$ such that $|t- t'| > \frac{\delta}{2}$ and $\zeta (t) = \eta^{\ast} (t')$. Note that $\zeta (t) = \eta^{\ast} (t') $ must lie in $\mathbb{D}$. The reason for this is as follows. Let $x_{0} = \zeta (t) = \eta^{\ast} (t')$ and suppose that $x_{0} \notin \mathbb{D}$. The fact that $\zeta$ is a bijection 
then guarantees that $x_{0} \in \partial \mathbb{D}$ and that 
$\widehat{\zeta}_{x_{0}} = \widehat{\zeta}$. Namely, $x_{0}$ must be the endpoint of $\widehat{\zeta}$ that lies in 
$\partial \mathbb{D}$. This implies that $t= t' = U$ which contradicts $|t- t'| > \frac{\delta}{2}$. Thus, we have $x_{0} = \zeta (t) = \eta^{\ast} (t') \in \mathbb{D}$. Consequently, it holds that 
\begin{equation*}
\{ Z_{0} > \delta \} \subset \bigcup_{p=1}^{\infty} J_{p, \delta},
\end{equation*}
with
\begin{equation}\label{CONTRA4}
J_{p, \delta} := \Big\{ \text{ There exist $t, t' \in [0, U]$ such that $|t- t'| > \frac{\delta}{2}$ and $\zeta (t) = \eta^{\ast} (t') \in   \mathbb{D}_{1 - p^{-1} }$ } \Big\}.
\end{equation} 
Combining this with \eqref{CONTRA2}, it follows that there exist $\delta \in (0, 1)$, $c_{0} \in (0, 1)$ and $p > 1$ such that 
\begin{equation}\label{CONTRA3}
P \big( J_{p, \delta} \big) \ge \frac{c_{0}}{4}.
\end{equation}
By Proposition \ref{prop:regularity},
we can pick some $0<\epsilon < \delta^2/10^5$ and $0<a< \epsilon ^{3L}/10^4$ such that the event
\begin{equation}\label{eq:notequal}
J_{p, \delta}  \cap \{W_{n_{l_{j}},a} \cap V_{n_{l_{j}},\epsilon} \cap F_{n_{l_{j}}, \epsilon} \mbox{ occurs for a sub-sub-sequence } n_{l_{j}},\;j = 1, 2, \cdots \}       
\end{equation}
happens with probability at least $c_{0}/8$.
From now on till the end of this proof, we always work on this event. 

We will first consider the case that $t > t'$. 

\noindent \underline{{\bf Case 1}:  \  $t - t' > \frac{\delta}{2}$ }:  \\

\noindent This case is easy. Let $x_{0} = \zeta (t) = \eta^{\ast} (t') \in \mathbb{D}_{1 - p^{-1}} $. 
We write
\begin{equation}\label{eq:Adef}
A:= \widehat{\zeta}_{x_{0}} = \zeta [0, t] = \eta^{\ast} [0, t'],
\end{equation}
and then define $A_{\epsilon}$ as follows. Recall the covering $B_i$, $i=1,\ldots,M$, of $\mathbb{D}$ from Definition \ref{def:minimal}. 
We write 
\begin{equation}\label{eq:IAJA}
I_{A} : = \{ 1 \le i \le M \ | \ A \cap B_{i} \neq \emptyset \}\mbox{ and set }J_{A} := \{ 1 \le j \le M \ | \ B_{j} \cap B_{i} \neq \emptyset \text{ for some } i \in I_{A} \}     
\end{equation}
and define
\begin{equation}\label{Aprime}
A_{\epsilon} = \bigcup_{i \in I_{A} \cup J_{A}} B_{i}.
\end{equation}
We note that if a point $y$ satisfies $\text{dist} (y, A) < \frac{\epsilon}{2}$ then $y$ is contained in $A_{\epsilon}$. Also, note that $A_{\epsilon } \downarrow A$ as $\epsilon \downarrow 0$, and that $\mu (A) = t'$. Therefore, by the monotone convergence theorem, taking $\epsilon_{1} > 0 $ sufficiently small, we have 
\begin{equation*}
\mu (A_{\epsilon_{1}} ) \le t' + \frac{\delta}{10}.
\end{equation*}
Furthermore, since $\mu_{n_{l}} $ converges weakly to $\mu$, by Proposition \ref{basic1} and taking $l$ sufficiently large, it follows that 
\begin{equation}\label{saigo}
\mu_{n_{l}} (A_{\epsilon_{1}} ) \le t' + \frac{\delta}{5}.
\end{equation}
On the other hand, since $\rho (\eta_{n_{l}}, \zeta ) \to 0 $ as $l \to \infty$, taking $l$ sufficiently large, we have $\eta_{n_{l}} [0, t] \subset A_{\epsilon_{1}}$. This implies that  $\mu_{n_{l}} (A_{\epsilon_{1}} ) \ge t  > t' + \frac{\delta}{2}$. This contradicts \eqref{saigo}. Therefore, we make a contradiction for the first case.

\medskip
\noindent\underline{{\bf Case 2}:  \  $t' - t > \frac{\delta}{2}$ }:  \\

\noindent 
Recall the definition of $A$ in \eqref{eq:Adef}. Write $y_{0} = \zeta (t')$ and set $A' = \zeta [0, t'] = \widehat{\zeta}_{y_{0}}$. Since $t' > t$, we have $A \subsetneq A'$.

We know that $\displaystyle \lim_{j \to \infty } \rho (\eta_{n_{l_{j}}} , \zeta ) = 0$. Thus, taking $j$ sufficiently large, we have 
$$\rho (\eta_{n_{l_{j}}} , \zeta )  < \frac{a}{2000} \cdot \min \{ 1,  t_{\zeta} \}$$ where the constant $a$ is defined above \eqref{eq:notequal}.
Using this, we can find $0 \le s < s' \le t_{\eta_{n_{l_{j}}}}$ such that the following conditions hold:
\begin{align}\label{0726-q-1}
&\text{(I) \ $ |s - t | < \frac{a}{1000}$ \ \  and \ \  $| s' - t' | < \frac{a}{1000}$;} \notag \\
&\text{(II) \  $\big| \eta_{n_{l_{j}}} (s) - \zeta (t) \big| \le \frac{a}{1000}$ \ \ and  \ \  $\big| \eta_{n_{l_{j}}} (s') - \zeta (t') \big| \le \frac{a}{1000}$;} \notag \\
&\text{(III) \  $\text{dist} \big( \eta_{n_{l_{j}}} [0, s], \zeta [0, t] \big) \le \frac{a}{1000}$ \ \ and  \ \  $\text{dist} \big( \eta_{n_{l_{j}}} [s', t_{n_{l_{j}}}], \zeta [t', t_{\zeta}] \big) \le \frac{a}{1000}$.}
\end{align}
Since we have assumed the occurrence of $V_{n_{l_{j}},\epsilon}$, an argument similar to the proof of \eqref{as4} now implies 
\begin{equation}\label{0705-c-1}
\eta_{n_{l_{j}}} [s, s'] \not\subset B ( x_{0}, 4 \epsilon),
\end{equation}
where we recall that $x_{0} = \zeta (t) = \eta^{\ast} (t')$.

Recall that we have assumed 
 $F_{n_{l_{j}},\epsilon}$. This implies that 
\begin{equation}\label{0705-c-2}
{\rm dist} \Big( \eta_{n_{l_{j}}} [s', t_{n_{l_{j}}}] , A \Big) > \frac{\epsilon^{L}}{2}.
\end{equation}
Otherwise, 
by our choice of $a$, the distance between $\eta_{n_{l_{j}}} [0, s]$ and $\eta_{n_{l_{j}}} [s', t_{n_{l_{j}}}]$ is less than $\frac{ 2 \epsilon^{L}}{3}$ (recall that $A= \zeta[0, t]$, see \eqref{eq:Adef} for this). Therefore, we can find $s_{1} \in [0, s]$ and $s_{2} \in [s', t_{n_{l_{j}}}]$ such that the distance between  $ \eta_{n_{l_{j}}} ( s_{1} ) $ and $\eta_{n_{l_{j}}} (s_{2})$  is smaller than $\frac{ 2 \epsilon^{L}}{3}$. However, since $\eta_{n_{l_{j}}}  [s, s'] \subset \eta_{n_{l_{j}}} [s_{1}, s_{2} ]$, it follows from \eqref{0705-c-1} that $\eta_{n_{l_{j}}} [s_{1}, s_{2} ] \not\subset B \big( \eta_{n_{l_{j}}} (s_{1}), \epsilon \big)$. This means that $\eta_{n_{l_{j}}} $ has a $(\epsilon^{L}, \epsilon)$-quasi-loop at $\eta_{n_{l_{j}}} (s_{1})$, which contradicts the event $F_{n_{l_{j}},\epsilon}$. So, the inequality \eqref{0705-c-2} must hold.

\begin{figure}[h]
\begin{center}
\includegraphics[scale=0.55]{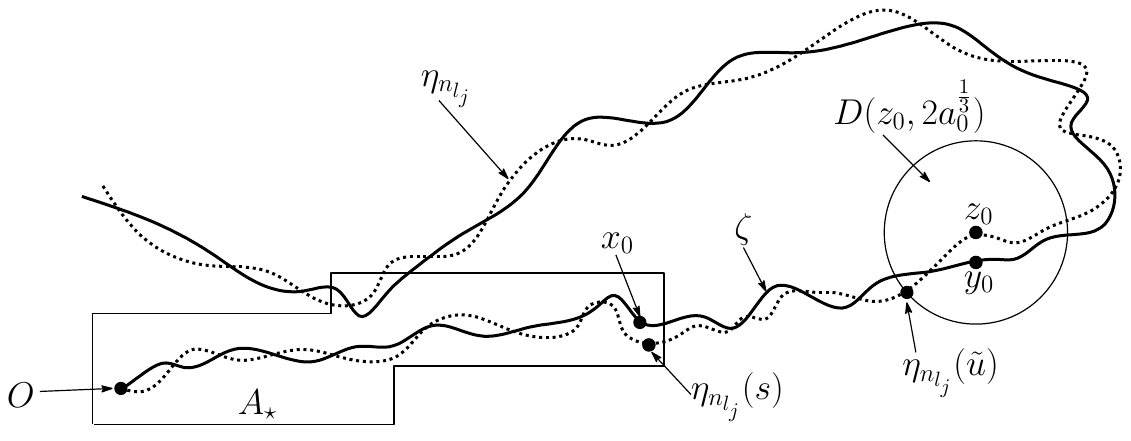}
\caption{Illustration for case 2. The solid curve stands for $\zeta$ while the dotted curve represents $\eta_{n_{l_{j}}}$. $x_{0} = \zeta (t)$, $y_{0} = \zeta (t')$ and $z_{0} = \eta_{n_{l_{j}}} (s')$. In the case 2, we make a contradiction by generating a quasi-loop at some point lying in $A_{\star}$.
}\label{0706-f} 
\end{center}
\end{figure}

With this in mind, let $z_{0} = \eta_{n_{l_{j}}} (s') $ and set
\begin{equation*}
\widetilde{u} = \max \Big\{ v \le s' \ \Big| \ \eta_{n_{l_{j}}} (v) \in \partial D \big( z_{0}, 2 a^{\frac{1}{3}} \big) \Big\}.
\end{equation*}
Combining \eqref{0705-c-2} with our choice of $a$, it holds that 
\begin{equation*}
{\rm dist} \Big( D \big( z_{0}, 2 a^{\frac{1}{3}} \big), A \Big) \ge \frac{\epsilon^{D}}{2} - 2 a^{\frac{1}{3}} \ge  \frac{\epsilon^{D}}{3}.
\end{equation*}
Furthermore, it follows from the event $W_{n_{l_{j}},a}$ that $s' - \widetilde{u}  > a$ because $| \eta_{n_{l_{j}}} ( \widetilde{u}  ) - \eta_{n_{l_{j}}} (s' ) | \ge 2 a^{\frac{1}{3}}$.
Consequently, letting 
\begin{equation*}
N_{\star} := \text{ the number of lattice points in $A_{\star}$ hit by $\eta_{n_{l_{j}}} [0, s' ]$},
\end{equation*}
where $A_{\star}:=A_{10^{-2}\cdot a}$ is defined in \eqref{Aprime},
we see that (recall \eqref{eq:fndef} for the defintion of $f_\cdot$)
\begin{equation}\label{tarinai}
N_{\star}   \leq \big( s' - a \big) f_{n_{l_j}}
\end{equation}
Since $A = \zeta[0, t] =  \eta^{\ast} [0, t']$, it follows that $\mu (A) = t'$. Because $A \subset A_{\star}$, we have $\mu (A_{\star}) \ge t'$. Furthermore, using our assumption that $\mu_{n_{l_{j}}} $ converges weakly to $\mu$, it follows from Proposition \ref{basic1} that taking $j$ sufficiently large, 
\begin{equation}\label{oki}
\mu_{n_{l_{j}}} (A_{\star}) \ge t' -  \frac{a}{100}.
\end{equation}
This implies that the number of lattice points in $A_{\star}$ hit by $\eta_{n_{l_{j}}} [0, t_{n_{l_{j}}} ]$ is bounded below by 
$$
\big( t' -  \frac{a}{100} \big) f_{n_{l_j}} \ge \big( s' -  \frac{a}{50} \big) f_{n_{l_j}}, 
$$
where we used the fact that $|s' - t' | < 10^{-3} \cdot a $ (see the condition (I) for this). 
Combining this with \eqref{tarinai},  it follows that 
\begin{equation}\label{help}
\eta_{n_{l_{j}}} [s', t_{{n_{l_{j}}}}]  \cap A_{\star} \neq \emptyset;
\end{equation}
see Figure \ref{0706-f} for this.
In particular, we have  $\text{dist} \big( \eta_{n_{l_{j}}} [s', t_{{n_{l_{j}}}}] , A \big) \le \frac{a}{3}$ since we take $A_{\star} = A_{10^{-2} \cdot a }$. 
But, this contradicts \eqref{0705-c-2} because of our choice of $a$ above \eqref{eq:notequal}. So we also make a contradiction for the second case.

This finishes the proof of Proposition \ref{end}.
\end{proof}

\section{The scaling limit of infinite loop-erased random walk}\label{sec:8}
In this section, we will consider the infinite loop-erased random walk and show that it has a scaling limit in the natural parametrization as stated in Theorem \ref{4th}.


We now briefly recap the setup. 
Let $S^{(n)}$ be the simple random walk on $2^{-n} \mathbb{Z}^{3}$ started at the origin. 
We consider the loop-erasure of $S^{(n)} [0, \infty )$ and then call the random simple path $$\gamma^{\infty}_{n} = \text{LE} \big( S^{(n)} [0, \infty ) \big)$$ the infinite loop-erased random walk (ILERW). Note that  $\lim_{k \to \infty} |\gamma^{\infty}_{n}  (k)| = \infty$ almost surely. Through linear interpolation, we can regard 
$
\gamma^{\infty}_{n} 
$
as a random element of $\cal C$. We assume the linear interpolation throughout this section. We also write (see \eqref{eq:fndef} for the definition of $f_n$)
\begin{equation}\label{irelw-scaling}
\eta^{\infty}_{n} (t) = \gamma^{\infty}_{n} \big(f_n t \big),  \ \    t \ge 0
\end{equation}
for the rescaled process which is also a random element of $\cal C$.

The following theorem then restates Theorem \ref{4th}. 

\begin{thm}\label{thm:4rep}
 There exists a random element $\eta^{\infty}$ of the space $({\cal C}, \chi)$ such that $\eta^{\infty}_{n}$ converges weakly to $\eta^{\infty}$ as $n \to \infty$ with respect to the metric $\chi$.
\end{thm}

Before diving into the proof, we need some preparation. Recall that $\mathbb{B}_{r}$ stands for the cube of side length $2r$ centered at the origin, see Section \ref{metric}. Given a curve $\lambda$,  the first time that $\lambda$ exits from $\mathbb{B}_{r}$ is denoted by $\tau^{\lambda}_{r}$, and we set $\lambda_{r} $ for the truncation of $\lambda$ up to $\tau^{\lambda}_{r}$ as defined in \eqref{trunc}. (For convenience of symbols, we sometimes also use $\lambda^{r}$ to denote the truncation in this section.)
It follows from Theorem \ref{3rd} that $\eta_{n}$ converges weakly to $\eta$ with respect to the metric $\rho$, where $\eta_{n}$ and $\eta$ are given by \eqref{etan-1st} and Theorem \ref{3rd} respectively. For each $r \in (0,1/10)$ we write $\eta_{n, r}$  for the truncation of $\eta_{n}$  up to $\tau_{r}^{\eta_{n}}$. Set $\eta^{r}$ for the truncation of $\eta$ up to $\tau^{\eta}_{r}$.  In order to prove Theorem \ref{thm:4rep}, we need the following proposition which ensures the convergence of these truncated LERWs.

\begin{prop}\label{trunc-conv}
Let $r \in (0,1/10)$. With the notation above, $\eta_{n, r}$ converges weakly to $\eta^{r}$ as $n \to \infty$ with respect to the metric $\rho$.
\end{prop}

We defer the proof of Proposition \ref{trunc-conv} to Section \ref{lem-trunc}.  We will prove Theorem \ref{thm:4rep} assuming this proposition.

\begin{figure}[h]
\begin{center}
\includegraphics[scale=0.6]{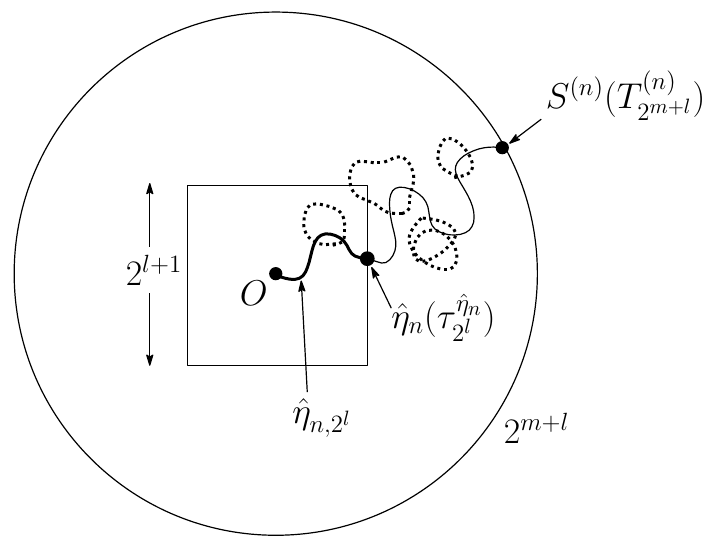}
\caption{Illustration for $\hat{\eta}_{n} =  \text{LE} \big( S^{(n)} [0, T^{(n)}_{2^{m+l}} ] \big)$ which is drawn by the solid curve. Loops erased in the loop-erasing procedure for $S^{(n)}  [0, T^{(n)}_{2^{m+l}}]$ are represented by dotted closed curves. The cube stands for $\mathbb{B}_{r} $ and the ball represents $B (2^{m+l})$. The first time that $\hat{\eta}_{n}$ exits from $\mathbb{B}_{r} $ is denoted by $\tau^{\hat{\eta}_{n}}_{2^{l}}$. The thick solid curve lying in the cube stands for $\hat{\eta}_{n, 2^{l}} $.}\label{0707-fig}
\end{center}
\end{figure}

\begin{proof}[Proof of Theorem \ref{thm:4rep}]
Fix two integers $l, m \ge 1$ with $2^{-m} <1/10$. 
We write $\eta_{n, 2^{l}}^{\infty}$ for the truncation of $\eta_{n}^{\infty}$ up to $\tau^{\eta^{\infty}_{n}}_{2^{l}}$ as defined in \eqref{trunc}.
It suffices to show that for each $l \ge 1$, $\{  \eta_{n, 2^{l}}^{\infty} \}_{n \ge 1}$ converges weakly to some random curve $\eta^{\infty, 2^{l}}$ as $n \to \infty$ with respect to the metric $\rho$ of \eqref{rho-metric}, since the collection of the limiting curves $\{ \eta^{\infty, 2^{l}} \}_{l \ge 1}$ clearly forms a consistent family.

Let $s= T^{(n)}_{2^{m+l}}$ be the first time that $S^{(n)}$ exits from $D (2^{m+l})$. We set 
\begin{equation*}
\hat{\eta}_{n} (t)  =  \text{LE} \big( S^{(n)} [0, s ] \big) \big( f_n t \big), \ \ \ 0 \le t \le f_n^{-1} \cdot \text{len} \Big( \text{LE} \big( S^{(n)} [0, s ] \big) \Big)
\end{equation*}
for the time rescaled curve obtained from the loop-erasure of $S^{(n)}$ with the truncation up to $s$. Similar to the infinite LERW case, we write   $\hat{\eta}_{n, 2^{l}} $ for the truncation of $\hat{\eta}_{n}$ up to $\tau^{\hat{\eta}_{n}}_{2^{l}}$ as defined in \eqref{trunc}, see Figure \ref{0707-fig} for the setup. It then follows from Lemma \ref{asymp-indep} that  for all $l, m \ge 1$, $n \ge 1$ and a curve $\lambda$, one has
\begin{equation}\label{mas-approx}
P ( \eta_{n, 2^{l}}^{\infty} = \lambda ) = P ( \hat{\eta}_{n, 2^{l}} = \lambda ) \big\{ 1 + O ( 2^{-m} ) \big\},
\end{equation}
where the constant coming from $ O (2^{-m} )$ is universal which can be chosen uniformly in $l, m \ge 1$, $n \ge 1$ and a curve $\lambda$. Thus, if we write $\mu_{n, 2^{l}}^{\infty}$ and $\hat{\mu}_{n, 2^{l}} $ for the probability measures induced by $\eta_{n, 2^{l}}^{\infty} $ and $\hat{\eta}_{n, 2^{l}} $ respectively, the total variation distance between $\mu_{n, 2^{l}}^{\infty}$ and $\hat{\mu}_{n, 2^{l}}$ is bounded above by $C 2^{-m}$ for some universal constant $C \in (0, \infty)$.
However, by Proposition \ref{trunc-conv}, we know that $\{ \hat{\mu}_{n, 2^{l}} \}_{n \ge 1}$ is a convergent sequence as $n \to \infty$ 
with respect to the Prokhorov metric. Since $2^{-m} < 1/10$ is arbitrary 
and the Prokhorov metric is bounded above by the total variation distance (see \cite{Gibbs} for this), this implies that $\{ \mu_{n, 2^{l}}^{\infty} \}_{n \ge 1}$ is a Cauchy sequence with respect to the Prokhorov metric. On the other hand,  by using Corollary \ref{0626} and \eqref{mas-approx} again,  we have the tightness of $\{ \eta_{n, 2^{l}}^{\infty} \}_{n \ge 1}$ with respect to the metric $\rho$. Thus, combining these facts, 
we see that $\{ \mu_{n, 2^{l}}^{\infty} \}_{n \ge 1}$ is a convergent sequence with respect to the Prokhorov metric. This finishes the proof.
\end{proof}

\begin{rem}\label{rem:lastrem}
 It would be also very interesting to ask what properties  $\eta^\infty$ 
satisfies. For instance, it is very natural to expect the law of $ r^{-1} \eta^\infty \big( r^{\beta} \bullet \big) $ is same as that of $\eta^{\infty} ( \bullet )$ for every $r > 0$.
\end{rem}

\subsection{Proof of Proposition \ref{trunc-conv}}\label{lem-trunc}
\begin{figure}[h]
\begin{center}
\includegraphics[scale=0.75]{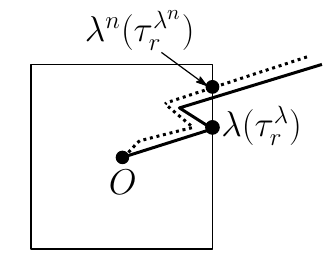}
\caption{The thick solid polygonal line stands for $\lambda$ while the dotted polygonal line represents $\lambda^{n}$. The cube stands for $\mathbb{B}_{r}$.}\label{counterex}
\end{center}
\end{figure}

In this section, we will prove Proposition \ref{trunc-conv}. We start with the following preparatory consideration. Recall that $\big( {\cal C} ( \overline{\mathbb{D}} ), \rho \big)$ is the space of continuous curves lying in $\overline{\mathbb{D}}$ defined in Section \ref{metric}. 
 Suppose that $\lambda^{n}$ ($n \ge 1$) and $\lambda $  are continuous simple curves satisfying 
 \begin{itemize}
 \item $\lambda^{n} (0) = 0$, $\lambda^{n} [0, t_{\lambda^{n}} ) \subset \mathbb{D}$ and $  \lambda^{n} (t_{\lambda^{n}} ) \in \partial \mathbb{D}$ for each $n \ge 1$ where $t_{\lambda^{n}} \in (0, \infty)$ stands for the time duration of $\lambda^{n}$, 
 
 \item $\lambda (0) = 0$, $\lambda [0, t_{\lambda} ) \subset \mathbb{D}$ and $  \lambda (t_{\lambda} ) \in \partial \mathbb{D}$ where $t_{\lambda} \in (0, \infty)$ stands for the time duration of $\lambda$,
 
 \item $\rho (\lambda^{n}, \lambda) \to 0$ as $n \to \infty$.
 
 \end{itemize}
Fix $r \in (0, 1/10 )$.   Is it true that $\rho ( \lambda^{n}_{r} , \lambda_{r} ) \to 0$ as $n \to \infty$? (Here we recall that $\lambda_{r}$ stands for the truncation of $\lambda$ up to $\tau^{\lambda}_{r}$ defined in \eqref{trunc}. The time duration of $\lambda_{r}$ is equal to $\tau_{r}^{\lambda}$.) The answer is ``no'' in general unfortunately. A counterexample for this question is illustrated in Figure \ref{counterex} where $\lambda$ has a ``reflection'' when it hits $\partial \mathbb{B}_{r}$, and 
consequently $\tau^{\lambda^{n}}_{r}$ does not converge to $\tau^{\lambda}_{r}$. However, any counterexample is limited to this kind of case. Namely, if $\lambda$ does not have a reflection at $\partial \mathbb{B}_{r}$ in the sense that for all $\epsilon  > 0$ there exists $\delta > 0$ such that $\tau^{\lambda}_{r + \delta} - \tau^{\lambda}_{r} \le \epsilon$, then one can prove that $\rho ( \lambda^{n}_{r} , \lambda_{r} ) \to 0$ as $n \to \infty$. 

 Therefore, Proposition \ref{trunc-conv} boils down to the following proposition. Recall that $\tau^{\eta}_{r}$ stands for the first time that the scaling limit $\eta$ exits from $\mathbb{B}_{r}$.
 
 \begin{prop}\label{trun-con-1}
 Let $r \in (0, 1/10)$. For all $\epsilon \in (0,1)$, there exists $\delta  > 0$ such that 
 \begin{equation}\label{tru-con-2}
 P \big( \tau^{\eta}_{r + \delta} - \tau^{\eta}_{r} \le \epsilon \big) \ge 1 - \epsilon.
 \end{equation}
 \end{prop}

Proposition \ref{trun-con-1} bears some similarity with Proposition \ref{HC-upper}. To prove it, we will first define some unwanted ``bad'' events $F_1$ through $F_5$ and control the probabilities of these events. We will prove Proposition \ref{trun-con-1} at the end of this section after handling these events.


Fix  $r \in (0, 10^{-1} )$. For $k \ge 10$, we write 
\begin{equation}\label{0716-b-1}
\delta = \delta_{k, r} = r 2^{- 2^{k^{11}}}.
\end{equation}
Let $A = A_{k, r} = \overline{\mathbb{B}}_{r + r 2^{-k}} \setminus \mathbb{B}_{r}$. Then we can find a collection of cubes $\{ B_{i} \}_{i=1}^{N_{k}}$ in the form of 
\begin{equation}\label{0716-b-2}
B_{i} = B_\infty (x_{i},r 2^{-k-1}) \ \ \text{ for }  1 \le i \le N_{k}
\end{equation} 
satisfying the following conditions:
\begin{align}
&\text{$\partial \mathbb{B}_{r} \subset \bigcup_{i=1}^{N_{k}} B_{i}$, } \  \ \text{$N_{k} \asymp  2^{2k}$, }  \ \ \text{$\| x_{i} - x_{j} \|_{\infty} \ge r 2^{-k}$ and $B_{i} \cap B_{j} \cap 2^{-n} \mathbb{Z}^{3} = \emptyset$ for all $i$, $j$ with $i \neq j$,} \notag \\
&\text{and the $L^{\infty}$ distance between $x_{i}$ and $ \mathbb{B}_{r}$ is equal to $r 2^{-k-1}$ for each $1 \le i \le N_{k}$.}      \label{0716-3} 
\end{align} 
Note that $x_{i}$ denotes the center of $B_{i}$ which lies in $\partial \mathbb{B}_{r + r 2^{-k-1}}$ for all $1 \le i \le N_{k}$. 
Thus, we have ${\rm dist} \big( B_{i}, \{ 0 \} \cup \partial \mathbb{D} \big) =  {\rm dist} \big( 0, B_{i} \big) \in [r, 2 r]$ because we assume $10 r <  1$.  With Definition \ref{0710-1} in mind, we set
\begin{equation}\label{0716-c-1}
\widehat{B}_{i} = B_\infty(x_{i},r 2^{-k-1} - \delta).
\end{equation} 
See \eqref{0716-b-1} for $\delta = \delta_{k, r}$. 

Applying  Proposition \ref{PROP2} to the case that $B =  B_{i}$ and $\widehat{B} =  \widehat{B}_{i}$,  we have 
\begin{equation}\label{0717}
P \Big( \eta_{n} \text{ hits } B_{i}  \text{ but } \eta_{n} \text{ does not hit } \widehat{B}_{i} \Big) \le  C_{r} 2^{ -k^{11}} \ \ \  \text{ for all } \  n \ge n_{k, r}
\end{equation}
where the constant $n_{k, r}$ depends only on $k$ and $r$ while the constant $C_{r}$ depends only on $r$.

Since $N_{k}  \le C 2^{2k}$ by \eqref{0716-3}, we have that for all  $n \ge n_{k, r}$
\begin{equation}\label{0717-e-1}
P \big(F_1\big) \le C_{r} \, 2^{-k^{10}}\quad \mbox{ where }\quad F_{1} := \bigcup_{i=1}^{N_{k}} \Big\{  \eta_{n} \text{ hits } B_{i}  \text{ but } \eta_{n} \text{ does not hit } \widehat{B}_{i}    \Big\}.
\end{equation}

We next deal with the presence of a quasi-loop for the LERW. Recall that ${\rm QL} (\bullet,  \ \bullet \, ; \bullet)$ is the set of quasi-loops as defined in Section \ref{qloops}. Using \eqref{eq:noquasiloop}, 
 it follows that there exist universal constants $1 \le M < \infty$ and $C < \infty$ such that for all $n$ and $k$
\begin{equation}\label{f2event}
P \big( F_2 \big) \le C 2^{- \frac{k}{M} }\quad \mbox{ where }\quad F_{2} := \Big\{ {\rm QL} \big( 100 r \cdot 2^{-k},  \, 100 r^{\frac{1}{M}} \cdot 2^{-\frac{k}{M}}; \,  \eta_{n} \big) \neq \emptyset \Big\}.
\end{equation}

Recall that $\tau^{\eta_{n}}_{r}$ is the first time that $\eta_{n}$ exits from $\mathbb{B}_{r} $. We write  
\begin{equation}\label{0717-d-3}
u = \inf \Big\{ t \ge \tau^{\eta_{n}}_{r} \ \Big| \  \big| \eta_{n} (t) - \eta_{n} \big( \tau^{\eta_{n}}_{r} \big) \big| \ge 100 r^{\frac{1}{M}} \cdot 2^{-\frac{k}{M}} \Big\},
\end{equation}
where the constant $M$ is chosen such that the inequality \eqref{f2event} holds. Using $\delta = \delta_{k,r}$ as defined in \eqref{0716-b-1}, we define the event $F_{3}$ by 
\begin{equation}\label{f3event}
F_{3} := \Big\{ \eta_{n} \big[  \tau^{\eta_{n}}_{r}, u \big] \cap \partial  \mathbb{B}_{r + \delta} = \emptyset \Big\}.
\end{equation}
Namely, the event $F_{3}$ says that the diameter of $\eta_{n}$ from its first hitting of $\partial  \mathbb{B}_{r}$ to that of $ \partial  \mathbb{B}_{r + \delta}$ is much bigger than $\delta$.
The next lemma shows that the event $F_{3}$ is unlikely to happen.

\begin{figure}[h]
\begin{center}
\includegraphics[scale=0.6]{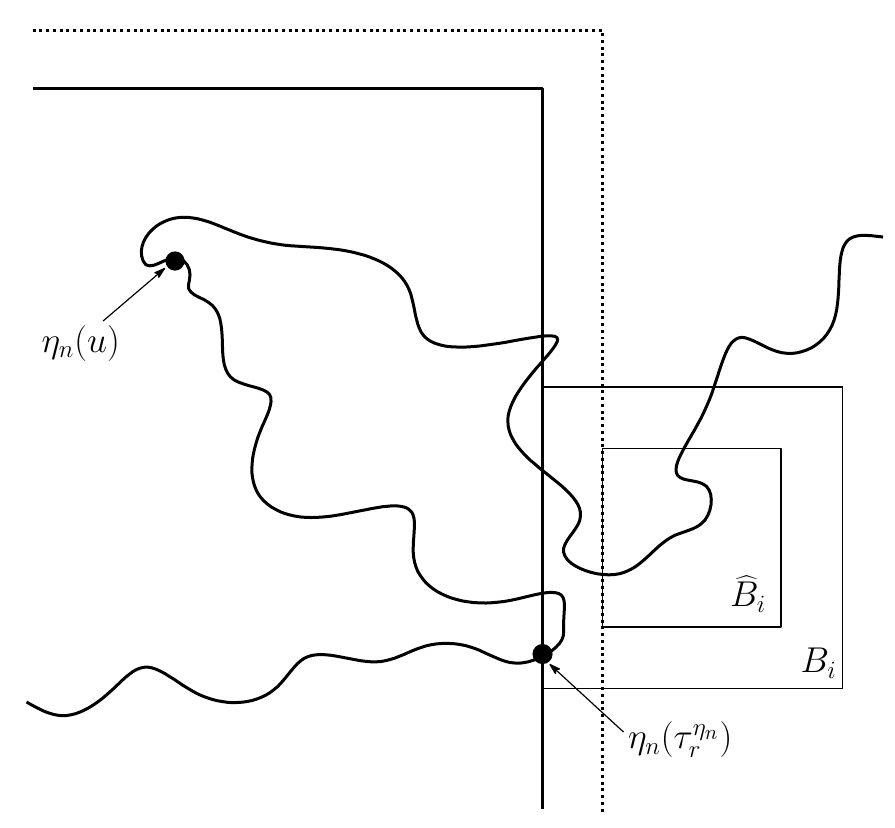}
\caption{Illustration for the event $F_{1}^{c} \cap F_{3}$. The thick solid 
polygonal line stands for $\partial \mathbb{B}_{r}$, while the thick dotted 
polygonal line represents $\partial \mathbb{B}_{r+ \delta}$. In this picture, $\eta_{n}$ is depicted by the thick solid curve. On the event $F_{1}^{c} \cap F_{3}$, $\eta_{n}$ has a quasi-loop.}\label{10-24-2}
\end{center}
\end{figure}

\begin{lem} It follows that
\begin{equation}\label{hobo}
P (F_{3} ) \le C_{r} 2^{-\frac{k}{M}}  \ \  \text{ for all } \  n \ge n_{k, r}
\end{equation}
where the constant $n_{k, r}$ depends only on $k$ and $r$, the constant $C_{r}$ depend only on $r$, and  
the constant $M$ is chosen such that the inequality \eqref{f2event} holds.
\end{lem}
\begin{proof}

Note that by \eqref{0717-e-1} and \eqref{f2event}
\begin{equation*}
P (F_{3} ) = P \big( F_{3} \cap (F_{1} \cup F_{2} )  \big) + P \big( F_{3} \cap  F_{1}^{c} \cap F_{2}^{c}   \big) \le  C_{r} 2^{-\frac{k}{M}} + P \big( F_{3} \cap  F_{1}^{c} \cap F_{2}^{c}   \big). 
\end{equation*}
We will show that $F_{3} \cap  F_{1}^{c} \cap F_{2}^{c} = \emptyset$ by contradiction.
So suppose that $F_{3} \cap  F_{1}^{c} \cap F_{2}^{c}$ occurs. Since $\{ B_{i} \}_{i=1}^{N_{k}}$ is a covering of $\partial \mathbb{B}_{r}$ by \eqref{0716-3}, there is a cube $B_{i}$ such that 
\begin{equation*}
\eta_{n} \big( \tau^{\eta_{n}}_{r}  \big) \in \partial B_{i}.
\end{equation*} 
Note that $\widehat{B}_{i} \subset \mathbb{B}_{r + \delta}^{c}$ by our construction of $\widehat{B}_{i}$, see \eqref{0716-3}  and \eqref{0716-c-1} for this. Since $F_{3}$ occurs, we see that 
\begin{equation*}
\eta_{n} \big[  \tau^{\eta_{n}}_{r}, u \big] \cap \mathbb{B}_{r + \delta}^{c} = \emptyset,
\end{equation*}
which implies that
\begin{equation*}
\eta_{n} \big[  \tau^{\eta_{n}}_{r}, u \big] \cap \widehat{B}_{i} = \emptyset.
\end{equation*}
On the other hand, since $\eta_{n}$ hits $B_{i}$ and $F_{1}^{c}$ occurs, it follows that $\eta_{n}$ must hit $\widehat{B}_{i}$. Note that $\eta_{n} [0,  \tau^{\eta_{n}}_{r}] \cap \widehat{B}_{i}  = \emptyset$ since $ \tau^{\eta_{n}}_{r}$ is the first time that $\eta_{n}$ hits $B_{i}$. Consequently, it follows that $\eta_{n} [u, t_{n} ]  \cap \widehat{B}_{i} \neq \emptyset $ where $t_{n}$
is the time duration of $\eta_{n}$. By definition of $u$ as in \eqref{0717-d-3}, the distance between $\eta_{n} ( \tau^{\eta_{n}}_{r} )$ and $\eta_{n} ( u)$ is equal to $100 r^{\frac{1}{M}} \cdot 2^{-\frac{k}{M}}$.
Thus, $\eta_{n}$ has a $(100 r \cdot 2^{-k},  100 r^{\frac{1}{M}} \cdot 2^{-\frac{k}{M}} )$-quasi-loop at $\eta_{n} ( \tau^{\eta_{n}}_{r})$ (see Figure \ref{10-24-2} for this). This contradicts the event $F_{2}^{c}$. Therefore, we conclude that  $F_{3} \cap  F_{1}^{c} \cap F_{2}^{c} = \emptyset$. Thus, the claim \eqref{hobo} follows.
\end{proof}

Recall that  the constant $M$ is chosen such that the inequality \eqref{f2event} holds. We will next consider  a collection of cubes $\{ Q_{j} \}_{j=1}^{L_{k, r}}$ of side length 
\begin{equation}\label{0717-g-1}
\epsilon_{k, r} := 100 r^{\frac{1}{M}} \cdot 2^{-\frac{k}{M}}
\end{equation}
which satisfies the following: 
\begin{align}\label{0717-h-1}
&\text{$\bullet$ $Q_{j}$ is the form of $\big\{ y \in \mathbb{R}^{3} \ \big| \ \| y - z_{j} \|_{\infty} \le \frac{\epsilon_{k, r}}{2} \big\}$ for each $1 \le j \le L_{k, r}$}, \notag \\
&\text{$\bullet$ each $z_{j}$ is an element of $\epsilon_{k} \mathbb{Z}^{3}$ with $z_{j} \neq z_{j'}$ when $j \neq j'$}, \notag \\
&\text{$\bullet$ $Q_{j} \cap Q_{j'} \cap 2^{-n} \mathbb{Z}^{3} = \emptyset$ for all $j \neq j'$,} \notag \\
&\text{$\bullet$ every $z_{j}$ lies in $\frac{4}{5} \cdot \mathbb{D}$ \ \  and  \ \ $\{ Q_{j} \}_{j=1}^{L_{k, r}}$ is a covering of $\frac{3}{4} \cdot \mathbb{D}$.}
\end{align}
Note that $L_{k, r} \asymp    \epsilon_{k, r}^{-3} $.
We write $X_{j}^{n, \epsilon_{k, r}}$ for the number of points in $Q_{j} \cap 2^{-n} \mathbb{Z}^{3}$ hit by $\eta_{n} $. Let
\begin{equation*}
F_{4} := \Big\{ X^{n, \epsilon_{k, r}}_{j} > \sqrt{\epsilon_{k, r}} \, 2^{\beta n} \text{ for some } j = 1, 2, \cdots , L_{k, r}  \Big\}.
\end{equation*}
Since the side length of each cube $Q_{j}$ is $\epsilon_{k, r}$, $X_{j}^{n, \epsilon_{k, r}}$ is typically bounded above by $C (\epsilon_{k, r} 2^{n} )^{\beta}$. Thus, the event $F_{4}$ says that there exists some cube $Q_{j}$ such that  $X_{j}^{n, \epsilon_{k, r}}$ is much bigger than this typical upper bound because $\beta > 1$, which has already been shown in \eqref{card} to be an very unlikely event, more precisely
\begin{equation}\label{0718}
P \big( F_4 \big) \le C  \,  \epsilon_{k, r}^{-3} \, e^{- \frac{c}{\sqrt{\epsilon_{k, r}}} } \le C e^{- \frac{c'}{\sqrt{\epsilon_{k, r}}} }.
\end{equation}


Finally, we define the eventual ``bad'' event $F_{5}$ by 
\begin{equation}\label{f5event}
F_{5}:= \Big\{ \tau^{\eta_{n}}_{r + \delta} - \tau^{\eta_{n}}_{r} \ge 30 \cdot \sqrt{\epsilon_{k, r}} \Big\},
\end{equation}
where we recall that $\delta= \delta_{k, r}$ and $\epsilon_{k, r}$ are as defined in \eqref{0716-b-1} and \eqref{0717-g-1}.
We have the following lemma.
\begin{lem} We have 
\begin{equation}\label{hotondo}
P (F_{5} ) \le C_{r} \epsilon_{k, r} \ \ \text{ for all } \  n \ge n_{k, r}
\end{equation}
where the constant $n_{k, r}$ depends only on $k$ and $r$ while the constant $C_{r}$ depends only on $r$. 
\end{lem}
\begin{proof}
Note that by \eqref{hobo}, \eqref{0717-g-1} and \eqref{0718} 
\begin{equation}\label{eq:h2}
P (F_{5} ) = P \big( F_{5} \cap (F_{3} \cup F_{4} )  \big) + P \big( F_{5} \cap  F_{3}^{c} \cap F_{4}^{c}   \big) \le  C_{r} \epsilon_{k, r} + P \big( F_{5} \cap  F_{3}^{c} \cap F_{4}^{c}   \big). 
\end{equation}
We will show that $F_{5} \cap  F_{3}^{c} \cap F_{4}^{c} = \emptyset$ by contradiction.
So suppose that $F_{5} \cap  F_{3}^{c} \cap F_{4}^{c}$ occurs. We can find a cube $Q_{j}$ satisfying $\eta_{n} \big( \tau^{\eta_{n}}_{r} \big) \in  Q_{j}$. We write $Q_{j_{1}}, \cdots ,  Q_{j_{q}}$ for the set of cubes with $Q_{j} \cap Q_{j_{p}} \neq \emptyset$ for $p = 1, \cdots , q$. By our assumption described in \eqref{0717-h-1}, we have $q \le 27$.

Since the event $F_{3}^{c}$ occurs, we see that 
\begin{equation*}
\eta_{n} \big[  \tau^{\eta_{n}}_{r},  \tau^{\eta_{n}}_{r + \delta} \big] \subset \bigcup_{p=1}^{q} Q_{j_{p}}.
\end{equation*}
On the other hand, since the event $F_{5}$ occurs, the number of points in $2^{-n} \mathbb{Z}^{3}$ hit by $\eta_{n} \big[  \tau^{\eta_{n}}_{r},  \tau^{\eta_{n}}_{r + \delta} \big] $ is bounded below by $30 \cdot  \sqrt{\epsilon_{k, r}} \,  2^{\beta n}$. This implies that there exists some $1 \le p \le q$ such that $X^{n, \epsilon_{k, r}}_{j_{p}} \ge \frac{10}{9} \cdot \sqrt{\epsilon_{k, r}}  \, 2^{\beta n}$. But this contradicts the event $F_{4}^{c}$. Therefore, we conclude that $F_{5} \cap  F_{3}^{c} \cap F_{4}^{c} = \emptyset$. The inequality \eqref{hotondo} then follows from \eqref{eq:h2}.
\end{proof}

We are now ready to prove Proposition \ref{trun-con-1}. 

\begin{proof}[Proof of Proposition \ref{trun-con-1}]
We will prove the following claim.

\vspace{2mm}

\noindent{\bf Claim}: \textit{Fix $r \in (0, 1/10)$ and take $\epsilon \in (0, \frac{1}{2} )$. There exist $\delta > 0$ and $N < \infty$ such that for all integer $n \ge N$ }
\begin{equation}\label{fin-1-1}
P \Big(  \tau^{\eta_{n}}_{r + \delta} -\tau^{\eta_{n}}_{r} \le \epsilon \Big) \ge 1 - \epsilon.
\end{equation}
Proposition \ref{trun-con-1} is an easy consequence of this claim and Theorem \ref{3rd}.\end{proof}
\begin{proof}[Proof of \eqref{fin-1-1}]
Fix $r \in (0, 1/10)$ and  $\epsilon \in (0, \frac{1}{2} )$. To prove \eqref{fin-1-1},  we choose $k \ge 10$ sufficiently large so that $30 \cdot \sqrt{\epsilon_{k, r}}  < \epsilon$ and $C_{r} \epsilon_{k, r} < \epsilon$ where $\epsilon_{k, r}$ is as defined in \eqref{0717-g-1}
 and the constant $C_{r}$ depending only on $r$ comes from \eqref{hotondo}. For this choice of $k$, with \eqref{0716-b-1} in mind, we let $\delta := \delta_{k, r} = r 2^{-2^{k^{11}}}$. Then it follows from \eqref{hotondo} that 
\begin{equation*}
P \Big(  \tau^{\eta_{n}}_{r + \delta} -\tau^{\eta_{n}}_{r} \le \epsilon \Big) \ge 1 - \epsilon  \ \  \text{ for all } \  n \ge n_{k, r}
\end{equation*}
where the constant $n_{k, r}$ depends only on $k$ and $r$. Namely, it depends only on $\epsilon$ and $r$ in this case.
So, by letting $N = n_{k, r}$, we finish the proof.
\end{proof}

\end{document}